\documentclass[12pt]{article}
\usepackage{amsfonts}
\usepackage{stmaryrd}
\usepackage{stmaryrd}
\usepackage{graphicx}
\usepackage{graphics}
\usepackage{epstopdf}
\usepackage{epsfig}
\usepackage{amsmath}
\usepackage{amsthm}
\usepackage{amssymb}
\usepackage{float}
\usepackage{cite}
\usepackage{multirow}
\usepackage{verbatim}
\usepackage{fancyhdr}
\usepackage{subfigure}
\usepackage{color}
\usepackage{mathtools}
\usepackage{sectsty}
\usepackage[title]{appendix}
\usepackage{threeparttable}
\usepackage{dcolumn}
\usepackage{booktabs}
\usepackage{indentfirst}
\usepackage{setspace}
\usepackage{bm}
\usepackage{enumerate}
\usepackage[labelfont=bf,labelsep=period]{caption}
\usepackage{setspace}

\newcolumntype{d}[1]{D{.}{.}{#1}}
\theoremstyle{definition}

\numberwithin{equation}{section}

\makeatletter
\renewcommand{\maketag@@@}[1]{\hbox{\m@th\normalsize\normalfont#1}}%
\makeatother

\allowdisplaybreaks
\begin{document}
\pagenumbering{arabic}
\baselineskip=1.4pc

\vspace*{0.5in}

\begin{center}

{{\bf \Large
A new type of simplified inverse Lax-Wendroff boundary treatment I: hyperbolic conservation laws}}

\end{center}

\vspace{.03in}
\vspace{.03in}

\centerline{
Shihao Liu \footnote{
	Department of Mathematics and Linné FLOW Centre,
	KTH Royal Institute of Technology, 100 44 Stockholm, Sweden. 
	E-mail: shihaoli@kth.se.},
Tingting Li \footnote{
	School of Mathematics and Statistics, 
	Henan University, Kaifeng, Henan 475004, China. 
	Research is supported by NSFC grant 11801143. 
	Email: ltt901120@henu.edu.cn},
Ziqiang Cheng \footnote{
	The School of Mathematics, 
Hefei University of Technology, Hefei, Anhui 230026, China. 
Research supported by NSFC grant 12201169.
E-mail: czq10491@hfut.edu.cn. },	
Yan Jiang \footnote{
	School of Mathematical Sciences,
University of Science and Technology of China, 
Hefei, Anhui 230026, China.  E-mail: jiangy@ustc.edu.cn. 
Research is supported in part by NSFC grant 12271499 and Cyrus Tang Foundation. },
Chi-Wang Shu\footnote{
	Division of Applied Mathematics, Brown University, 
	Providence, RI 02912, USA.
E-mail: chi-wang\_shu@brown.edu. 
Research is supported in part by NSF grant DMS-2309249.},
Mengping Zhang\footnote{
	School of Mathematical Sciences,
University of Science and Technology of China, Hefei,
Anhui 230026, China.  E-mail: mpzhang@ustc.edu.cn.
Research is supported in part by NSFC grant 12126604.}}

\vspace{.1in}

\noindent
\textbf{Abstract:}
In this paper, we design a new kind of high order inverse Lax-Wendroff (ILW) boundary treatment for solving hyperbolic conservation laws with finite difference method on a Cartesian mesh.
This new ILW method decomposes the construction of ghost point values near inflow boundary into two steps: interpolation and extrapolation. 
At first, we impose values of some artificial auxiliary points through a polynomial interpolating the interior points near 
the boundary. 
Then, we will construct a Hermite extrapolation based on those auxiliary point values and the spatial derivatives at boundary obtained via the ILW procedure. 
This polynomial will give us the approximation to the ghost point value.
By an appropriate selection of those artificial auxiliary points, high-order accuracy and stable results can be achieved. 
Moreover, theoretical analysis indicates that comparing with the original ILW method, especially for higher order accuracy, the new proposed one would require fewer terms using the relatively complicated ILW procedure and thus improve computational efficiency on the premise of maintaining accuracy and stability. 
We perform numerical experiments on several benchmarks, including one- and two-dimensional scalar equations and systems. The robustness and efficiency of the proposed scheme is numerically verified.

\bigskip


\pagenumbering{arabic}

\baselineskip=1.4pc

\section{Introduction}
In this paper, we will propose a new high order accuracy boundary treatment based on finite difference methods with fixed Cartesian mesh for hyperbolic conservation laws. For the problems on complex domain under such mesh, there are often two main difficulties. First, the computational stencil of high order finite difference scheme is often relatively wide, thus we need to evaluate the values at several ghost points near the boundary. Secondly, the physical boundary often does not happen to be on the grid points, so we need to design an algorithm to introduce the boundary conditions into our boundary scheme. If the boundary scheme is not well designed, it may bring the so-called ``cut-cell' problem, i.e., requiring the extremely small time step to ensure the stability and resulting in low computational efficiency.

A common treatment is to use body-fitted grid. That is to establish appropriate body-fitted coordinates so that the grid points lie on the physical boundary, and solve the partial differential equation in the new coordinate system. Therefore, the boundary conditions can be given directly on the grid points. The advantage of this method is that it can accurately meet the given boundary conditions. Its disadvantage is that the generation of body-fitted grid is extremely difficult. The quality of the grid directly determines the computational efficiency and accuracy. Especially for problems with moving boundary, the management of a moving grid is generally complex, which will increase the computational cost greatly. In addition, the governing equation needs to be changed during computation. The transformed PDE is often more complex than the original equation, which will also increase the computational cost.

For non body-fitted mesh methods, many scholars have also proposed some methods, such as embedded boundary method \cite{EBM1, EBM2, EBM3, EBM4, EBM5, EBM6}, immersed boundary method \cite{IBM1, IBM2, IBM3, IBM4, IBM5, IBM6}, ILW (inverse Lax-Wendroff) method \cite{TS1, TS2, TS3, TWSN} and so on. 
In this paper, we concentrate on the ILW method, which can achieve arbitrary high order accuracy and avoid the cut-cell problem effectively. 

The prototype of the earliest ILW method comes from the simulation of pedestrian flow \cite{Huang, Xiong}. The pedestrian walking direction can be determined by solving an Eikonal equation. They deal with the boundary conditions by transforming the normal derivative into the tangential derivative. 
Later, this method was extended to hyperbolic conservation law equations for the inflow boundary conditions by Tan and Shu \cite{TS1}. They repeatedly used the partial differential equation to transform the normal derivatives into time derivatives and tangential derivatives (different from the original Lax-Wendroff scheme, which transformed the time derivative into spatial derivative, which is the meaning of ``inverse"), and imposed values of ghost points near the boundary by a Taylor expansion.

After the ILW method was proposed, many scholars have done a series of work, which have greatly developed this method. To reduce the heavy algebra of the original ILW method for nonlinear systems, especially in the high-dimensional cases, the simplified ILW (SILW) method was proposed in \cite{TWSN},
in which the high order spatial derivatives were given as extrapolation instead of the relatively tedious ILW procedure.
Later, Lu et al. \cite{Lu2} proposed an ILW method to deal with ``sonic point" by evaluating the flux values at ghost points, so it can handle problems with changing wind direction. 
\cite{Ding} redefined the concept of ``conservation" for the finite difference scheme, and gave an ILW method satisfying conservation in the new sense. 
Recently, Li et al. \cite{LZSZ2023} employed the ILW boundary treatment with the fast sweeping method to capture the steady state of hyperbolic conservation laws.
In addition to hyperbolic conservation laws, the (S)ILW method was also applied to other types of equations, such as convection-diffusion equation \cite{Lu1,Li2,Li3} and Boltzmann equation \cite{FY}. 
Moreover, the boundary treatment has been successfully applied to the moving boundary problem. Tan and Shu \cite{TS2} extended the ILW method to simulate the compressible inviscid fluid containing moving circumferential wall. 
Along the same lines, by redefining the material derivative on the boundary,  \cite{Cheng} extended the method to deal with the arbitrary motion of the boundary, and used it to simulate the interaction between shock wave and rigid body with complex geometry. ILW method for convection-diffusion equations on moving domain was proposed by Liu et al. \cite{Liu1}, in which a unified algorithm was designed for pure convection, convection-dominated, convection-diffusion, diffusion-dominated and pure diffusion cases. 
The three-dimensional cases were concerned in \cite{Liu2} with simulation of the interaction between inviscid/viscous fluid and three-dimensional rigid body. 
Besides, references \cite{Li1,Li2,Li3,VS} have analyzed the linear stability of ILW and SILW methods, which provide guidelines for us to design stable ILW boundary treatments.

In this paper, we will design a new type of SILW method for conservation laws to further reduce the computational complexity.
The new method decomposes the construction of the ghost point values into two steps: interpolation and extrapolation.
At first, we approximate values on some artificial auxiliary points through the interpolating polynomial based on interior points near boundary. Then, we will construct a Hermite extrapolation based on those auxiliary points values and the spatial derivatives at the boundary obtained through the ILW procedure. This extrapolating polynomial will give us the approximation of those ghost point values.
Linear stability analysis will be performed, requiring maintaining stability with the same Courant-Friedrichs-Lewy (CFL) number as the periodic boundary case for any boundary locations.
Moreover, under the premise of stability and same order of accuracy,
we aim to reduce the necessary terms obtained via the ILW procedure comparing with the original SILW method by 
an appropriate selection of those artificial auxiliary points.
Thus, the new method can improve the computational efficiency, especially for high-dimensional systems.

The organization of this paper is as follows. In Section 2, we will give the description of the new ILW method for one-dimensional scalar conservation laws, and use the eigenvalue analysis method to perform the linear stability analysis. 
In Section 3, we will extend this algorithm to systems and high-dimensional cases. 
High order accuracy and robustness of our algorithm will be showed through numerical tests in Section 4. 
Conclusion remarks will be given in Section 5. 

\section{The one-dimensional scalar conservation laws}
Consider the one-dimensional scalar hyperbolic conservation law in the following form:
\begin{equation}
\left\{\begin{split}
&u_t+f(u)_x=0, \quad x \in (-1,1), \, t>0,\\
&u(-1,t)=g_l(t), \quad t>0, \\
&u(x,0)=u_0(t), \quad x\in[-1,1].
\end{split}
\right.
\label{eq:1dcl}
\end{equation}
We assume that $f'(u(-1,t))>0$, such that the left boundary $x=-1$ is an inflow boundary, where a boundary condition needs to be given. We also assume that $f'(u(1,t))>0$. Hence the right boundary $x=1$ is an outflow boundary, where no boundary condition is required.

Suppose the domain is divided by the uniform mesh:
\begin{equation}\label{eq:mesh_1D}
	-1+C_a\Delta x=x_0<\cdots<x_{N}=1-C_b\Delta x
\end{equation}
with mesh size $\Delta x=2/(C_a+C_b+N)$ and $C_a,C_b\in[0,1)$. Note that we have allowed the physical boundary $x =\pm 1$ not coinciding with grid points 
when $C_a, C_b\neq0$. Here, a uniform algorithm for both body-fitted and non body-fitted grids is under consideration. 

We use the framework of method of lines (MOL) to construct a semi-discrete scheme on the interior point $x_j, \, j=0,1,2,\cdots,N$:
\begin{equation} \label{semid}
\frac{d}{dt}u_j=L_h(u)_j,
\end{equation}
where,
$$L_h=-\frac{1}{\Delta x}(\hat{f}_{j+1/2}-\hat{f}_{j-1/2}) \approx -f(u)_x|_{x_{j}}$$
is the spatial discretization operator. Here, $u_j(t)$ is the numerical approximation to the exact solution $u(x_j,t)$, and $\hat{f}_{j+1/2}$ is the numerical flux. In this paper, we will use an upwind-biased conservative finite difference scheme to construct $\hat{f}_{j+1/2}$, such as the WENO scheme \cite{Jiang}.

After the spatial discretization, the semi-discrete scheme \eqref{semid} is a system of ordinary differential equations. For time discretization, we use the total variation diminishing (TVD) Runge-Kutta (RK) scheme \cite{SO}. From time level $t^n$ to $t^{n+1}$, the third order TVD RK scheme is given as 
\begin{equation}\label{eq:RK}
\begin{split}
&u^{(1)}_{j} = u^n_{j}  + \Delta tL_h(u^n)_{j} , \\
&u^{(2)}_{j}  = \frac{3}{4}u^n_{j} +\frac{1}{4}u^{(1)}_{j} +\frac{1}{4}\Delta tL_h(u^{(1)})_{j} , \\
&u^{n+1}_{j}  = \frac{1}{3}u^n_{j} +\frac{2}{3}u^{(2)}_{j} +\frac{2}{3}\Delta tL_h(u^{(2)})_{j} .
\end{split}
\end{equation}
In particular, \cite{CGAD} pointed out that the boundary conditions in the intermediate stages of the above RK scheme should be modified as follows to avoid order reduction:
\begin{equation}
\begin{split}
&u^{n} \sim g_l(t_n), \\
&u^{(1)} \sim g_l(t_n)+\Delta t g_l'(t_n), \\
&u^{(2)} \sim g_l(t_n)+\frac{1}{2}\Delta t g_l'(t_n)+\frac{1}{4}\Delta t^2 g_l''(t_n).
\end{split}
\end{equation}


Note that for a high order finite difference scheme, a wide computational stencil is generally required. Hence, it is inevitable that some points in the computational stencil are not in our computational domain,
$$x_{-p} = x_0 - p\Delta x, \quad x_{N+p}=x_{N}+p\Delta x, \quad p=1, 2, \cdots.$$ 
Therefore, we can regard the boundary treatment problem as construction of the ghost point values. In the following, we will first review the original (S)ILW method proposed by Tan et al. \cite{TS1,TWSN}.
And then, a new SILW method will be proposed to improve the computational efficiency on the premise of maintaining accuracy and stability. Linear stability analysis will be given to demonstrate the advantage of the new proposed method.

\subsection{Review of the original (S)ILW method}

The main idea of the original inverse Lax-Wendroff method for hyperbolic conservation laws \cite{TS1} is to convert the spatial derivatives into the time derivatives through the PDE and boundary conditions at the inflow boundary. 
At the outflow boundary, the spatial derivatives of each order are approximated by extrapolation. 
After that, the values of the ghost points outside the computational domain are obtained by Taylor
expansion at the boundary. More specifically, 
the ghost points near outflow boundaries, such as the right boundary $x=1$ in our example problem
(\ref{eq:1dcl}), can be obtained by extrapolation directly. We can choose the traditional Lagrange extrapolation with appropriate accuracy when the solution is smooth near the boundary, or least square extrapolation / WENO type extrapolation \cite{TS1,TWSN, Lu2} when the solution contains discontinuities near the boundary.

To be specific, for the inflow boundary, such as the left boundary $x=-1$ in our example problem (\ref{eq:1dcl}), to ensure our boundary treatment has $d$-th order accuracy, the value of the ghost points near $x=-1$ is obtained by Taylor expansion:
\begin{equation}\label{taylor}
u_j=\sum_{k=0}^{d-1}\frac{(x_j+1)^k}{k!}u^{*(k)},\quad j=-1,-2,\cdots
\end{equation}
where, $u^{*(k)}$ is the approximation of $\partial^{(k)}_x u|_{x=-1}$ with at least $(d-k)$-th order accuracy.
Using PDE and boundary condition repeatedly, we have that 
\begin{equation} \label{eq:ilw_scalar}
\begin{split}
&\partial^{(0)}_xu|_{x=-1}=g_l(t),\\
&\partial^{(1)}_xu|_{x=-1}=\frac{g_l'(t)}{-f'(g_l(t))},\\
& \partial^{(2)}_xu|_{x=-1}=\frac{f'(g_l(t))g_l''(t)-2f''(g_l(t))g_l'(t)^2}{f'(g_l(t))^3},\\
&...
\end{split}
\end{equation}
Thus, we can set
$$u^{*(k)}=\partial^{(k)}_xu|_{x=-1}.$$

To avoid the very heavy algebra of the above original ILW method when calculating the high order space derivatives, 
the simplified ILW (SILW) method was proposed in \cite{TWSN}. Specifically, $u^{*(0)},u^{*(1)} $ are constructed by the original ILW procedure, i.e., converting the spatial derivatives into the time derivatives through the PDE and boundary condition. The higher order spatial derivatives $u^{*(k)}$, $2 \leq k \leq d-1 $,  are extrapolated from the interior points directly.
This method can greatly improve the computational efficiency, especially for high-dimensional systems. 
However, \cite{Li1} analyzed the linear stability of the SILW method through the eigenvalue method,
showing that the SILW method is stable with the same CFL number as the periodic boundary case for any $C_a\in[0,1)$ only when $d=3$, 
but need less time step or even unstable for $d>3$. In order to guarantee the stability, more spatial derivatives need to be constructed via the ILW procedure at the boundary.

Suppose that for a $d$-th order scheme, $u^{*(k)}$ is obtained through the ILW procedure if $k \leq k_d-1$, or by extrapolation if $k_d \leq k \leq d-1 $. For different $d$, \cite{Li1} used the eigenvalue analysis to find out the minimum $k_d$, denoted by $(k_d)_{min}$, to make sure the scheme is stable for all $C_a\in[0,1)$. 
Values of $(k_d)_{min}$ for a variety of $d$ are shown in Table \ref{TAB:SILW_STABILITY}. It can be seen that the $(k_d)_{min}$ is still large for high order schemes.  As a result, the computational efficiency is still low for higher
order methods, especially for high dimensional systems.

\begin{table}
	\caption{The table of $(k_d)_{min}$ for original SILW method.}
	\centering
		\begin{tabular}{c|ccccccc|}
			\hline
			 $d$ & 3 & 5  & 7 & 9 & 11 & 13 \\
			 \hline
			 $(k_d)_{min}$& 2 & 3  & 4 & 6 & 8 & 10\\
			\hline
		\end{tabular}
\label{TAB:SILW_STABILITY}
\end{table}

In summary, the above SILW method can be divided into the following two steps, i.e.,  ``interpolation" and ``extrapolation":\\

\noindent
\textbf{Algorithm 1. Original SILW method for 1D scalar cases}
\begin{itemize}
\item[Step 1.] Construct a polynomial $p(x)$ 
 of degree at most $d-1$ interpolating interior points $\{x_0, \cdots, x_{d-1}\}$,
and obtain the approximation of spatial derivatives of each order on the boundary with $k_d\geq (k_d)_{min}$,
$$u^{*(k)} \approx \partial_x^{k} p|_{x=-1}, \quad k=k_d,\cdots, d-1.$$

\item[Step 2.] Construct the extrapolation polynomial $q(x)$ of degree at most $d-1$ satisfying 
$$q^{(k)}(-1) = u^{*(k)}, \quad  k=0,\cdots, d-1,$$ 
where, $u^{*(k)}$ for $k<k_d$ are obtained by the ILW procedure (\ref{eq:ilw_scalar}), and the others are obtained by Step 1. Actually, in this case, $q(x)$ is the Taylor expansion polynomial.
%
Then, we can get the ghost point values
$$u_{j}=q(x_{j}), \quad j=-1,-2,\ldots .$$
\end{itemize}

In particular, if the solution of the equation has discontinuities near the boundary, the WENO extrapolation technique can be utilized to prevent spurious oscillations. More details about WENO extrapolation can be found in \cite{TS1,TWSN,Lu2}. 

On the other hand, it is observed that the information utilized to construct the extrapolation polynomial $q(x)$ mainly consists of two parts. One includes the first $k_d-1$ spatial derivatives at the boundary obtained by the ILW procedure, and the other part is the polynomial $p(x)$ interpolating interior points $\{x_0, \cdots, x_{d-1}\}$. 
Thus, we employ $(k_d+d)$ information to construct a $d$-th order approximation, which is an underdetermined system 
of equations, hence the way of construction is not unique. 
Therefore, we would like to explore this freedom and find another way to obtain the approximation from 
the same set of information, and expect the new algorithm to be more efficient with $k_d$ as small as possible under the premise of stability and accuracy.

\subsection{A new SILW method}
In the following, we will describe our new SILW method for the scalar conservation law (\ref{eq:1dcl}), asking for smaller $(k_d)_{min}$ for the same $d$.
The key difference between the new method and the original one is that the extrapolation polynomial $q(x)$ will employ the values on some artificial auxiliary points in the computational domain,  
\begin{align*}
u(-1+k\alpha\Delta x,t)\approx  p(-1+k\alpha\Delta x),\quad k = 1, 2 \ldots 
\end{align*}
instead of using the high order derivatives of the interpolation polynomial $p(x)$ at the boundary. Here, $\alpha\geq0$ is a parameter to be determined such that the $(k_d)_{min}$ could be as small as possible.

Specifically, we summarize the procedure of our new SILW method with $d$-th order accuracy as follows:\\

\noindent
\textbf{Algorithm 2. New SILW method for 1D scalar cases}
\begin{itemize}
\item[Step 1.]
Obtain the interpolating polynomial $p(x)$ of degree at most $d-1$ based on the points $\{x_{0},\cdots,x_{d-1}\}$. Let
$$u_{k*} = p(-1+k\alpha\Delta x),\quad 1 \leq k \leq d-k_d.$$
\item[Step 2.]Construct the extrapolation polynomial $q(x)$ of degree at most $d-1$ to satisfy the following conditions:
\begin{align*}
\begin{aligned}
&q^{(k)}(-1)=\partial^{(k)}_xu|_{x=-1},\quad 0\leq k \leq k_d-1,\\
&q(-1+k\alpha\Delta x)=u_{k*},\quad 1 \leq k \leq d-k_d,
\end{aligned}
\end{align*}
where, $\partial^{(k)}_xu|_{x=-1}$ is obtained by the ILW procedure.
Let the ghost point values be the values of the extrapolation polynomial $q(x)$ at the corresponding points:
$$u_{j}=q(x_{j}) \quad j=-1,-2,\ldots.$$
\end{itemize}


Notice that, for the problems with changing wind direction (i.e. $f'(g_l(t))=0$ in the scalar case), the above inverse Lax-Wendroff procedure may contain zero denominator. \cite{TS1,Lu2} gave different methods to deal with such situations, which can be utilized in our method as well.

In the next subsection, we will show that through adjusting the parameter $\alpha$, our new SILW method could improve computational efficiency and stability comparing with the original SILW method.

\subsection{Linear stability analysis}

Here, we will give the stability analysis of the fully discrete schemes using the eigenvalue spectrum visualization.
We consider the case of $d=2k-1\, (k=2,3,4,5,6,7)$ and assume that $f'(u)>0 $. The conservative linear upwind scheme is used for spatial discretization. That is, $L_h$ in the scheme (\ref{semid}) is in the following form:
\begin{itemize}
\item[$d=3$:]
$$ L_h(u)_j=-\frac{1}{\Delta x} \left( \frac{1}{6}f_{j-2}-f_{j-1}+\frac{1}{2}f_{j}+\frac{1}{3}f_{j+1} \right),$$

\item[$d=5$:]
$$ L_h(u)_j=-\frac{1}{\Delta x} \left( -\frac{1}{30}f_{j-3}+\frac{1}{4}f_{j-2}-f_{j-1}+\frac{1}{3}f_{j}+\frac{1}{2}f_{j+1}-\frac{1}{20}f_{j+2} \right),$$

\item[$d=7$:]
$$ L_h(u)_j=-\frac{1}{\Delta x} \left( \frac{1}{140}f_{j-4}-\frac{7}{105}f_{j-3}+\frac{3}{10}f_{j-2}-f_{j-1}+\frac{1}{4}f_{j}+\frac{3}{5}f_{j+1}-\frac{1}{10}f_{j+2}+\frac{1}{105}f_{j+3} \right),$$

\item[$d=9$:]
\begin{align*}
\begin{aligned}
L_h(u)_j=-\frac{1}{\Delta x}& \left( -\frac{1}{630}f_{j-5}+\frac{1}{56}f_{j-4}-\frac{2}{21}f_{j-3}+\frac{1}{3}f_{j-2}-f_{j-1}+\frac{1}{5}f_{j} \right.\\
&\left. +\frac{2}{3}f_{j+1}-\frac{1}{7}f_{j+2}+\frac{1}{42}f_{j+3}-\frac{1}{504}f_{j+4} \right),
\end{aligned}
\end{align*}

\item[$d=11$:]
\begin{align*}
\begin{aligned}
L_h(u)_j=-\frac{1}{\Delta x}& \left( \frac{1}{2772}f_{j-6}-\frac{1}{210}f_{j-5}+\frac{5}{168}f_{j-4}-\frac{5}{42}f_{j-3}+\frac{5}{14}f_{j-2}-f_{j-1}+\frac{1}{6}f_{j} \right.\\
&\left. +\frac{5}{7}f_{j+1}-\frac{5}{28}f_{j+2}+\frac{5}{126}f_{j+3}-\frac{1}{168}f_{j+4}+\frac{1}{2310}f_{j+5} \right.),
\end{aligned}
\end{align*}

\item[$d=13$:]
\begin{align*}
\begin{aligned}
& L_h(u)_j=-\frac{1}{\Delta x} \left( -\frac{1}{12012}f_{j-7}+\frac{1}{792}f_{j-6}-\frac{1}{110}f_{j-5}+\frac{1}{24}f_{j-4}-\frac{5}{36}f_{j-3}+\frac{3}{8}f_{j-2} \right.\\
& \left. -f_{j-1}+\frac{1}{7}f_{j}+\frac{3}{4}f_{j+1}-\frac{5}{24}f_{j+2}+\frac{1}{18}f_{j+3}-\frac{1}{88}f_{j+4}+\frac{1}{660}f_{j+5}-\frac{1}{10296}f_{j+6} \right).
\end{aligned}
\end{align*}
\end{itemize}

In particular, for the linear case $f(u)=u$, 
the semi-discrete scheme (\ref{semid}) 
coupled with the numerical boundary treatments
can be written in the matrix-vector form, 
\begin{equation*}
\frac{d\bm{U}}{dt}= \frac{1}{\Delta x} \bm{Q}\bm{U}  +\bm{B},
\end{equation*}
where, $\bm{U}=(u_0,u_2,...,u_N)^T$, $\bm{Q}$ is the coefficient matrix of the spatial discretization, and $\bm{B}$ is a vector corresponding to the inflow boundary condition. 
In our analysis we suppose the linear equation is imposed with homogeneous Dirichlet condition, i.e. $g_l(t)=0$, which leads $\bm{B}=\bm{0}$. 

References \cite{VS, Li1} pointed out that we only need to care about the fixed eigenvalues of the matrix $\bm{Q}$ with different grid number 
for stability analysis. If we use the third-order TVD RK time discretization \eqref{eq:RK}, the stability region can be expressed as
\begin{equation}
|z(\mu)|\leq 1, \quad z(\mu)=1+\mu+\frac{\mu^2}{2}+\frac{\mu^3}{6},
\end{equation}
where, $\mu=s\frac{\Delta t}{\Delta x}$, and $s$ is the fixed eigenvalue of $\bm{Q}$ 
which would not change as $N$ varies. Notice that $z$ may not exist or there may be more than one. If there are more than one $z$, we consider the largest $|z(\mu)|$. 
In order to avoid the cut-cell problem, we expect the time step would not be effected by boundary treatment.
Hence, we discuss stability on the premise of maximum CFL number $\frac{\Delta t}{\Delta x} = (\lambda_{cfl})_{max}$, where, $(\lambda_{cfl})_{max}$ is the maximum CFL number for the corresponding Cauchy problem or periodic boundary case, and their specific values, which can be calculated by Fourier analysis, are shown in the Table \ref{tab:cfl}.

\begin{table}[htb!]
	\caption{The maximum CFL number for Cauchy problem.}
	\centering
		\begin{tabular}{c|ccccccc|}
			\hline
			 $d$& 3 & 5  & 7 & 9 & 11 & 13 \\\hline
			 $(\lambda_{cfl})_{max}$ & 1.62 & 1.43  & 1.24 & 1.12 & 1.04 & 0.99\\
			\hline
		\end{tabular}
	\label{tab:cfl}
\end{table}

\begin{figure}[htb!]
\centering
\subfigure[$\alpha=0.60$]{
		\includegraphics[width=0.47\textwidth]{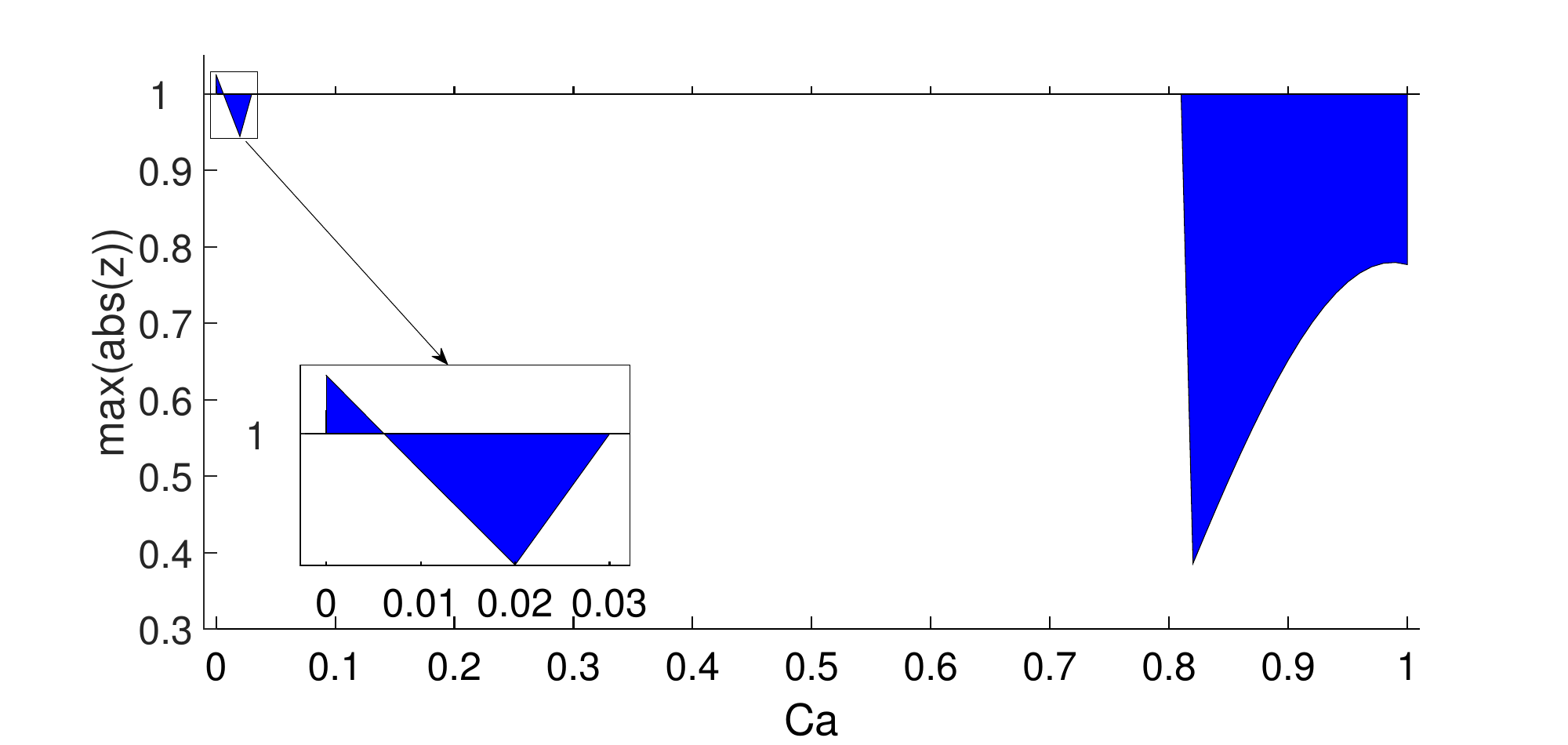}}
\subfigure[$\alpha=0.61$]{
		\includegraphics[width=0.47\textwidth]{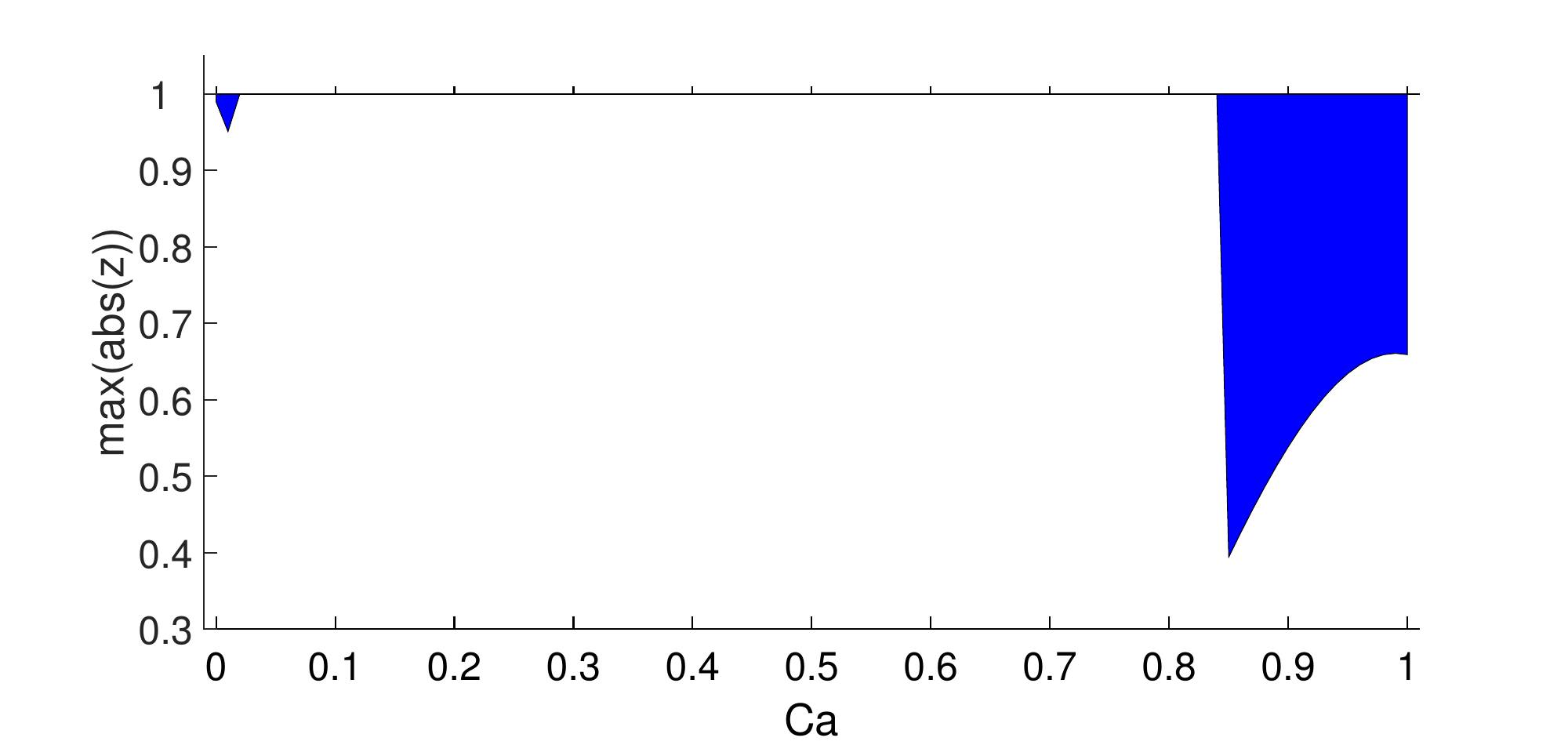}}
\caption{The third order scheme with the new SILW procedure and $k_d=2$. The horizontal axis represents $C_a$ and the vertical axis represents the largest $|z(\mu)|$.
}
\label{fig:stability_add}
\end{figure}

We select several groups of different $\alpha\in[0,10]$ and $k_d$ for linear stability analysis.
By using the software Matlab, we show the 
$\max|z(\mu)|$ for $C_a\in[0,1)$ with different $\alpha$. 
In particular, to obtain $z$, we compute eigenvalues of the matrix $\bm{Q}$ with different $N$. The candidate eigenvalues $s$ is the common value 
when $N$ changes over this range.
In Figure \ref{fig:stability_add}, we plot the results of the third order scheme $d=3$ coupling the new SILW method 
with $k_d=2$ as an example. When $\alpha=0.60$, $\max|z(\mu)|>1$ as $C_a$ approaches 0. 
However, $\max|z(\mu)|\leq1$ for all $C_a\in[0,1)$ if $\alpha=0.61$. This indicates that we should take 
$\alpha \geq 0.61$ to guarantee the third order scheme is stable with $k_d=2$. More cases are placed in Appendix 
\ref{sec:append}. 
Finally, the minimum $k_d$ and the corresponding appropriate range of $\alpha$ for 
different $d$ are shown in Table \ref{TAB:NEWSILW_STABILITY}.
As can be seen from the Table \ref{TAB:NEWSILW_STABILITY}, compared with the original SILW method, we can construct stable boundary treatments with smaller $k_d$ by adjusting $\alpha$. 
This advantage will be more prominent as $d$ increases.

\begin{table}[htb!]
	\caption{Linear stability analysis results of the new SILW method.}
	\centering
		\begin{tabular}{c|ccccccc|}
			\hline
			 $d$&  3 & 5  & 7 & 9 & 11 & 13 \\\hline
			 $(k_d)_{min}$ & 2 & 2  & 2 & 3 & 3 & 4\\
			 $\alpha$ & [0.61,10] &[0.92,5.11]&[1.34,1.99]&[1.29,2.43]&[1.42,1.70]&[1.49,2.08]\\
			\hline
		\end{tabular}
	\label{TAB:NEWSILW_STABILITY}
\end{table}


Next, we want to verify the results of the above stability analysis numerically. Consider the following problem:
\begin{equation}\label{ex1}
\left\{\begin{aligned}u_t+u_x&=0,\quad -1< x< 1,\, t>0,\\
u(x,0)&=0.25+0.5\sin(\pi x),\quad -1\le x\le 1,\\
u(-1,t)&=0.25+0.5\sin(\pi t),\quad t>0.\end{aligned}\right.
\end{equation}
The exact solution is
\begin{equation*}
u(x,t)=0.25+0.5\sin(\pi(x-t)).
\end{equation*}

\begin{figure}[htb!]
	\centering
	\subfigure[$\alpha$ between 0.61 and 10. ]{
		\includegraphics[width=0.45\textwidth]{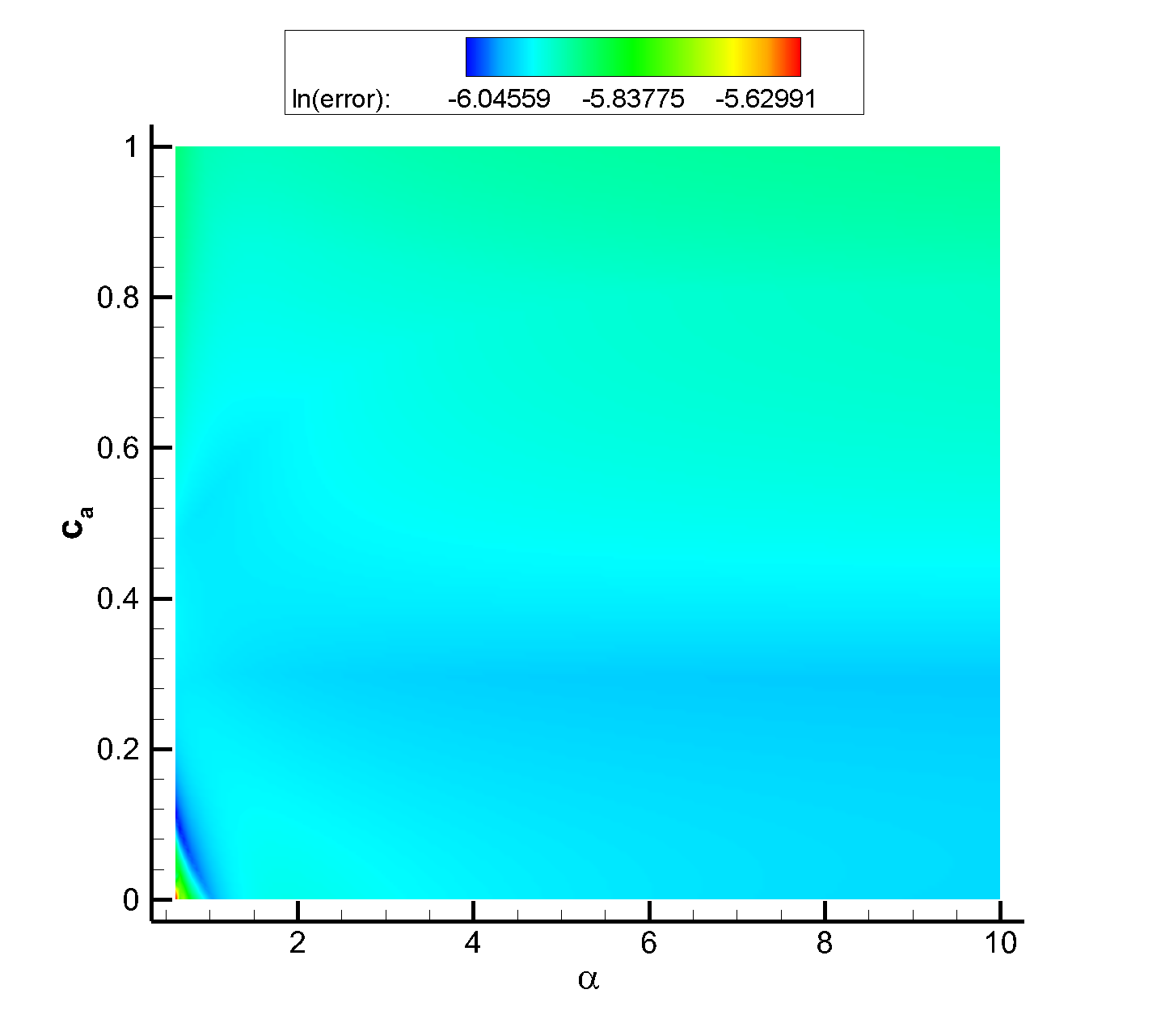}}
	\subfigure[$\alpha$ near the left boundary $\alpha=0.61$. ]{
		\includegraphics[width=0.45\textwidth]{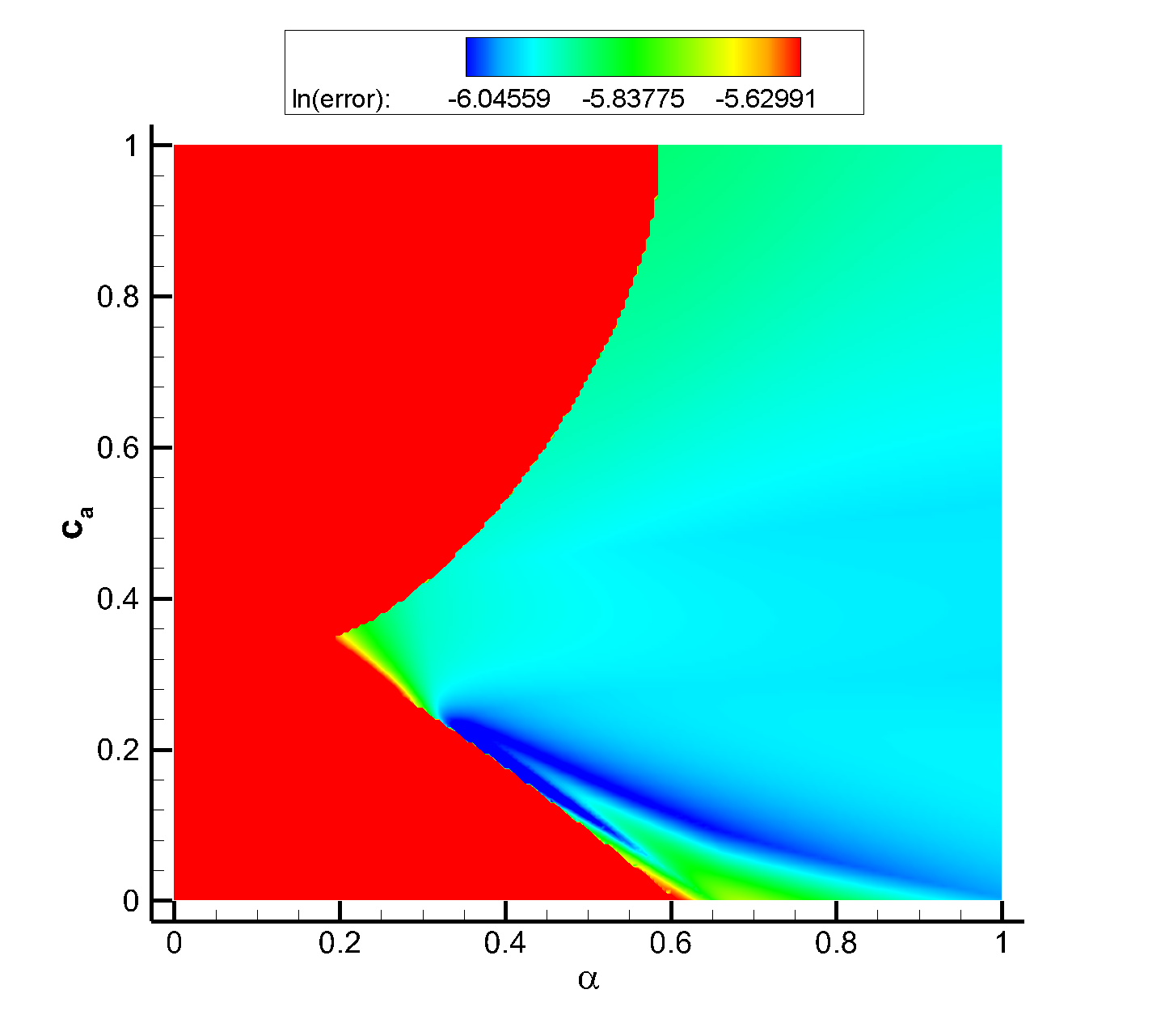}}
	\caption{Error contour of the linear problem with the third-order boundary treatment. }
	\label{fig:logerror1}
\end{figure}

\begin{figure}[thb!]
	\centering
	\subfigure[$\alpha$ between 0.92 and 5.11 . ]{ 
		\includegraphics[width=0.45\textwidth]{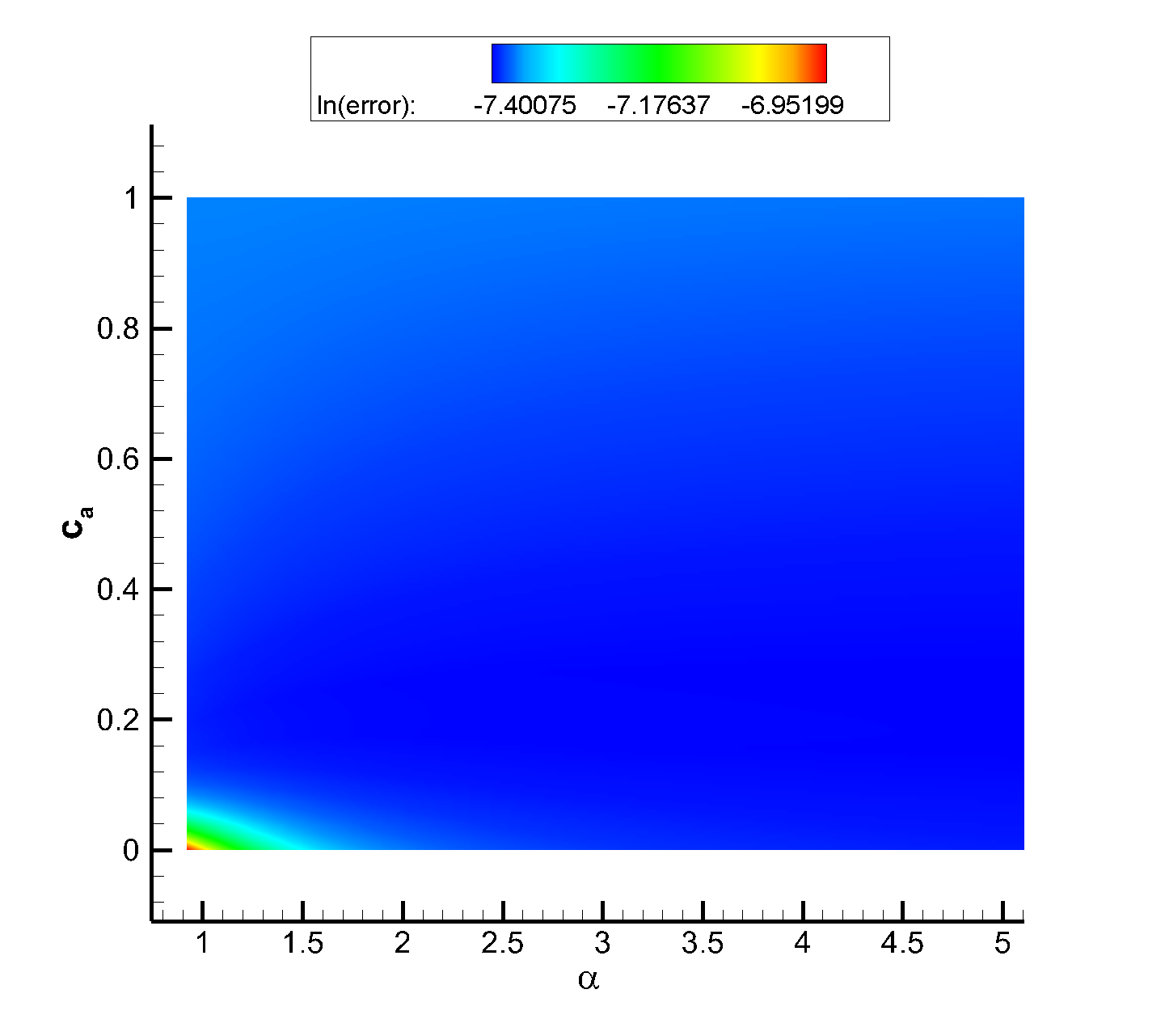}}
	\subfigure[$\alpha$ near the left boundary $\alpha=0.92$. ]{
		\includegraphics[width=0.45\textwidth]{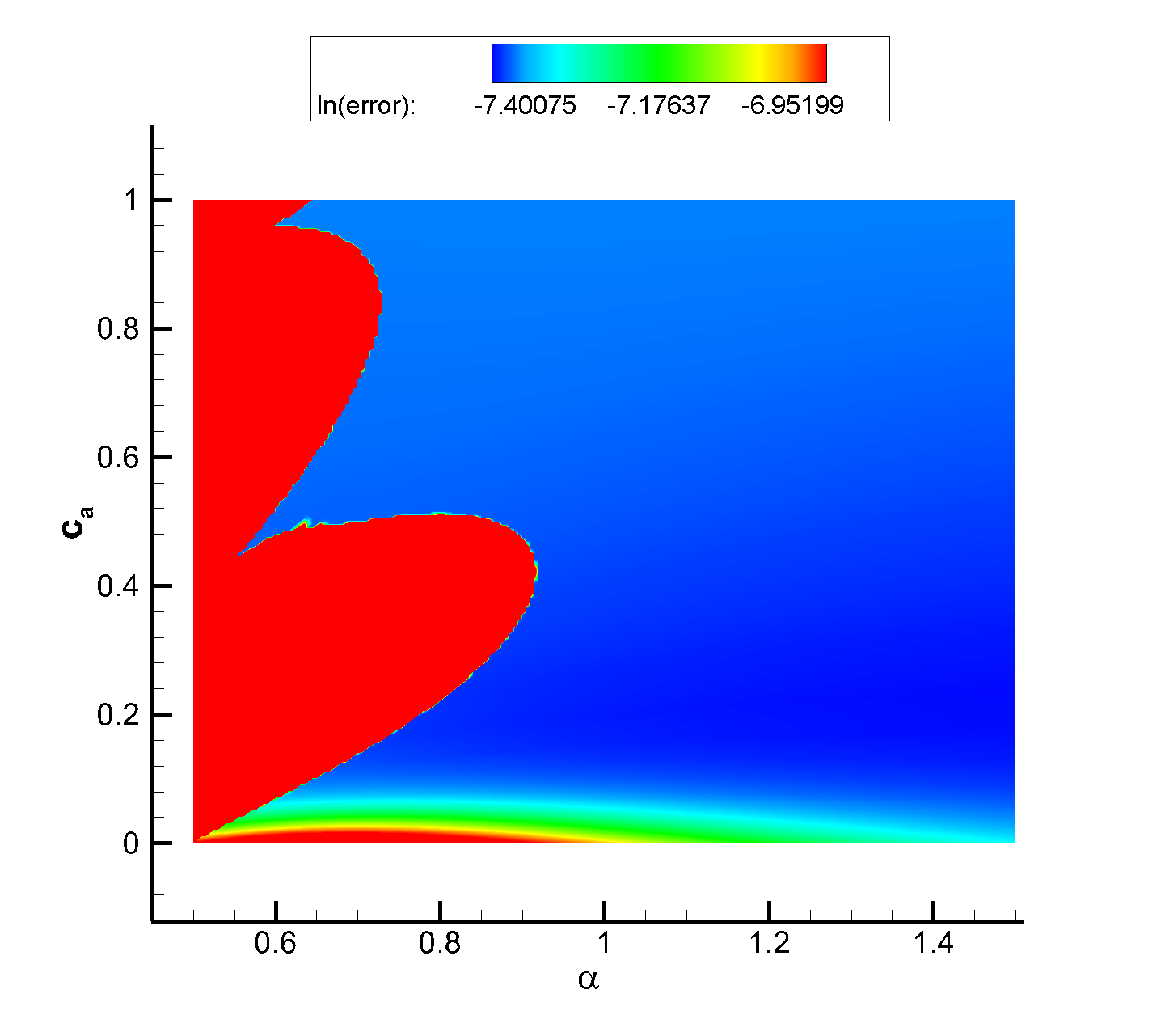}}
	\subfigure[$\alpha$ near the right boundary $\alpha=5.11$. ]{
		\includegraphics[width=0.45\textwidth]{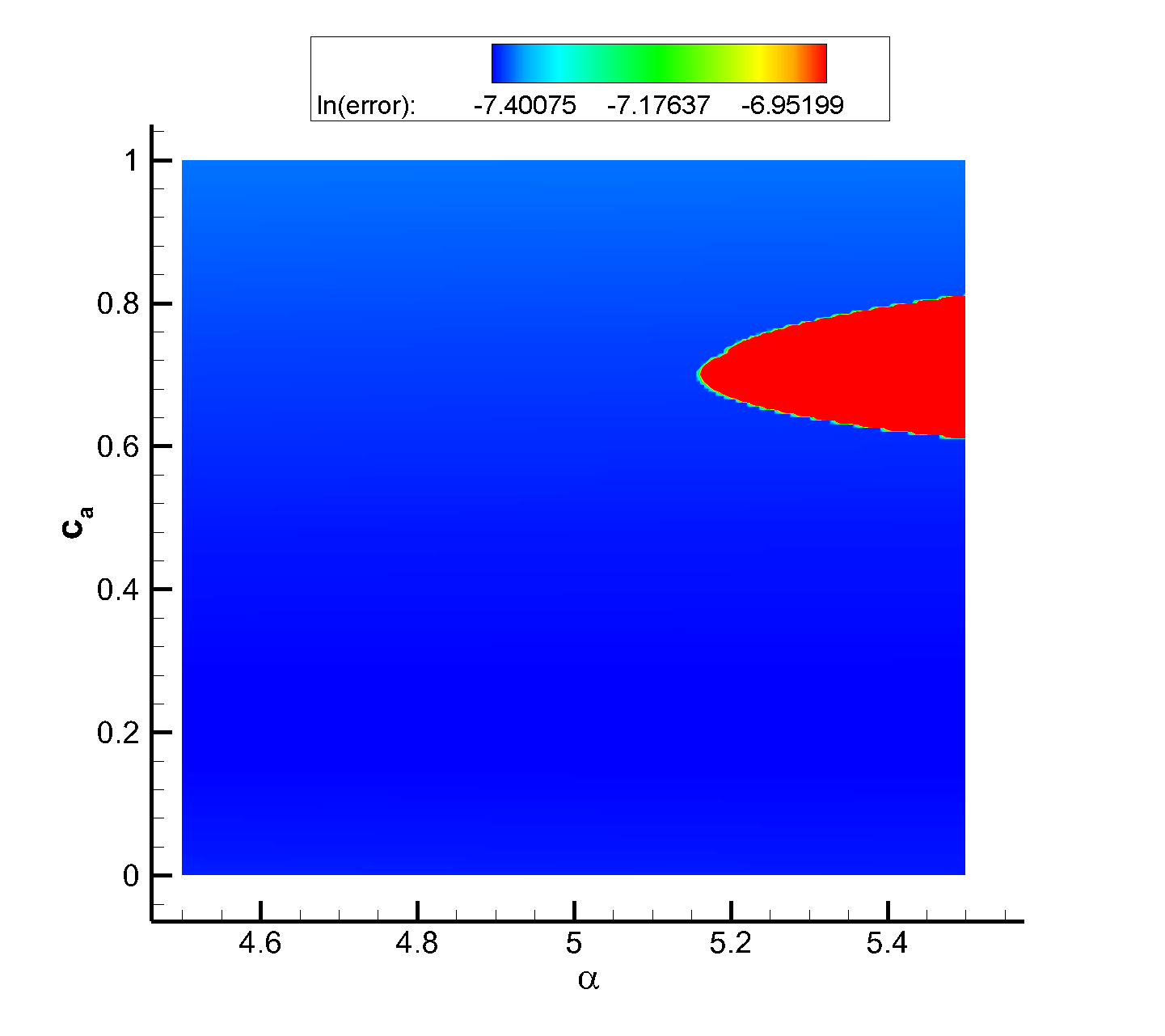}}
	\caption{Error contour of the linear problem with the fifth-order boundary treatment. }
	\label{fig:logerror2}
\end{figure}

\noindent
We use the $d$-th order upwind scheme for spatial discretization and the third-order TVD RK scheme for time discretization. Let $t_{end}=30$, $N=200$.
Take time step 
\begin{equation}\label{eq:time_step}
\Delta t=(\lambda_{cfl})_{max}\Delta x.
\end{equation}
We test the problem with $\alpha$ located in or out of the range given in Table \ref{TAB:NEWSILW_STABILITY}. The numerical results are shown in Table \ref{tab:test_result}. It can be observed that when $\alpha$ falls in the range, the scheme will be stable for all tested $C_a$. Otherwise, if $\alpha$ is out of the range, we can always find one $C_a$ such that the scheme is unstable.
	Moreover, we try to study how the error changes with varying $\alpha$ and $C_a$. 
	We fix $\Delta x = 1/25$ and take the time step as \eqref{eq:time_step}.
	Error contours of the third order and fifth order cases are plotted in Figure \ref{fig:logerror1} and Figure \ref{fig:logerror2} respectively. 
From the numerical results, we can observe that when we choose $\alpha$ inside the stable region given in Table \ref{TAB:NEWSILW_STABILITY}, the numerical error stays relatively low for any $C_a\in[0,1]$. If we choose $\alpha$ outside the interval, there would be some $C_a$ making large error or even blow up (i.e. the red regions in the figures). The results also verify the correctness of the stable interval given in Table \ref{TAB:NEWSILW_STABILITY}.

\begin{table}[htb!]
	\caption{Numerical verification results of linear stability analysis}
	\centering
	\begin{tabular}{c|c|c|c}
		\hline
		$d$ & Stable $\alpha$ in Table \ref{TAB:NEWSILW_STABILITY}& $\alpha$  &  Result\\\hline
		\multirow{2}{1cm}{3} & \multirow{2}{3cm}{[0.61,10] }  & 0.60 & Unstable for $C_a=10^{-6}$\\
		& & 1.00 & Stable for all tested $C_a$\\ \hline 
		\multirow{3}{1cm}{5} & \multirow{3}{3cm}{[0.92,5.11] }
		 & 0.91 & Unstable for $C_a=0.38$\\
		& & 1.00 & Stable for all tested $C_a$\\ 
		& & 5.12 & Unstable for $C_a=0.70$\\ \hline 
		\multirow{3}{1cm}{7} & \multirow{3}{3cm}{[1.34,1.99] } & 1.33 &  Unstable for $C_a=0.40$\\
		& & 1.50 & Stable for all tested $C_a$\\
		& & 2.00 & Unstable for $C_a=0.40$\\  \hline 
		\multirow{3}{1cm}{9} & \multirow{3}{3cm}{[1.29,2.43] }
		& 1.28 & Unstable for $C_a=0.85$\\
		& & 1.50 & Stable for all tested $C_a$\\
		& & 2.44 & Unstable for $C_a=0.03$\\ \hline 
		\multirow{3}{1cm}{11} & \multirow{3}{3cm}{[1.42,1.70]  }
		& 1.41 & Unstable for $C_a=0.93$\\
		&  & 1.50 & Stable for all tested $C_a$\\
		&  & 1.71 & Unstable for  $C_a=0.01$\\	\hline 
		\multirow{3}{1cm}{13} & \multirow{3}{3cm}{[1.49,2.08]  }
		& 1.48 & Unstable for $C_a=1-10^{-6}$\\
		& & 1.75 & Stable for all tested $C_a$\\
		& & 2.09 & Unstable for $C_a=1-10^{-6}$\\
		\hline
	\end{tabular}
	\label{tab:test_result}
\end{table}

%

\section{The new SILW method for systems}
\label{sec:system}
\setcounter{table}{0}
\setcounter{figure}{0}

\subsection{The new SILW method for one-dimensional Euler equation}
We consider the following one-dimensional compressible Euler equation:
\begin{equation}\label{1dEuler}
	\bm{U}_t+\bm{F}(\bm{U})_x=\bm{0}, \quad x\in(-1,1),\, t>0,
\end{equation}
where
\begin{equation*}
{\bm{U}}=
\left(
\begin{array}{c}
U_1   \\
U_2 \\
U_3
\end{array}
\right)=
\left(
\begin{array}{c}
\rho   \\
\rho u \\
E
\end{array}
\right),
\quad
{\bm{F}}({\bm{U}})=
\left(
\begin{array}{c}
\rho u \\
\rho u^2 +p \\
u(E+p)
\end{array}
\right).
\end{equation*}
Here, $ \rho $, $ u $, $ p $ and $ E $ represent the density, velocity, pressure and total energy per volume, respectively. In order to close the system, we give the following equation of state of ideal gas
\begin{equation*}
E=\frac{p}{\gamma -1} +\frac{1}{2} \rho u^2.
\end{equation*}
Here, $\gamma$ is the adiabatic constant, which equals to 1.4 for an ideal polytropic gas.

We consider the boundary treatment of left boundary $x=-1$ as an example. The original Euler equation \eqref{1dEuler} can be rewritten into the following nonconservative form:
\begin{equation*}
\bm{U}_t + \bm{A}(\bm{U}) \bm{U}_x = \bm{0},
\end{equation*}
where, $\bm{A}(\bm{U})=\bm{F}'(\bm{U})$ is the Jacobi matrix, 
\begin{align} \label{eqn_1d_Jacobi}
\bm{A}(\bm{U}) =
\begin{pmatrix}
0 & 1 & 0 \\
 \dfrac{1}{2} (\gamma -3) u^2 & (3 - \gamma) u & \gamma - 1 \\
 \dfrac{1}{2} (\gamma - 1) u^3 - u  H & H - (\gamma - 1) u^2 & \gamma u
\end{pmatrix}\
=\begin{pmatrix}
\bm{a}_1(\bm{U})\\
\bm{a}_2(\bm{U})\\
\bm{a}_3(\bm{U})
\end{pmatrix}\, , 
\end{align}
with the total enthalpy $H = (E + p)/\rho$.  The Jacobian matrix is diagonalizable
$$\bm{A}(\bm{U})=\bm{F}'(\bm{U}) = \bm{ R \Lambda L }.$$
Here, $\bm{\Lambda} = {\rm diag}(u - c,\, u,\, u + c)$,
$c = \sqrt{\gamma p/ \rho}$ is the speed of sound, $\bm{R}$ and $\bm{L}={\bm{R}}^{-1}$ are matrixes as follows:
\begin{align*}
\begin{aligned}
 \bm{R}(\bm{U}) & \, =
\begin{pmatrix}
1 & 1 & 1 \\
u-c & u & u+c \\
H - u c & \dfrac{1}{2} u^2 & H + u c
\end{pmatrix},
\end{aligned}
\end{align*}
\begin{align*}
\begin{aligned}
 \bm{L}(\bm{U}) & \, = \dfrac{1}{c^2}
\begin{pmatrix}
\dfrac{1}{2} uc + \dfrac{1}{4} (\gamma - 1) u^2 & - \dfrac{1}{2} (\gamma - 1) u -
\dfrac{1}{2} c
& \dfrac{1}{2 } (\gamma - 1) \\
c^2 - \dfrac{1}{2} (\gamma - 1) u^2 & (\gamma - 1) u & 1 - \gamma \\
- \dfrac{1}{2} u c + \dfrac{1}{4} (\gamma - 1) u^2 &  - \dfrac{1}{2} (\gamma - 1) u
+ \dfrac{1}{2} c
& \dfrac{1}{2} (\gamma - 1)
\end{pmatrix} \
=\begin{pmatrix}
\bm{l}_1(\bm{U}) \\
\bm{l}_2(\bm{U}) \\
\bm{l}_3(\bm{U}) \\
\end{pmatrix} \, .
\end{aligned}
\end{align*}

The number of boundary conditions we need to give is determined by the sign of the eigenvalues of Jacobi matrix $\bm{A}(\bm{U})$ at the boundary. Specifically, on the left boundary $x=-1$, it can be divided into the following cases:
\begin{itemize}
\item [Case 1.]
$u - c > 0$: three boundary conditions need to be given;
\item [Case 2.]
$u - c \leq 0$, $\,$ $u > 0$: two boundary conditions need to be given;
\item [Case 3.]
$u \leq 0$, $\,$ $ u + c > 0$: only one boundary condition needs to be given;
\item [Case 4.]
$u + c \leq 0$: no boundary conditions are required.
\end{itemize}

Here, we take the case 2 as an example to describe our algorithm. Suppose two boundary conditions are given at the left boundary,
$$U_1(-1,t)=g_1(t),\quad \text{and} \quad U_2(-1,t)=g_2(t).$$

\noindent
Again, we use the uniform mesh \eqref{eq:mesh_1D}
and employ the finite difference methods to get the semi-discrete scheme:
\begin{equation} \label{semid_system}
\frac{d}{dt}\bm{U}_j=-\frac{1}{\Delta x}(\hat{\bm{F}}_{j+1/2}-\hat{\bm{F}}_{j-1/2}), 
\quad j=0,\cdots, N.
\end{equation}
Here, $\bm{U}_j(t)$ is the approximation of $\bm{U}(x_j,t)$, and the numerical flux $\hat{\bm{F}}_{j+1/2}$ can be obtained by the WENO reconstruction. We take the fifth order scheme as an example to describe our boundary algorithm, and other high order schemes can be obtained similarly. For the fifth order WENO scheme, we need the values at three ghost points near the boundary $x = - 1 $, which are $\bm{U}_{-1}$,$\bm{U}_{-2}$ and $\bm{U}_{-3}$.

Similar to the case of the scalar equation, in order to obtain the ghost point values, we will construct the extrapolation polynomials $\bm{q}(x)$ near the boundary. It can be seen from Table \ref{TAB:NEWSILW_STABILITY} that, for the fifth order scheme, we have to get the point values and first order spatial derivatives on the boundary through the ILW procedure. 
To ensure the order of accuracy, the value $\bm{U}^{*(0)}$ and $\bm{U}^{*(1)}$ should be 5th and 4th order approximations of $\bm{U}|_{x=-1}$ and $\bm{U}_x|_{x=-1} $, respectively.

Specifically, we use the left characteristic matrix $\bm{L}=\bm{L}(\bm{U}^{ext,0})$ 
to do the characteristic projection $\bm{V}=\bm{L}\bm{U}$.
Here, $\bm{U}^{ext,0}$ is the approximation of $\bm{U}|_{x=-1}$ obtained via extrapolation directly, and 
$$\bm{V}=(V_1,V_2,V_3)^T$$
is the characteristic variable. In case 2, $V_1$ is the outflow variable, $V_2 $ and $V_3 $ are the inflow variables.
Combined with the boundary conditions, we can obtain the following linear system:
\begin{equation}
\begin{aligned}
U_1^{*(0)} =& g_1(t),\\
U_2^{*(0)} =& g_2(t),\\
\bm{l}_1 \cdot \bm{U}^{*(0)} =& V_1^{*(0)},
\end{aligned}
\end{equation}
where, $V_1^{*(0)}$ is the 5th order approximation of ${V}_1|_{x=-1}$ and can be extrapolated from the interior grid points. By solving the above system, we can get the value of $\bm{U}^{*(0)}$.

For $\bm{U}^{*(1)}$, applying the ILW procedure, we have
\begin{subequations}
	\begin{equation} 
	\begin{aligned}
\bm{a}_1(\bm{U}^{*(0)}) \cdot \bm{U}^{*(1)} = -g_1'(t),\\
\bm{a}_2(\bm{U}^{*(0)}) \cdot \bm{U}^{*(1)} = -g_2'(t),
	\end{aligned}
	\end{equation}
\text{This is combined with the outflow conditions,}
\begin{equation}
\bm{l}_1 \cdot \bm{U}^{*(1)} = V_1^{*(1)},
\end{equation}
\end{subequations}
where, $V_1^{*(1)}$ is the 4th order approximation of $({V}_1)_x|_{x=-1}$ which is obtained by extrapolating from the interior grid points. 
By combining and solving the above equations, we can get the value of $\bm{U}^{*(1)}$.

Then, we can apply the new scalar SILW method (Algorithm 2) on each component of $\bm{U}$ to reconstruct the corresponding polynomials $\bm{p}(x)$ and $\bm{q}(x)$. After that, the ghost points values can be defined as $\bm{U}_{j}=\bm{q}(x_{j})$.

Additionally, in many cases, boundary conditions are not directly given to the conserved variables. For example, the temperature or pressure might be given on the boundary. For these cases, when dealing with the boundary conditions, we need to convert the conservation equations into equations in terms of primitive variables. This process is used in \cite{TS3, Cheng, Liu1, Liu2}, which is very similar to the above process of conserved variables. We will not explain it in more details here.

\subsection{The new SILW method for two-dimensional Euler equation}
Consider the two-dimensional Euler equation as follows:
\begin{equation}\label{eq:euler_2d}
\frac{\partial \bm{\bm{U}}}{\partial t}+\frac{\partial \bm{\bm{F}(\bm{U})}}{\partial x}+\frac{\partial \bm{\bm{G}(\bm{U})}}{\partial y}=\bm{0}, \quad (x,y)^T \in \Omega,
\end{equation}
where,
\begin{equation}
\bm{U}=
\begin{pmatrix}
\rho\\
\rho u\\
\rho v\\
E\\
\end{pmatrix},\quad
\bm{F}(\bm{U})=
\begin{pmatrix}
\rho u\\
\rho u^2+p\\
\rho uv\\
u(E+p)\\
\end{pmatrix}
,\quad
\bm{G}(\bm{U})=
\begin{pmatrix}
\rho v\\
\rho uv\\
\rho v^2+p\\
v(E+p)\\
\end{pmatrix}.
\end{equation}
Here, $ \rho $, $\bm{u}=(u,v)^T$, $ p $ and $ E $ represent the density, velocity, pressure and total energy per volume, respectively. Moreover, the following equation of state of ideal gas is given to close the system, 
\begin{equation}\nonumber
E=\frac{p}{\gamma-1}+\frac{1}{2}\rho(u^2+v^2).
\end{equation}
Here, $\gamma$ is the adiabatic constant, which equals to 1.4 for an ideal polytropic gas.

We use a uniform non body-fitted Cartesian mesh to divide the domain
$$x_{i+1}=x_{i}+\Delta x, \quad y_{j+1}=y_{j}+\Delta y,$$
with mesh sizes $\Delta x$ and $\Delta y$ in $x$- and $y$-direction, respectively.
We discretize the equation into the following conservative semi-discrete scheme:
\begin{equation*}
\frac{d\bm{U}_{i,j}}{dt}+\frac{\hat{\bm{F}}_{i+\frac{1}{2},j}-\hat{\bm{F}}_{i-\frac{1}{2},j}}{\Delta x}+\frac{\hat{\bm{G}}_{i,j+\frac{1}{2}}-\hat{\bm{G}}_{i,j-\frac{1}{2}}}{\Delta y}=\bm{0},
\end{equation*}
where, $\bm{U}_{i,j}(t)$ is approximation to the exact solution $\bm{U}(x_i, y_j, t)$, $\hat{\bm{F}}_{i+\frac{1}{2},j}$ and $\hat{\bm{G}}_{i,j+\frac{1}{2}}$ are numerical fluxes, which can be obtained by WENO reconstruction.

\begin{figure}[htb!]
	\centering
	\includegraphics[width=0.5\linewidth]{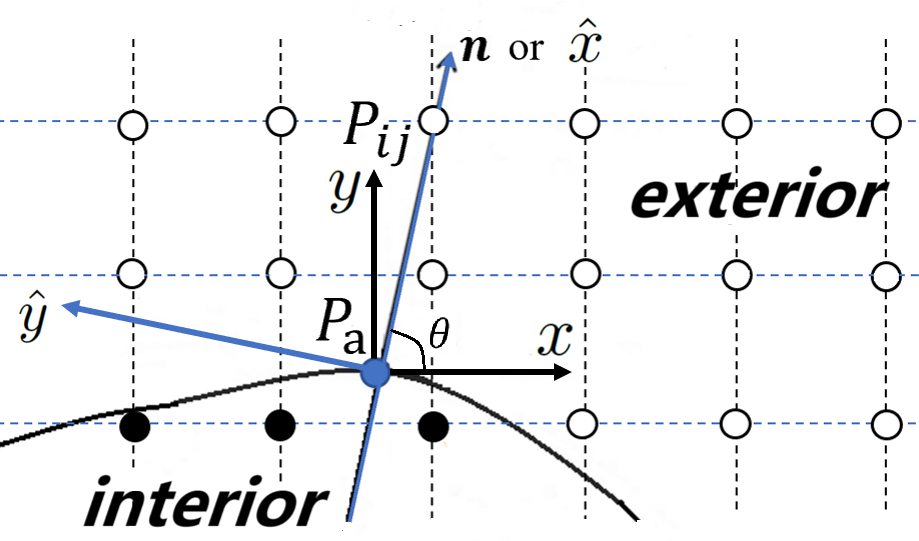}
	\caption{The local coordinate rotation diagram. }
	\label{fig:rotation}
\end{figure}

Following the idea in \cite{TS1, TWSN}, a local coordinate system is defined near the boundary to simplify the 2D boundary treatment into the 1D algorithm.
Suppose $P_{ij}=(x_i, y_j)$ is a ghost point near the boundary. At first, we find its foot point $P_a\in\partial \Omega$, so that the outward normal vector to $\partial \Omega$ at $P_a$ goes through $P_{ij}$, as shown in Figure \ref{fig:rotation}. 
The point $P_a$ is also known as boundary-intercept or body-intercept if $ \Omega$ is the boundary of a solid body in fluid–structure interaction problems.
Assume the unit normal vector from $P_a$ to $P_{ij}$ is $\mathbf{n}=(\cos\theta,\sin\theta)^T$. 
Then we perform a local coordinate rotation transformation at $P_a$,
such that $\hat{x}$-axis points in the normal direction to $\partial \Omega$ and the $\hat{y}$-axis points in the tangential direction to $\partial \Omega$,
\begin{equation*}
\begin{pmatrix}
 \hat{x}\\ \hat{y}
\end{pmatrix}
=
\begin{pmatrix}
\cos\theta&\sin\theta\\
-\sin\theta&\cos\theta
\end{pmatrix}
\begin{pmatrix}
x\\ y
\end{pmatrix}.
\end{equation*}

\noindent
In the new coordinate system, the equation \eqref{eq:euler_2d} can be rewritten as
\begin{equation}\label{eq:euler_2d_rot}
\frac{\partial \bm{\hat{\bm{U}}}}{\partial t}+\frac{\partial \bm{\bm{F}(\hat{\bm{U}})}}{\partial \hat{x}}+\frac{\partial \bm{\bm{G}(\hat{\bm{U}})}}{\partial \hat{y}}=\bm{0},
\end{equation}
where,
\begin{equation*}
\hat{\bm{U}}=
\begin{pmatrix}
\rho\\
\rho \hat{u}\\
\rho \hat{v}\\
E\\
\end{pmatrix}
=\begin{pmatrix}
\hat{U}_1\\
\hat{U}_2\\
\hat{U}_3\\
\hat{U}_4\\
\end{pmatrix}, \quad
\begin{pmatrix}
 \hat{u}\\ \hat{v}
\end{pmatrix}
=
\begin{pmatrix}
\cos\theta&\sin\theta\\
-\sin\theta&\cos\theta
\end{pmatrix}
\begin{pmatrix}
u\\ v
\end{pmatrix}.
\end{equation*}
Let
\begin{equation*}
\bm{A}(\hat{\bm{U}})=\bm{F}'(\hat{\bm{U}})=
\begin{pmatrix}
\bm{a}_1(\hat{\bm{U}})\\
\bm{a}_2(\hat{\bm{U}})\\
\bm{a}_3(\hat{\bm{U}})\\
\bm{a}_4(\hat{\bm{U}})
\end{pmatrix},
\quad
\bm{Res}=-\frac{\partial \bm{\bm{G}(\hat{\bm{U}})}}{\partial \hat{y}}=
\begin{pmatrix}
Res_1\\
Res_2\\
Res_3\\
Res_4
\end{pmatrix}.
\end{equation*}

\noindent
Then, the equations can be written in the following non conservative form
\begin{equation}\label{eq:euler_2d_non}
\hat{\bm{U}}_t + \bm{A}(\hat{\bm{U}}) \hat{\bm{U}}_x = \bm{Res}.
\end{equation}

\noindent
The original Euler equation is hyperbolic, so $\bm{A}(\hat{\bm{U}})$ is diagonalizable:
$$\bm{A}(\hat{\bm{U}}) = \bm{R}(\hat{\bm{U}}) \bm{\Lambda}(\hat{\bm{U}}) \bm{L}(\hat{\bm{U}}),$$
Here,
$$\bm{\Lambda}(\hat{\bm{U}})=diag(\hat{u}-c,\hat{u},\hat{u},\hat{u}+c)$$
\begin{equation*}
\bm{L}(\hat{\bm{U}})=
\begin{pmatrix}
\bm{l}_1(\hat{\bm{U}})\\
\bm{l}_2(\hat{\bm{U}})\\
\bm{l}_3(\hat{\bm{U}})\\
\bm{l}_4(\hat{\bm{U}})
\end{pmatrix}, \quad 
\bm{R}(\hat{\bm{U}}) = \bm{L}^{-1}(\hat{\bm{U}}).
\end{equation*}

The number of boundary conditions that should be given at the boundary point $P_a$ is related to the eigenvalues $\hat{u}-c,\hat{u},\hat{u},\hat{u}+c$ at this point. Specifically, it can be divided into the following situations:
\begin{itemize}
\item [Case 1:]
$\hat{u} - c \geq 0$ , no boundary condition is required;
\item [Case 2:]
$\hat{u} - c < 0$, $\,$ $\hat{u} \geq 0$ , only one boundary condition needs to be given;
\item [Case 3:]
$\hat{u} < 0$, $\,$ $ \hat{u} + c \geq 0$ , three boundary conditions need to be given;
\item [Case 4:]
$\hat{u} + c < 0$ , four boundary conditions need to be given.
\end{itemize}

We take case 2 as an example to describe our algorithm. Suppose the boundary condition given at the boundary point $P_a$ is
$$\hat{U}_2=g(t).$$

With the help of the local coordinate rotation transformation, we will apply the 1D boundary treatment along the $\hat{x}$-direction based on equation \eqref{eq:euler_2d_non}.
It can be seen from Table \ref{TAB:NEWSILW_STABILITY} that, for the fifth order scheme, we need to get the 0th and 1st order normal direction derivatives on the boundary from the ILW procedure when constructing the extrapolation polynomial $\bm{q}(s)$ along the $\hat{x}$-direction. That is, we need to get the values of $\hat{\bm{U}}^{*(0)}$ and $\hat{\bm{U}}^{*(1)}$, which are the 5th and 4th order approximations of $\hat{\bm{U}}|_{P_a}$ and $\hat{\bm{U}}_{\hat{x}}|_{P_a}$ respectively, through the ILW procedure.

Specifically, we use the left characteristic matrix $\bm{L}=\bm{L}(\hat{\bm{U}}^{ext,0})=(\bm{l}_1,\bm{l}_2,\bm{l}_3, \bm{l}_4)^T$ to do the characteristic projection $\bm{V}=\bm{L}\hat{\bm{U}}$. Here,
$\hat{\bm{U}}^{ext,0}$ is the extrapolation value at $P_a$, and $\bm{V}=(V_1,V_2,V_3,V_4)^T$ is the characteristic variable. For case 2, $V_2,V_3,V_4$ are the outflow variables, $V_1$ is the inflow variable. 
Similar to the 1D system, we can give a system combining the given boundary conditions and extrapolation on outgoing variables to get $\hat{\bm{U}}^{*(0)}$ and $\hat{\bm{U}}^{*(1)}$.
To be specific, for $\hat{\bm{U}}^{*(0)}$, we have that
\begin{equation}
\begin{aligned}
\hat{U}_2^{*(0)} &= g(t),\\
\bm{l}_2 \cdot \hat{\bm{U}}^{*(0)} &= V_2^{*(0)},\\
\bm{l}_3 \cdot \hat{\bm{U}}^{*(0)} &= V_3^{*(0)},\\
\bm{l}_4 \cdot \hat{\bm{U}}^{*(0)} &= V_4^{*(0)},
\end{aligned}
\end{equation}
where, $V_2^{*(0)}$,$V_3^{*(0)}$ and $V_4^{*(0)}$ can be extrapolated from the interior grid points. Then we can get the value of $\hat{\bm{U}}^{*(0)}$ by solving the above system.
For $\hat{\bm{U}}^{*(1)}$, apply the ILW procedure and we have
\begin{subequations}
	\begin{equation} 
	\bm{a}_2(\hat{\bm{U}}^{*(0)}) \cdot \hat{\bm{U}}^{*(1)} = -g'(t)+Res_2,
	\end{equation}
\text{where, $Res_2$ can be obtained from extrapolation as well. For the outflow variables,}
	\begin{equation}
	\begin{aligned}
		\bm{l}_2 \cdot \hat{\bm{U}}^{*(1)} = V_2^{*(1)},\\
		\bm{l}_3 \cdot \hat{\bm{U}}^{*(1)} = V_3^{*(1)},\\
		\bm{l}_4 \cdot \hat{\bm{U}}^{*(1)} = V_4^{*(1)},
	\end{aligned}
	\end{equation}
\end{subequations}
where, $V_2^{*(1)}$,$V_3^{*(1)}$ and $V_4^{*(1)}$ can be extrapolated from the interior grid points. Solving the above equations will give us the value of $\hat{\bm{U}}^{*(1)}$.

\begin{figure}[htb!]
	\centering
	\includegraphics[width=0.5\linewidth]{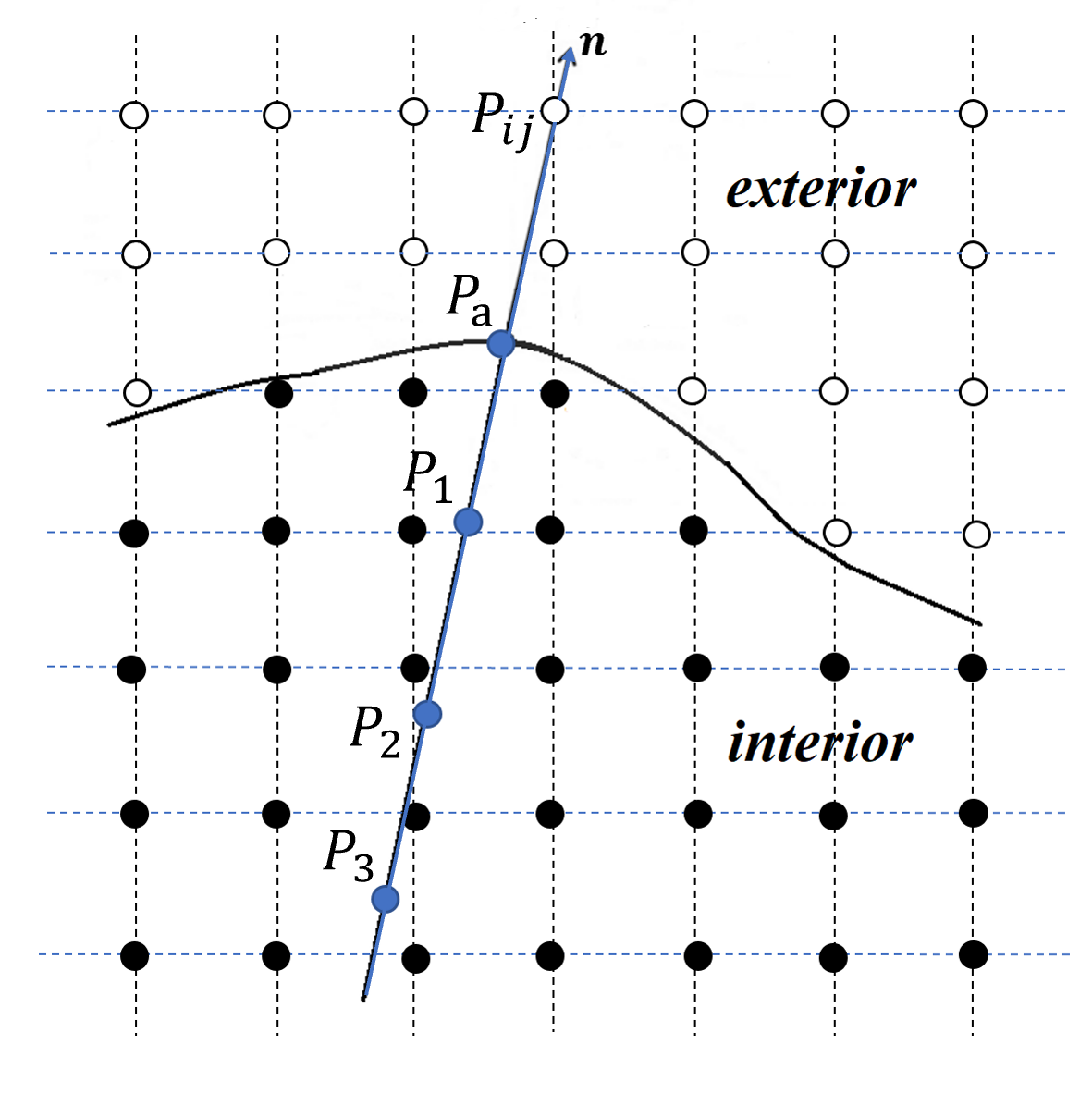}
	\caption{Two dimensional new SILW method diagram}
	\label{fig:2dsilw}
\end{figure}

Then, we can do the proposed SILW method along the $\hat{x}$-direction to construct the extrapolation polynomial $\bm{q}(s)$ and define 
$$\bm{U}_{ij}=\bm{q}(|P_{ij}-P_a|).$$
Note that due to the complex geometry, in the first step, we get the approximation polynomial $\bm{p}(x,y)$ of degree at most $4$ by least square method with the values at internal grid points near $P_a$ rather than interpolation. The artificial auxiliary points in 2D are defined as
$$P_k=P_a-k \alpha \delta \bm{n},\quad k=1, 2, 3.$$
As show in Figure \ref{fig:2dsilw}, $\{P_k\}_{k=1}^{3}$ are some non grid points in the interior area on the normal line.
Here, $\delta=\sqrt{\Delta x^2+\Delta y^2}$, and $\alpha$ can be chosen as any number in $[0.92\frac{\max(\Delta x,\Delta y)}{\delta},$$ 5.11\frac{\min(\Delta x,\Delta y)}{\delta}]$. We should point out that, in order to ensure the interval of $\alpha$ is not null, we would need to require that the difference between $\Delta x$ and $\Delta y$ is not too large.

It is also noted that for the non penetrating free slip boundary condition
$$\bm{u} \cdot \bm{n}=0,$$
we actually get $g(t)=0$ in the description above.
As before, for the problems with changing wind direction, the above inverse Lax-Wendroff procedure may involve solving an ill-conditioned linear algebraic system, which may ruin the accuracy or even lead to blowing up. There are two ways to deal with this problem. One is mentioned in \cite{TS1}, which adds additional extrapolation equations and solves a least squares problem whenever one of the eigenvalues is very closed to 0.
The other method is proposed in \cite{Lu2}, which evaluates the solution values and the flux values at ghost points separately. In this paper, we use the first technique in our numerical tests. 

\section{Numerical tests}\label{sec:num_tests}
\setcounter{table}{0}
\setcounter{figure}{0}

We take some numerical tests to show the efficiency and stability of our new proposed SILW method.
We use the third order and fifth order finite difference WENO schemes for the spatial discretization. Correspondingly, the new SILW boundary treatment with third order and fifth order accuracy will be coupled, respectively.
For all the one-dimensional numerical tests we take the parameter $\alpha = 1.0$, while for all the two-dimensional numerical tests we take $\alpha = 1.25$.
The third order TVD RK scheme \eqref{eq:RK} is employed for time discretization,
with the time step
$$\Delta t=\text{CFL}\frac{\Delta x^{k/3}}{a_x} $$
for one-dimensional problems, and 
$$\Delta t=\frac{\text{CFL}}{a_x/\Delta x^{k/3}+a_y/\Delta y^{k/3}}.$$
for two-dimensional problems. Here, the index $k/3$ help us to guarantee $k$-th order in time. 
$a_x=\max_{\bm{U}} |\lambda(\bm{F}'(\bm{U}))|$, $a_y=\max_{\bm{U}} |\lambda(\bm{G}'(\bm{U}))|$, and $\lambda$ is the eigenvalue of the Jacobian matrix. 
Throughout our numerical tests, the CFL number is taken as 0.6.

\vspace{0.3cm}
\noindent
\textbf{Example 1.}
At first, we consider the accuracy test of the new SILW on the one-dimensional Euler equation on the computational domain as $[-\pi,\pi]$. We choose suitable initial values and  boundary conditions such that the exact solution is:
\begin{equation}\label{ce1}
  \left\{\begin{aligned}
  \rho(x,t)&=1-0.2\sin(2t-x),\\
  u(x,t)&=2,\\
  p(x,t)&=2.\end{aligned}\right.
\end{equation}

%

In this example, we have $u\pm c\geq0$. Hence, all variables are outgoing on the right boundary $x=\pi$ and only extrapolation process with appropriate accuracy is used.
Meanwhile, three Dirichlet boundary conditions should be imposed on the left boundary $x=-\pi$ and the proposed SILW method works.
In this example, we fix $C_b=0.7$ and test with two
extreme choices $C_a=0.0001, 0.9999$, to verify the applicability of our algorithm to avoid the ``cut cell" instability.    
The computational errors about density $\rho$ at final time $t_{end}=1$ are shown in Tables \ref{eg1t1} - \ref{eg1t2}. We can see that for all cases, the schemes are stable and can achieve the designed order accuracy with mesh refinements.

\begin{table}[!htb]
\centering
{\caption{Example 1: errors and orders of accuracy of $\rho$ with third order scheme. }
\begin{tabular}{c|cccc|cccc}
\hline
& \multicolumn{4}{c|}{$C_a=0.0001, C_b=0.7$} 
& \multicolumn{4}{c}{$C_a=0.9999, C_b=0.7$} \\\cline{2-9}
$N$ &$L^1$ error &order &$L^{\infty}$ error& order &$L^1$ error &order &$L^{\infty}$ error& order \\
	\hline
20 &  1.67E-004 &   --      &
  5.73E-004 &   --      &  1.45E-004 &   --      &
  4.80E-004 &   --      \\
          40 &  1.58E-005 &   3.40     &
  4.78E-005 &   3.58      &  1.00E-005 &   3.85     &
  4.39E-005 &   3.45      \\
          80 &  2.07E-006 &   2.92     &
  5.30E-006 &   3.17     &	8.52E-007 &   3.56      &
  4.57E-006 &   3.26      \\
         160 &  2.69E-007 &   2.94     &
  7.30E-007 &   2.86     &	8.38E-008 &   3.34     &
  5.24E-007 &   3.12    \\
         320 &  3.43E-008 &   2.97     &
  9.50E-008 &   2.94     &	9.15E-009 &   3.19      &
  6.27E-008 &   3.06      \\
         640 &  4.32E-009 &   2.98     &
  1.21E-008 &   2.96     &  1.06E-009 &   3.10      &
  7.66E-009 &   3.03    \\
	\hline
  \end{tabular}
  \label{eg1t1}}
\end{table}

\begin{table}[!htb]
\centering
{\caption{Example 1: errors and orders of accuracy of $\rho$ with fifth order scheme. }
\begin{tabular}{c|cccc|cccc}
\hline
& \multicolumn{4}{c|}{$C_a=0.0001, C_b=0.7$} 
& \multicolumn{4}{c}{$C_a=0.9999, C_b=0.7$} \\\cline{2-9}
$N$ &$L^1$ error &order &$L^{\infty}$ error& order &$L^1$ error &order &$L^{\infty}$ error& order \\
	\hline
          20 &  9.33E-005 &    --    &
  1.76E-004 &   --     &  7.41E-005 &   --     &
  1.36E-004 &   --      \\
          40 &  2.99E-006 &   4.96     &
  6.09E-006 &   4.84      &  2.62E-006 &   4.81      &
  5.47E-006 &   4.63      \\
          80 &  9.28E-008 &   5.01     &
  1.99E-007 &   4.93      &  8.65E-008 &   4.92      &
  1.81E-007 &   4.91     \\
         160 &  2.88E-009 &   5.01    &
  6.06E-009 &   5.04     &  2.77E-009 &   4.96     &
  5.78E-009 &   4.97     \\
         320 &  8.89E-011 &   5.01     &
  1.77E-010 &   5.09      &  8.72E-011 &   4.99     &
  1.73E-010 &   5.05    \\
	\hline
  \end{tabular}
  \label{eg1t2}}
\end{table}

\vspace{0.3cm}
\noindent
\textbf{Example 2.}
Next, we consider the example given in \cite{LSQ} to test the accuracy of our method. The governing equation 
is still the one-dimensional compressible Euler equation, with following initial condition:
\begin{equation}\label{eg2IC}
  \left\{\begin{aligned}
  \rho(x,0)&=\frac{1+0.2\sin(x)}{2\sqrt{3}},\\
  u(x,0)&=\sqrt{\gamma}\rho(x,0),\\
  p(x,0)&=\rho(x,0)^{\gamma}.\end{aligned}\right.
\end{equation}
The computational domain is taken as $[0,2\pi]$.
We choose the parameter $\gamma=3$. Consequently, the exact solution is
\begin{equation*}
  \rho(x,t)=\frac{\mu(x,t)}{2\sqrt{3}}, \quad
  u(x,t)=\sqrt{\gamma}\rho(x,t), \quad
  p(x,t)=\rho(x,t)^{\gamma},
\end{equation*}
where $\mu(x,t)$ is the solution of the following Burgers' equation:
\begin{equation}\label{Burgers1d}
  \left\{
  \begin{aligned}
  &\mu_t+(\frac{\mu^2}{2})_x=0,\quad 0<x<2\pi, \quad t>0, \\
  &\mu(x,0)=1+0.2\sin(x), \quad 0\leq x\leq 2\pi.
  \end{aligned}\right.
\end{equation}
We take boundary conditions of the Euler equation from the exact solution of the initial value problem \eqref{Burgers1d} with periodic boundary conditions whenever needed. 
 
We consider the extrema situation and set $C_a=0.0001, C_b=0.9999$. The computational errors about the density $\rho$ and orders of accuracy at time $t_{end}=3.0$ are shown in Table \ref{eg2t1}, indicating that our methods can achieve the designed third order or fifth order accuracy.

\begin{table}[!htb]
\centering
{\caption{Example 2: errors and orders of accuracy of $\rho$. }
\begin{tabular}{c|cccc|cccc}
\hline
& \multicolumn{4}{c|}{third order scheme} 
& \multicolumn{4}{c}{fifth order scheme} \\\cline{2-9}
$N$ &$L^1$ error &order &$L^{\infty}$ error& order &$L^1$ error &order &$L^{\infty}$ error& order \\
	\hline
         40 &  1.30E-003  & --
         &  2.49E-003  & --	&  4.01E-004   &--
	             &  1.04E-003   &--\\
          80 &  2.04E-004  & 2.67
          &  6.20E-004  & 2.00	&  2.24E-005 &  4.15
	          &  8.03E-005   & 3.69\\
         160 &  2.84E-005  & 2.84
         &  8.50E-005  & 2.86	&  7.16E-007   & 4.97
	         &  2.98E-006   & 4.75\\
         320 &  3.60E-006  & 2.97
         &  9.93E-006  & 3.09	&  2.38E-008  & 4.90
	         &  1.00E-007  & 4.89\\
         640 &  4.52E-007  & 2.99
         &  1.18E-006  & 3.06	&  8.04E-010  & 4.89
	         &  2.90E-009  & 5.10\\
	\hline
  \end{tabular}
  \label{eg2t1}}
\end{table}
 
%

\begin{figure}
\centering
\subfigure[Third order scheme. ]{
		\includegraphics[width=0.45\textwidth]{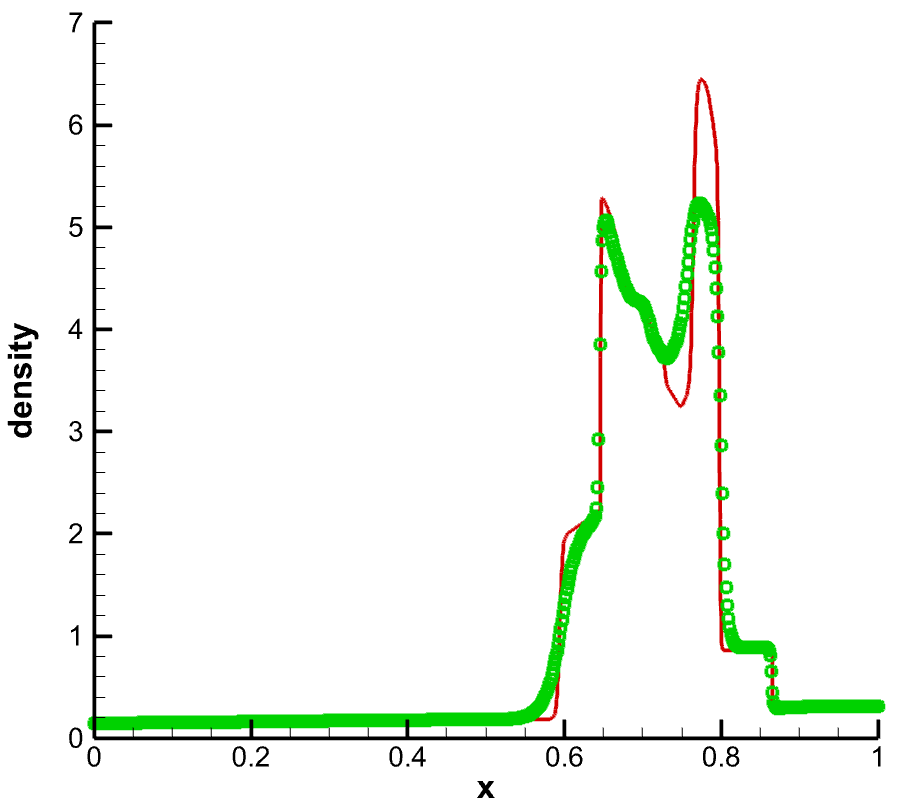}}
	\subfigure[Fifth order scheme. ]{
		\includegraphics[width=0.45\textwidth]{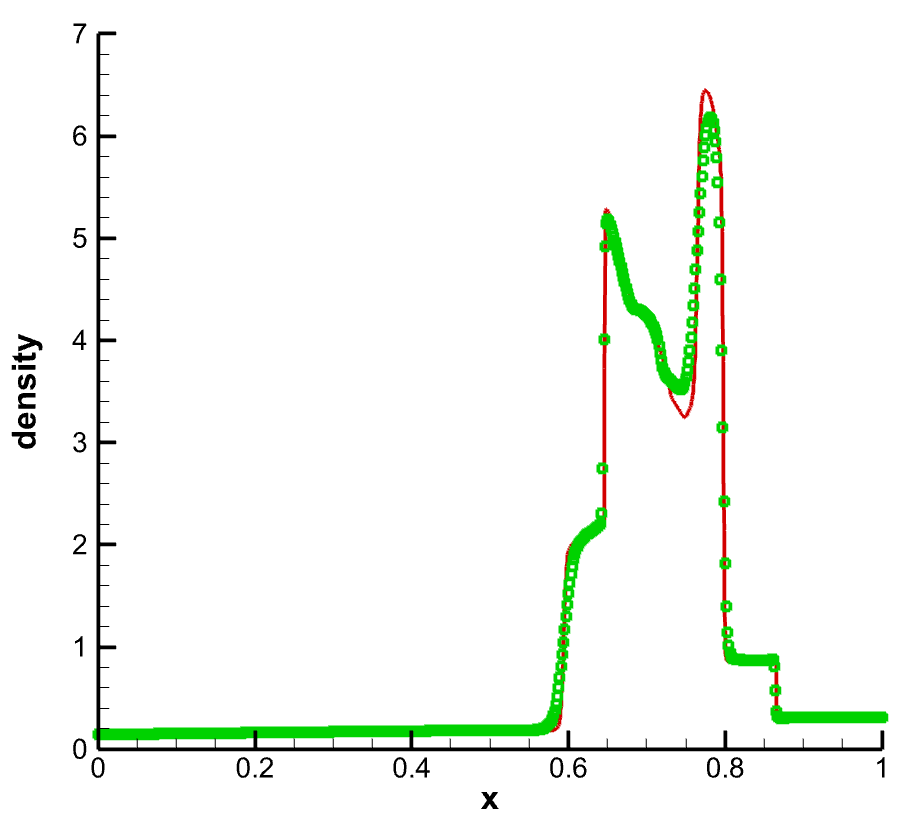}}
\caption{Example 3: Density profiles. $h=1/640$. The solid line represents the reference solution and the circle represents the numerical solution. }
\label{bla}
\end{figure}

\vspace{0.3cm}
\noindent
\textbf{Example 3.}
Now we consider the interaction of two blast waves \cite{TS1}. In this problem, multiple reflections occur between shock and rarefaction off the walls. The initial condition is
\begin{equation}\label{ce2}
  \bm{U}(x,0)
  =\left\{
  \begin{aligned}
  &\bm{U}_L,\quad x<0.1,\\
  &\bm{U}_M,\quad 0.1<x<0.9,\\
  &\bm{U}_R,\quad x>0.9.
  \end{aligned}\right.
\end{equation}
Here, $\rho_L=\rho_M=\rho_R=1$, $u_L=u_M=u_R=0$, $p_L=10^3,p_M=10^{-2}$, and $p_R=10^2$. 
Solid wall boundary conditions are used at $x=0$ and $x=1$. 
We take $t_{end}=0.038$ and $C_a=0.0001, C_b=0.7$. 
At the same time, we use a very dense grid with $\Delta x=1/2560$ and the original ILW method to obtain the reference solution. 
The numerical results are shown in Figure \ref{bla}. We can see that the new ILW method can distinguish the structure of the solution well, and higher order scheme has a better approximation to the complex structure.  

\begin{figure}[htb!]
	\centering
	\includegraphics[width=0.5\linewidth]{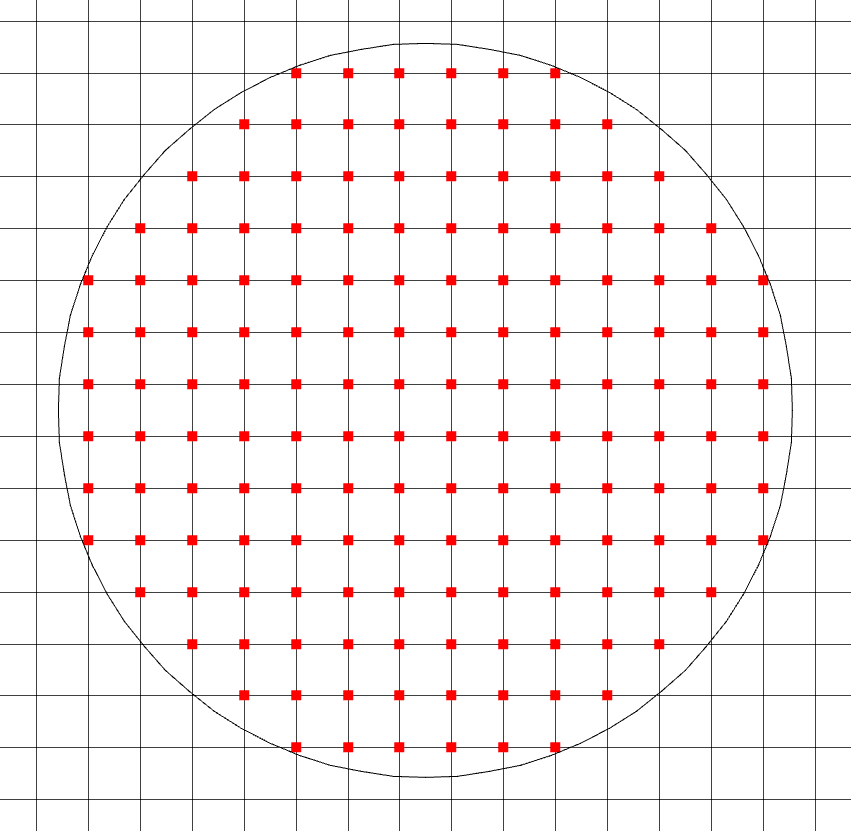}
	\caption{Example 4: Non body fitted Cartesian mesh. The red points are the interior points.}
	\label{fig:gridex4}
\end{figure}

\vspace{0.3cm}
\noindent
\textbf{Example 4.}
We consider the two-dimensional linear scalar equation on a disk:
\begin{equation}
u_t+u_x+u_y=0,\quad (x,y)^T \in \Omega=\{(x,y):x^2+y^2<0.5\}.
\end{equation}
The initial condition is given as $$u(x,y,0)=0.25+0.5\sin[\pi(x+y)],$$
and the boundary is given whenever needed such that the exact solution is 
$$u(x,y,t)=0.25+0.5\sin[\pi(x+y-2t)].$$
The domain is discretized by embedding the domain in a regular Cartesian mesh with $x_i=(i-\frac{1}{2})\Delta x, y_j=(j-\frac{1}{2})\Delta y$, and the non body fitted Cartesian mesh $h=\Delta x= \Delta y$. We show Figure \ref{fig:gridex4} as an example.
The final time is taken as $t_{end}=1.0$. The numerical results are given in Table \ref{eg4t1}, indicating that our schemes are stable and can achieve the designed order of accuracy.

%

\begin{table}[!htb]
\centering
\caption{Example 4: errors and orders of accuracy . }
\begin{tabular}{c|cccc|cccc}
\hline
& \multicolumn{4}{c|}{third order scheme} 
& \multicolumn{4}{c}{fifth order scheme} \\\cline{2-9}
$h$ &$L^1$ error &order &$L^{\infty}$ error& order &$L^1$ error &order &$L^{\infty}$ error& order \\
	\hline
         1/10 &  1.28E-004 &   --     
                     &  4.83E-004 &   --&  4.16E-004 &   --     
          &  1.42E-003 &   --\\     
                    1/20 &  1.33E-005 &  3.26    
                     &  4.54E-005 &   3.41 &  1.51E-005 &   4.77    
          &  1.38E-004 &   3.36\\     
                    1/40 &  1.43E-006 &  3.21    
                     &  5.84E-006 &   2.95 &  3.47E-007 &   5.44     
          &  5.54E-006 &   4.64\\    
                   1/80 &  1.46E-007 &   3.28   
                    &  6.68E-007 &   3.12  &  1.17E-008 &   4.89   
         &  2.71E-007 &   4.35\\   
                   1/160 &  1.19E-008 &   3.61    
                    &  8.82E-008 &  2.92&  4.47E-010 &   4.70     
         &  9.72E-009 &   4.80\\ 
	\hline
  \end{tabular}
  \label{eg4t1}
\end{table}

%

\vspace{0.3cm}
\noindent
\textbf{Example 5.}
We test the vortex evolution problem for the 2D Euler equation with $\gamma=1.4$. The mean flow is $\rho=u=v=p=1$ with following isentropic vortex perturbation centered at $(x_0, y_0)=(0,0)$ (perturbation in $(u, v)$ and temperature $T=p/\rho$ , no perturbation in the entropy $S=p/\rho^{\gamma}$):
\begin{equation}\label{ce3}
 \begin{aligned} 
 (\delta u,\delta v)=& \frac{\epsilon}{2\pi}e^{0.5(1-r^2)}(-\bar{y},\bar{x}), \\
  \delta T= & -\frac{(\gamma-1)\epsilon^2}{8\gamma \pi^2}e^{(1-r^2)} , \\
  \delta S= & 0.
  \end{aligned}
\end{equation}
where $(\bar{x},\bar{y})=(x-x_0,y-y_0)$, $r^2=\bar{x}^2 + \bar{y}^2$,  and the vortex strength $\epsilon=5$. 
It is clear that the exact solution of the corresponding Cauchy problem is just the passive convection of the
vortex with the mean velocity, 
i.e., $\bm{U}(x,y,t) = \bm{U}(x-t,y-t,0)$. 
Here, the computational domain is taken as $(-0.5,1)\times(-0.5,1)$ and the final time is taken as $t_{end}=1.0$. The boundary conditions are taken from the exact solution
whenever needed. 
We divide the domain with the uniform Cartesian mesh $x_i=(i-\frac{1}{2})h$ and $y_j=(j-\frac{1}{2})h$,
with mesh size $h=1.5/N$.
The numerical results in Table \ref{eg5t1} show that the schemes are stable and can reach the designed high order. 

\begin{table}[!htb]
\centering
{\caption{Example 5: errors and orders of accuracy of $\rho$ . }
\begin{tabular}{c|cccc|cccc}
\hline
& \multicolumn{4}{c|}{third order scheme} 
& \multicolumn{4}{c}{fifth order scheme} \\\cline{2-9}
$h$ &$L^1$ error &order &$L^{\infty}$ error& order &$L^1$ error &order &$L^{\infty}$ error& order \\
	\hline
         3/40 &  1.73E-004 &   --
                     & 3.05E-004 &   --& 3.39E-005 &   --
		 &  7.12E-005 &   --\\
                     3/80 &  2.17E-005 &   2.99
                     &  4.10E-005 &   2.89&  1.09E-006 &   4.95
		 &  2.33E-006 &   4.93\\
                     3/160 &  2.51E-006 &   3.10
                     &  4.93E-006 &   3.05&  3.46E-008 &   4.98
		 &  1.08E-007 &   4.43\\
                     3/320 &  2.99E-007 &   3.07
                     &  6.32E-007 &   2.96&  1.12E-009 &   4.94
		 &  4.71E-009&   4.51\\
                     3/640 &  3.63E-008 &   3.04
                     &  8.34E-008 &  2.92&  3.77E-011 &   4.89
		 &  1.90E-010 &   4.63\\
	\hline
  \end{tabular}
  \label{eg5t1}}
\end{table}

%
%

\vspace{0.3cm}
\noindent
\textbf{Example 6.}
Next, we consider the 2D version of Example 2 \cite{LSQ}. The governing equation is the two-dimensional compressible Euler equations with following initial condition:
\begin{equation}\label{eg6IC}
  \left\{\begin{aligned}\rho(x,y,0)&=\frac{1+0.2\sin(\frac{x+y}{2})}{\sqrt{6}},\\
  u(x,y,0)&=v(x,y,0)=\sqrt{\frac{\gamma}{2}}\rho(x,y,0),\\
  p(x,y,0)&=\rho(x,y,0)^{\gamma}.\end{aligned}\right.
\end{equation}
We choose the parameter $\gamma=3$, such that the exact solution is
\begin{equation*}
  \rho(x,y,t)=\frac{\mu(x,y,t)}{\sqrt{6}}, \quad
  u(x,y,t)=v(x,y,t)=\sqrt{\frac{\gamma}{2}}\rho(x,y,t), \quad
  p(x,y,t)=\rho(x,y,t)^{\gamma},
\end{equation*}
where $\mu(x,y,t)$ is the solution of the following 2D Cauchy problem for Burgers' equation:
\begin{equation}\label{Burgers2d}
  \left\{
  \begin{aligned}
  &\mu_t+(\frac{\mu^2}{2})_x+(\frac{\mu^2}{2})_y=0, \\
  &\mu(x,y,0)=1+0.2\sin(\frac{x+y}{2}).
  \end{aligned}\right.
\end{equation}
We consider following two computational domains:
\begin{subequations}
\begin{align}
\Omega = & [0,4\pi]\times[0,4\pi], \label{eq:domain1}\\
\Omega = & \{(x,y)|x^2+y^2<(1.5\pi)^2\},\label{eq:domain2}
\end{align}
\end{subequations}
and take boundary conditions from the exact solution whenever needed.
For the square domain \eqref{eq:domain1}, we use a grid similar to Example 5. And for the circular domain \eqref{eq:domain2}, we use a non body fitted grid similar to Example 4.
The numerical results at the final time $t_{end}=1$ are shown in Tables \ref{eg6t1} - \ref{eg6t2}. We can see that our schemes are stable and high order accuracy for all cases.

\begin{table}[!htb]
\centering
{\caption{Example 6: The errors and the orders of accuracy of $\rho$ on the square domain $\Omega = [0,4\pi]\times[0,4\pi]$. }
\begin{tabular}{c|cccc|cccc}
\hline
& \multicolumn{4}{c|}{third order scheme} 
& \multicolumn{4}{c}{fifth order scheme} \\\cline{2-9}
$h$ &$L^1$ error &order &$L^{\infty}$ error& order &$L^1$ error &order &$L^{\infty}$ error& order \\
	\hline
         4$\pi$/100 &  4.09E-004 &   --    
		  &  1.13E-005 &   --&  1.42E-005 &   --    
		   &  2.08E-006 &   --\\   
		  4$\pi$/150 &  1.22E-004 &   2.98   
		  &  3.51E-006 &   2.89& 1.93E-006 &   4.92     
		   &  2.89E-007 &   4.86\\  
		  4$\pi$/200 &  5.18E-005 &   2.98    
		  &  1.66E-006 &   2.60&  4.59E-007 &   4.99    
		   &  7.55E-008 &   4.67\\     
		  4$\pi$/250 &  2.66E-005 &   2.98     
		  &  9.14E-007 &   2.67&  1.49E-007 &  5.01   
		   &  2.77E-008 &   4.48\\       
		  4$\pi$/300 &  1.54E-005 &   2.99    
		  &  5.54E-007 &   2.74&  6.01E-008 &   5.01   
		   & 1.20E-008 &   4.56\\
		   4$\pi$/350
		   & 9.72E-006 & 2.99
		   & 3.61E-007 &   2.79
		   &  2.78E-008 &   5.00   
		    & 5.98E-009 &   4.56\\
		   4$\pi$/400 &  6.52E-006 &   2.99    
		   		  &  2.48E-007 & 2.81 &  1.43E-008 &   4.98   
		   		   & 3.43E-009 &   4.16\\   	\hline
  \end{tabular}
  \label{eg6t1}}
\end{table}

\begin{table}[!htb]
\centering
\caption{Example 6: The errors and the orders of accuracy of $\rho$ on the circular domain $\Omega = \{(x,y)|x^2+y^2<(1.5\pi)^2\}$. }
\begin{tabular}{c|cccc|cccc}
\hline
& \multicolumn{4}{c|}{third order scheme} 
& \multicolumn{4}{c}{fifth order scheme} \\\cline{2-9}
$h$ &$L^1$ error &order &$L^{\infty}$ error& order &$L^1$ error &order &$L^{\infty}$ error& order \\
	\hline
         4$\pi$/100 &  1.46E-004 &   --    
 &  1.14E-005 &   --&  1.48E-005 &   --     
&  8.59E-006 &   --\\    
4$\pi$/150 &  4.20E-005 &  3.08    
 &  4.28E-006 &   2.42&  2.13E-006 &   4.78     
&  1.28E-006 &   4.69\\       
4$\pi$/200 &  1.81E-005 &   2.91    
 &  1.79E-006 &   3.03&  6.65E-007 &  4.05    
&  3.51E-007 &   4.49\\     
4$\pi$/250 &  9.20E-006 &   3.05     
 &  8.90E-007 &   3.12&  2.21E-007 &   4.91   
&  1.16E-007 &   4.95\\    
4$\pi$/300 &  5.44E-006 &   2.87    
 &  5.23E-007 &   2.91&  8.42E-008 &   5.31   
&  4.72E-008 &   4.93\\
4$\pi$/350 &  3.40E-006 &   3.05    
 &  3.37E-007 &   2.85&  4.41E-008 &   4.20   
&  2.51E-008 &   4.10\\   
4$\pi$/400 &  2.23E-006 &   3.16    
 &  2.31E-007 &   2.84&  2.39E-008 &   4.59   
&  1.44E-008 &   4.19\\   \hline
  \end{tabular}
  \label{eg6t2}
\end{table}

\begin{figure}[htb!]
	\centerline{
	\includegraphics[width=0.5\textwidth]{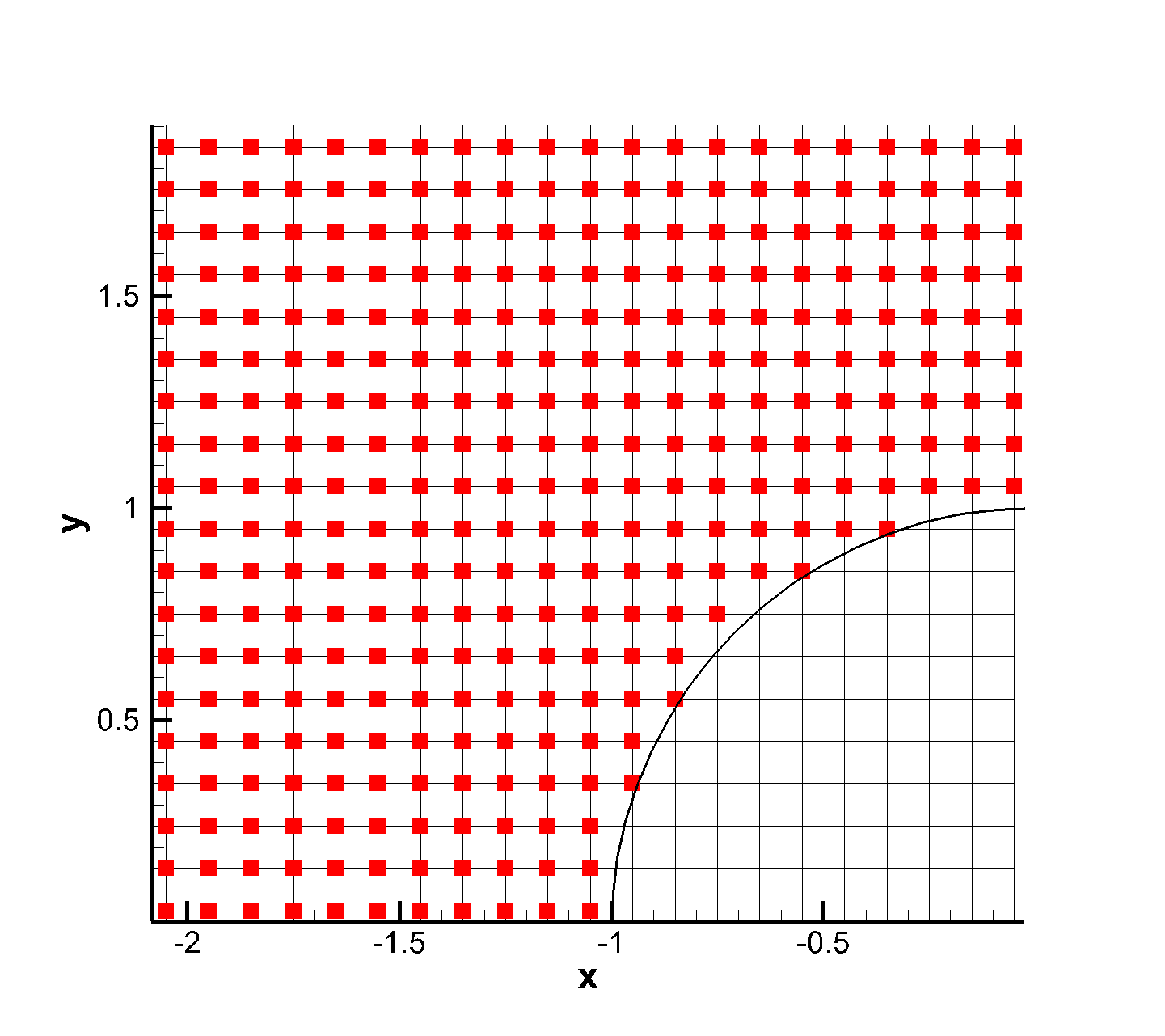}}
	\caption{Example 7: The non body-fitted Cartesian mesh near the cylinder boundary. The red points are the interior points.
	}
	\label{grideg7}
\end{figure}

\begin{figure}
	\subfigure[Pressure contour]{
	\includegraphics[width=1.0in,angle=0,scale=2.7]{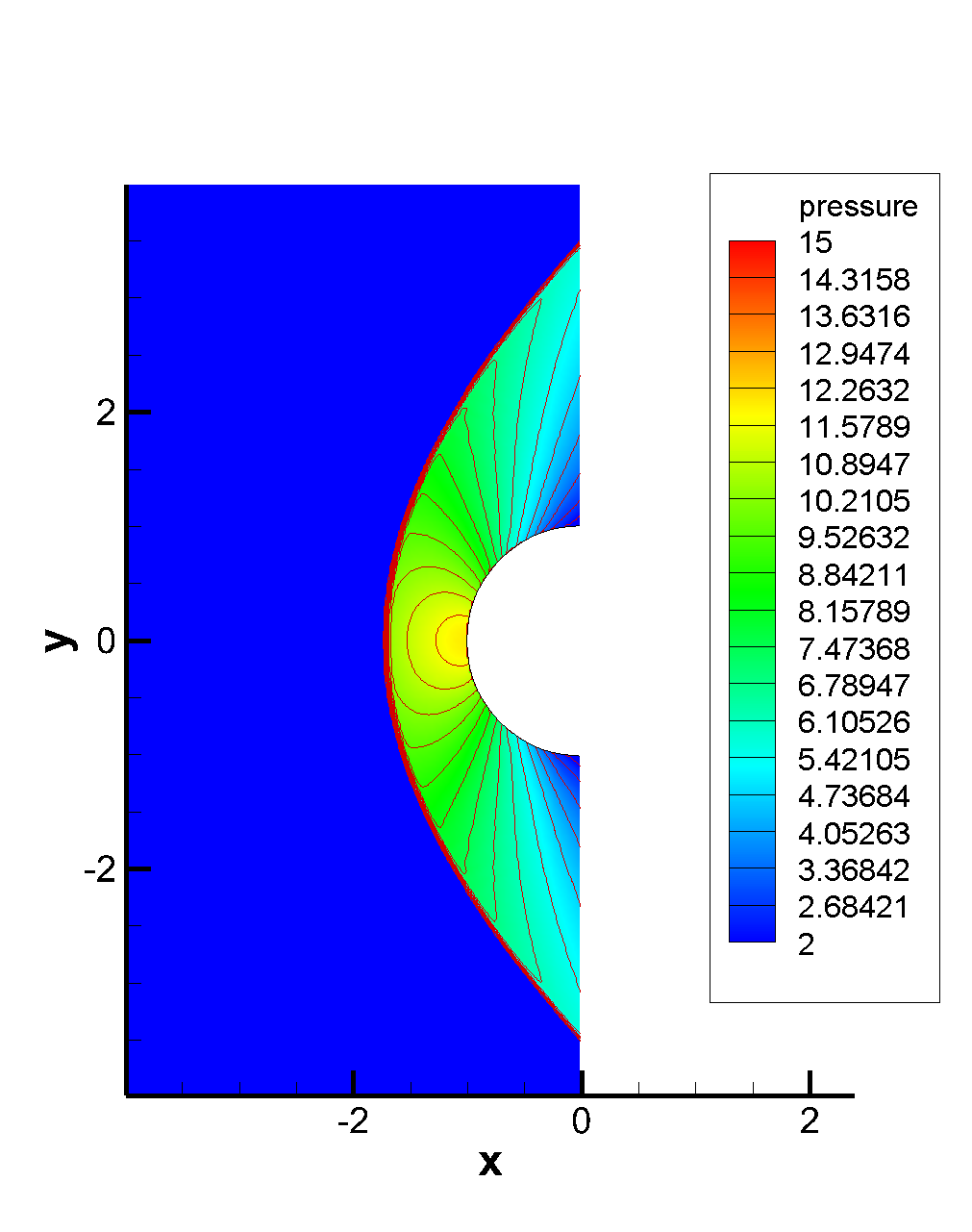}}
	\subfigure[Mach number contour]{
	\includegraphics[width=1.0in,angle=0,scale=2.7]{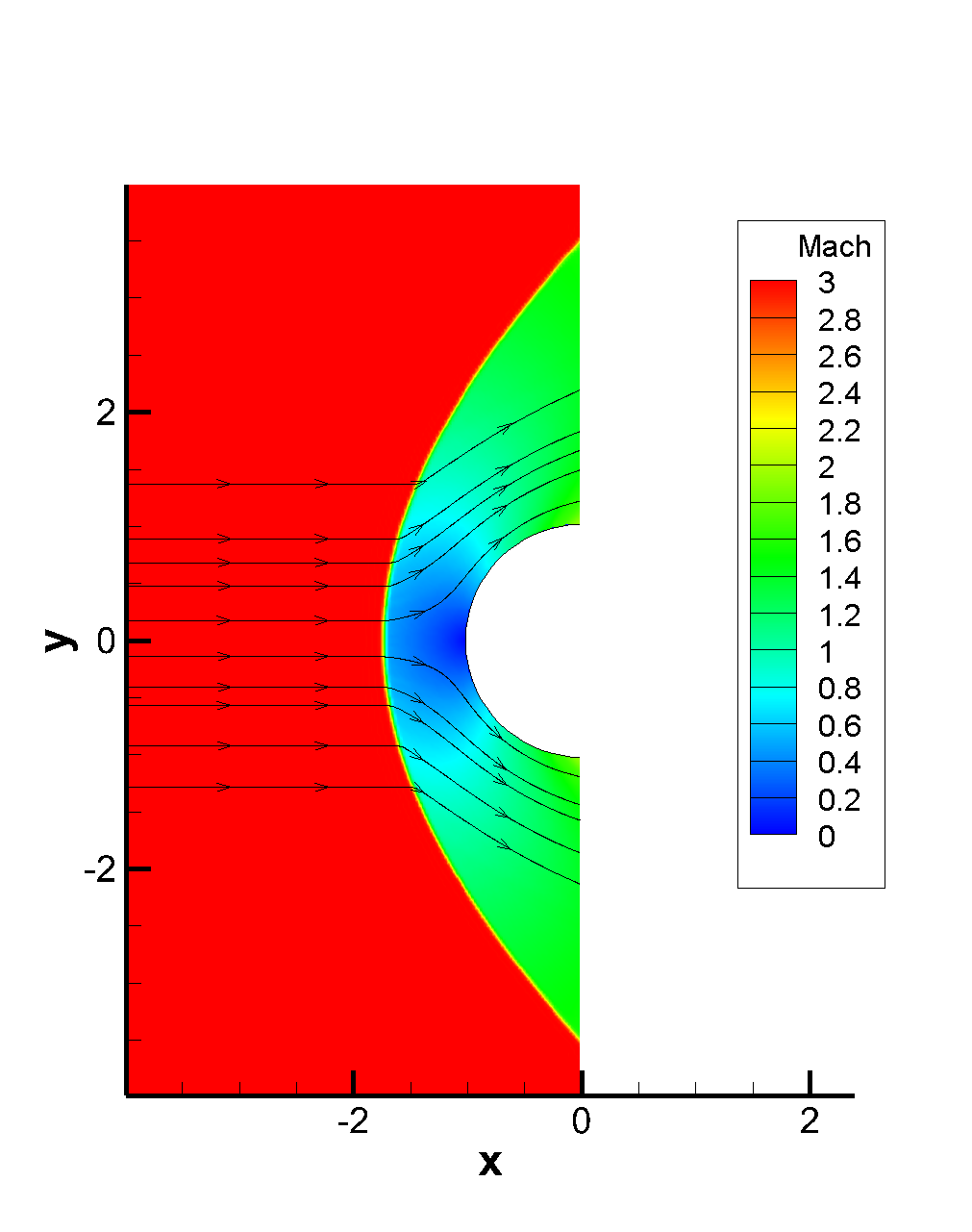}}
\caption{Example 7: Left: Pressure contour of flow past a cylinder. 
Right: Mach number contour with streamline. }
\label{eg7f2}
\end{figure}

\vspace{0.3cm}
\noindent
\textbf{Example 7.}
We consider a flow around a cylinder {\color{blue} \cite{TS1}}. The center of the cylinder is located at the origin, and the cylinder has a radius of 1. At the initial moment, a fluid with Mach 3 moves towards the cylinder. In consideration of
 the symmetry, we only consider the problem of an upper half plane. 
 For the lower boundary of the computation area at $y = 0$, we use the reflection technique; for the left boundary of the computation region at $x = - 4$, we give the inflow boundary condition. 
To show the full shape of the generated shock wave, we take the upper boundary $y = 6$ such that the shock would not reach the boundary.
Therefore, for the right boundary $x = 0$ and the upper boundary $y = 6$ of the computation area, the outflow boundary conditions are given. On the surface of a cylinder, our new ILW method is used to deal with a no-penetration boundary condition. 
As before, a uniform non body fitted Cartesian mesh with mesh size $\Delta x = \Delta y = h = 1/40$ is used, 
$$x_i = -4 + (i -1/2 )h,\quad y_j = (j - 1/2 )h, $$
which is shown in Figures \ref{grideg7}. 
Pressure field and Mach number field computed with fifth order scheme are plotted in Figure \ref{eg7f2}, which are comparable with those in \cite{TS1,TWSN,Lu2,Ding}.

%

\begin{figure}
	\subfigure[The computational region. ]{
		\includegraphics[width=0.5\textwidth]{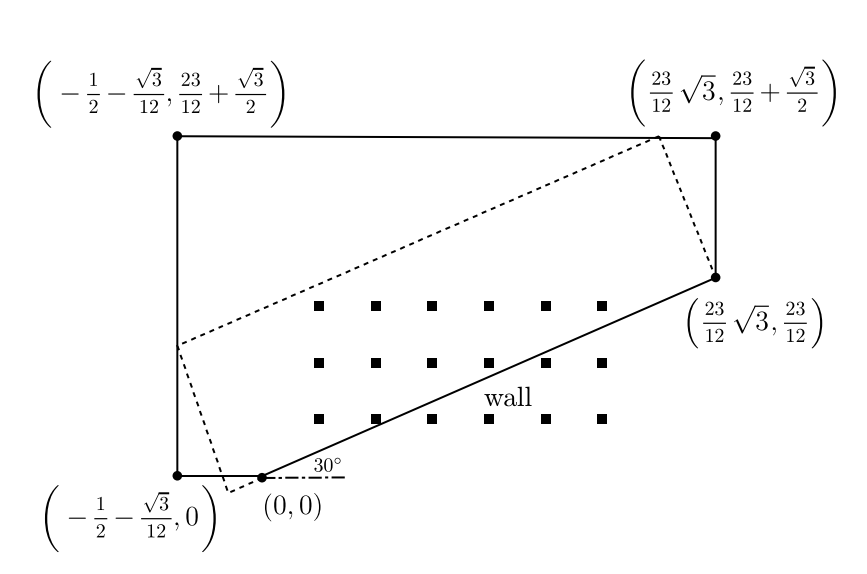}\label{db1.a}}
	\subfigure[Density contour. ]{
		\includegraphics[width=0.5\textwidth]{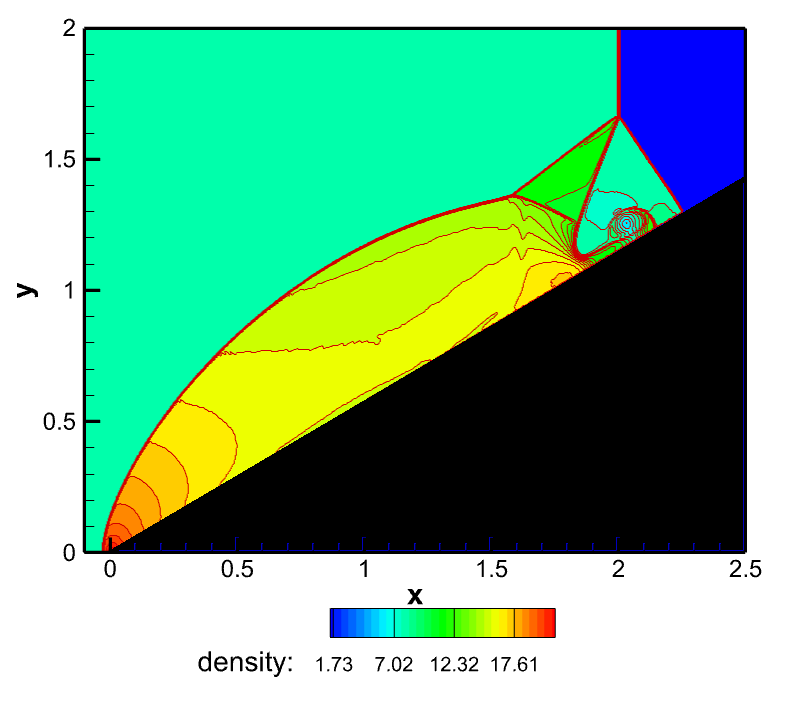}\label{db1.b}}
\caption{Example 8: Left: The computational region of the double Mach reflection problem. The dashed line indicates the computational domain used in \cite{Jiang, JYC}. Right: The density contour. 30 contours from 1.731 to 20.92. $\Delta x=\Delta y=1/320$.
}
\label{db1}
\end{figure}

\begin{figure}
	\subfigure[$\Delta x=\Delta y=1/320.$ ]{
		\includegraphics[width=0.5\textwidth]{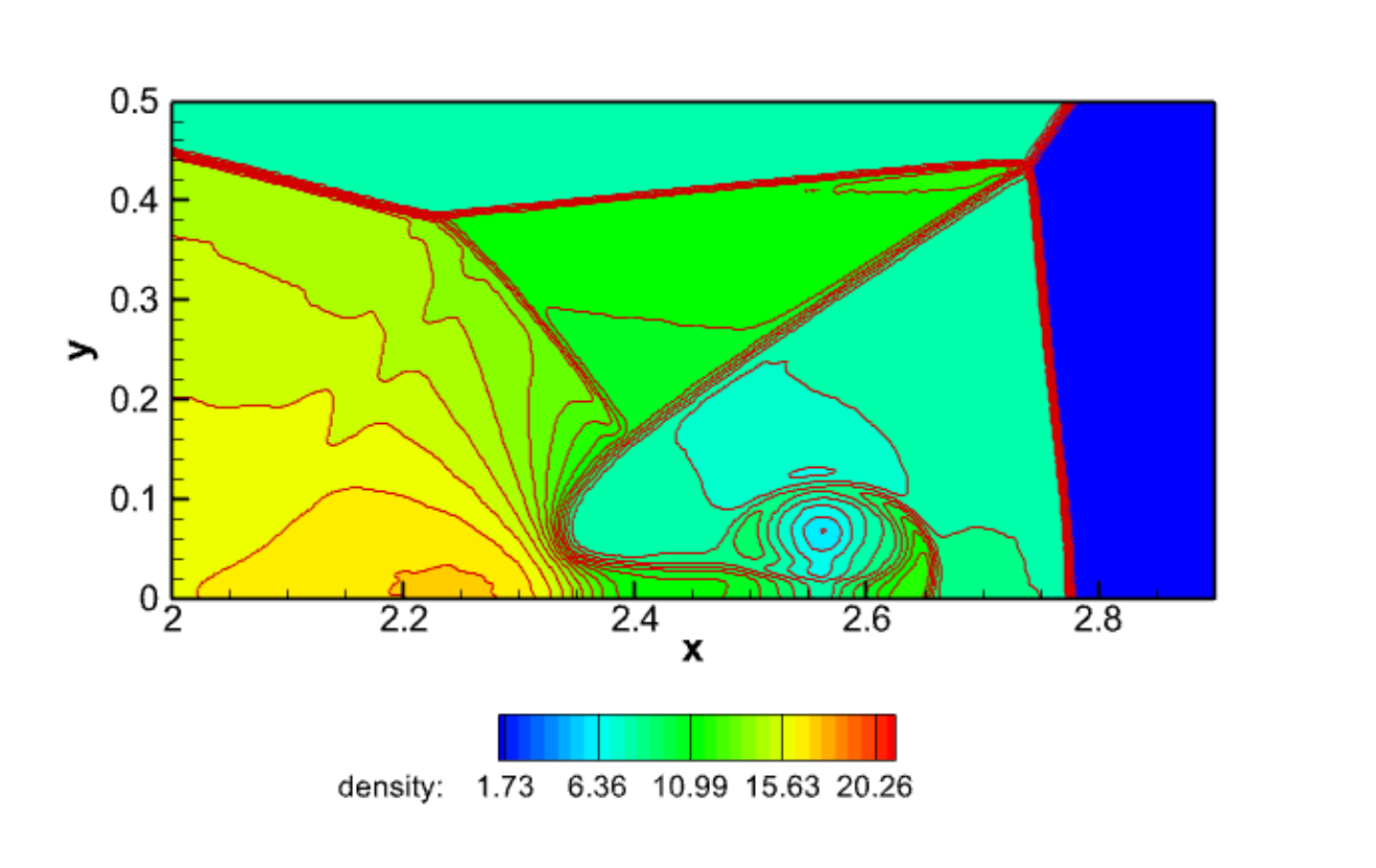}}
	\subfigure[$\Delta x=\Delta y=1/640.$ ]{
		\includegraphics[width=0.5\textwidth]{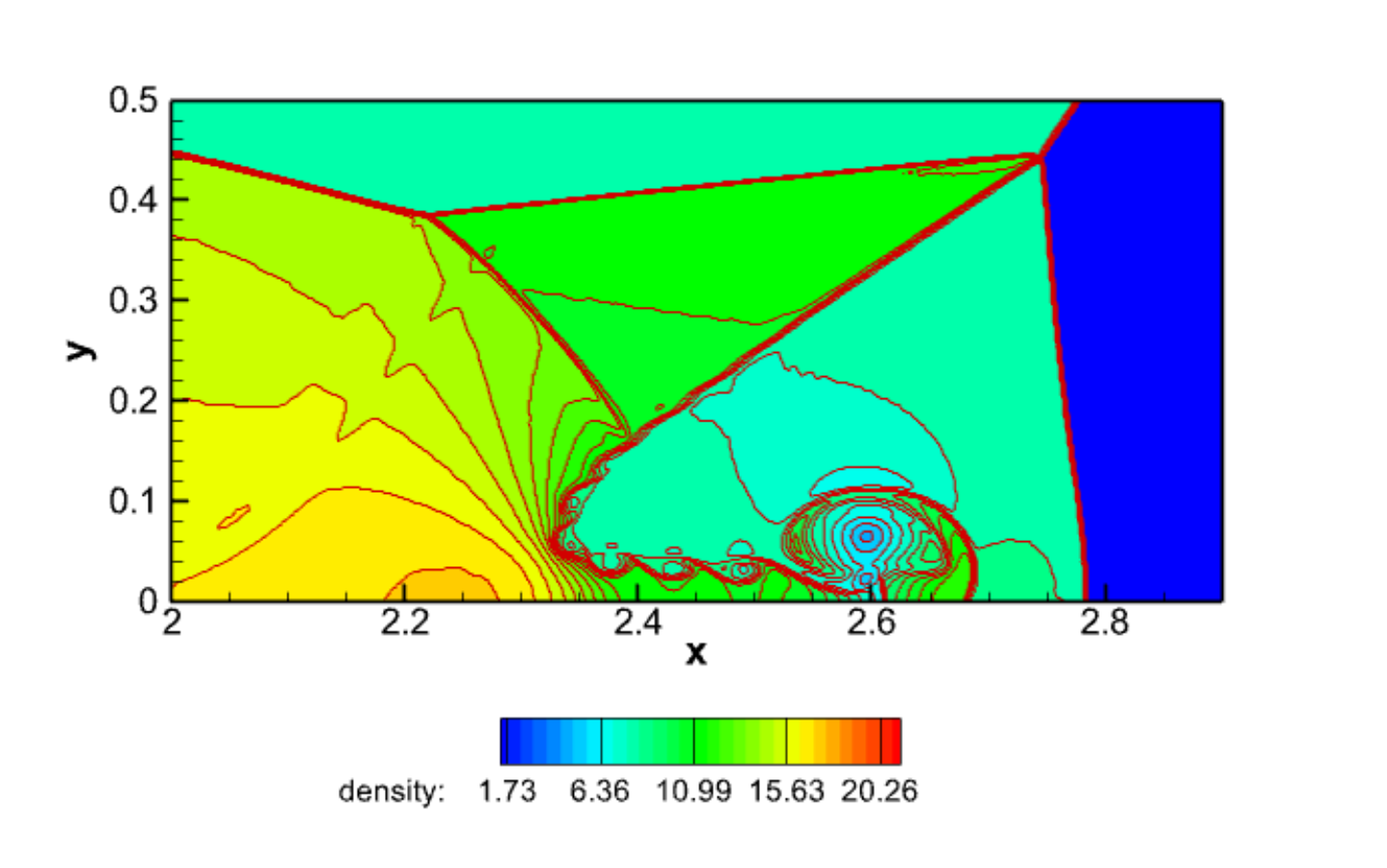}}
\caption{Example 8: Density contour on the local area. 30 contours from 1.731 to 20.92.}
\label{db2}
\end{figure}

\vspace{0.3cm}
\noindent
\textbf{Example 8.}
We consider the double Mach reflection problem with $\gamma=1.4$. At the initial moment, a horizontally moving Mach 10 shock wave passes through a wedge with an inclination angle of $30^{\circ}$. In common practice, the wedge is placed horizontally to apply reflective boundary conditions. Initially, the shock wave positioned at $x = 0$ forms an angle of $60^{\circ}$ with the wall. 
In \cite{Jiang, JYC}, the original double Mach number reflection problem is computed respectively. With the ILW method, people can also do numerical simulation on the original region \cite{TS1,TWSN}. Here, we use the new ILW method to simulate this problem. 
In detail, the computational region is the same as that in \cite{Jiang, JYC}, showed in Figure \ref{db1}(a).
At the top of the calculation area, we give the exact flow value according to the shock Mach number. At the left and right boundary, we give the supersonic inlet and outlet boundary conditions respectively. On the lower right boundary, the new ILW method is adopted. The discretization of space and time is consistent with the previous example. Figure \ref{db1}(b) shows the 
density contours at the time $t_{end} = 0.2$. The zoomed-in region near the double Mach stem is presented in Figure \ref{db2}. We rotate and translate the region for ease of comparison. It is observed that the new ILW method captures the shock wave well, and it is comparable with the previous results.

\section{Concluding remarks}
In this paper, we propose a new SILW method for conservation laws, which decomposes the procedure of the construction of
the ghost point values into two steps: interpolation and extrapolation.
First, we approximate some artificial auxiliary point values through a polynomial interpolating the interior points near boundary. Then, we construct a Hermite extrapolation polynomial based on those auxiliary point values and spatial derivatives at the boundary obtained through the ILW procedure. 
After that, we can get the approximation of the ghost point values. 
Through the linear stability analysis with the eigenvalue method, we can conclude that our new SILW method requires fewer terms using the ILW procedure, 
thus is more efficient than the original SILW method on the premise of ensuring the stability. 
This improvement is more noticeable for higher order schemes.
We then extend our new SILW method to systems and two-dimensional cases, and carry out a series of numerical experiments. The numerical results domenstrate that our new SILW method is stable and can achieve the expected accuracy. 
In the future, we are going to extend this new SILW method to deal with the initial-boundary value problems of diffusion equations and convection-diffusion equations.

\begin{appendices}
\section{More results about linear stability analysis} \label{sec:append}
\setcounter{table}{0}
\setcounter{figure}{0}

\begin{figure}[htb!]
\subfigure[Fifth order scheme with $\alpha=0.91$]{
		\includegraphics[width=0.5\textwidth]{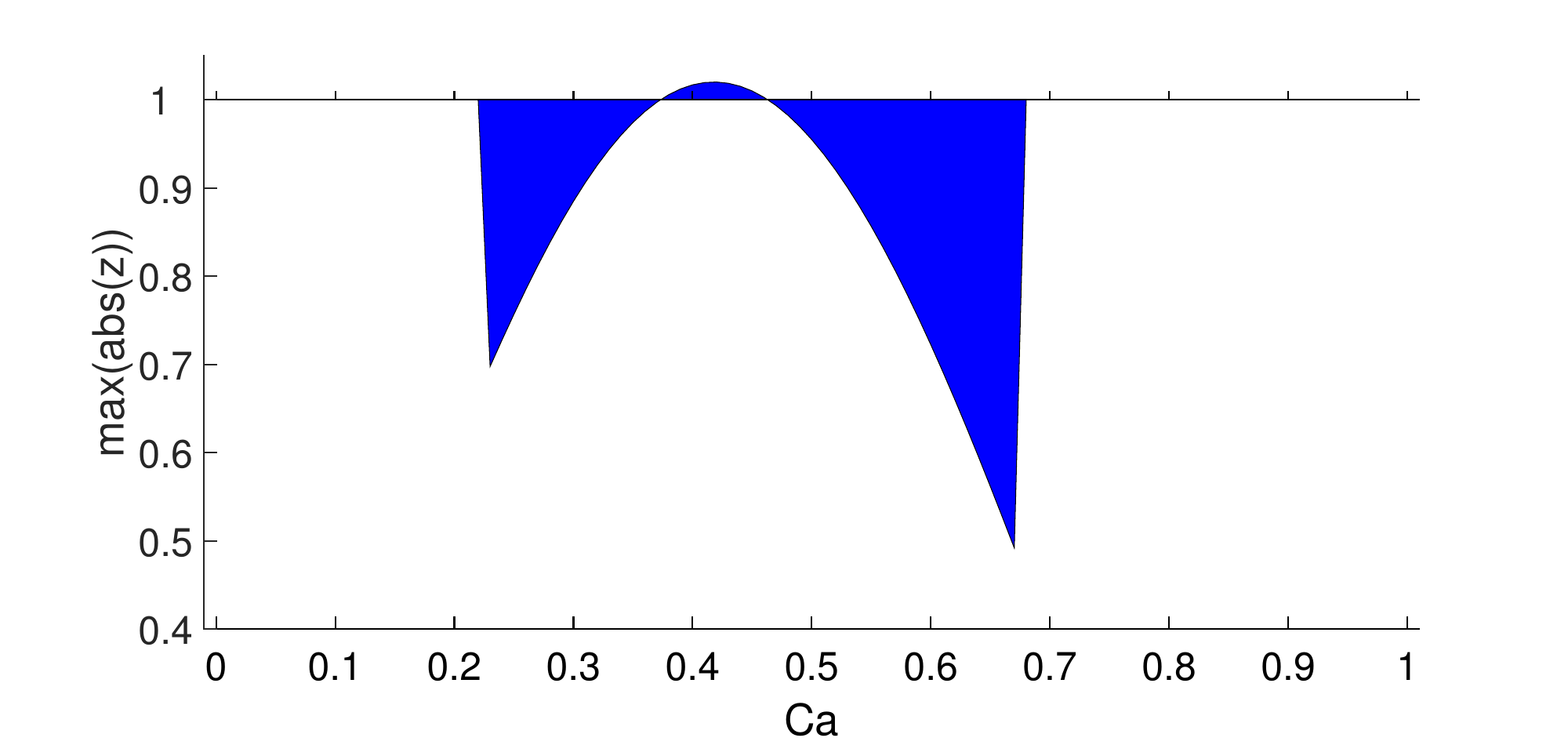}}
\subfigure[Fifth order scheme with $\alpha=0.92$]{
		\includegraphics[width=0.5\textwidth]{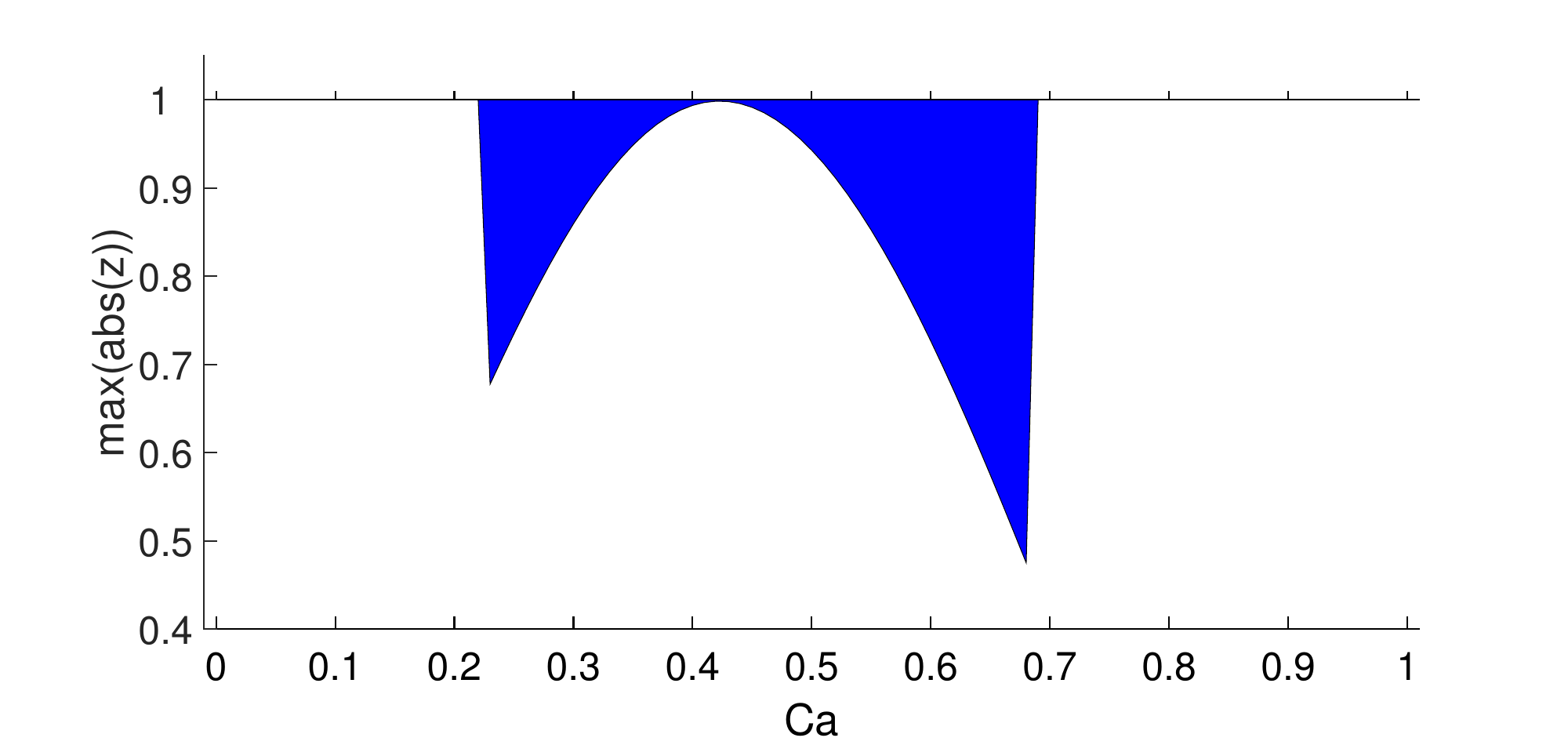}}
\subfigure[Fifth order scheme with $\alpha=5.11$]{
		\includegraphics[width=0.5\textwidth]{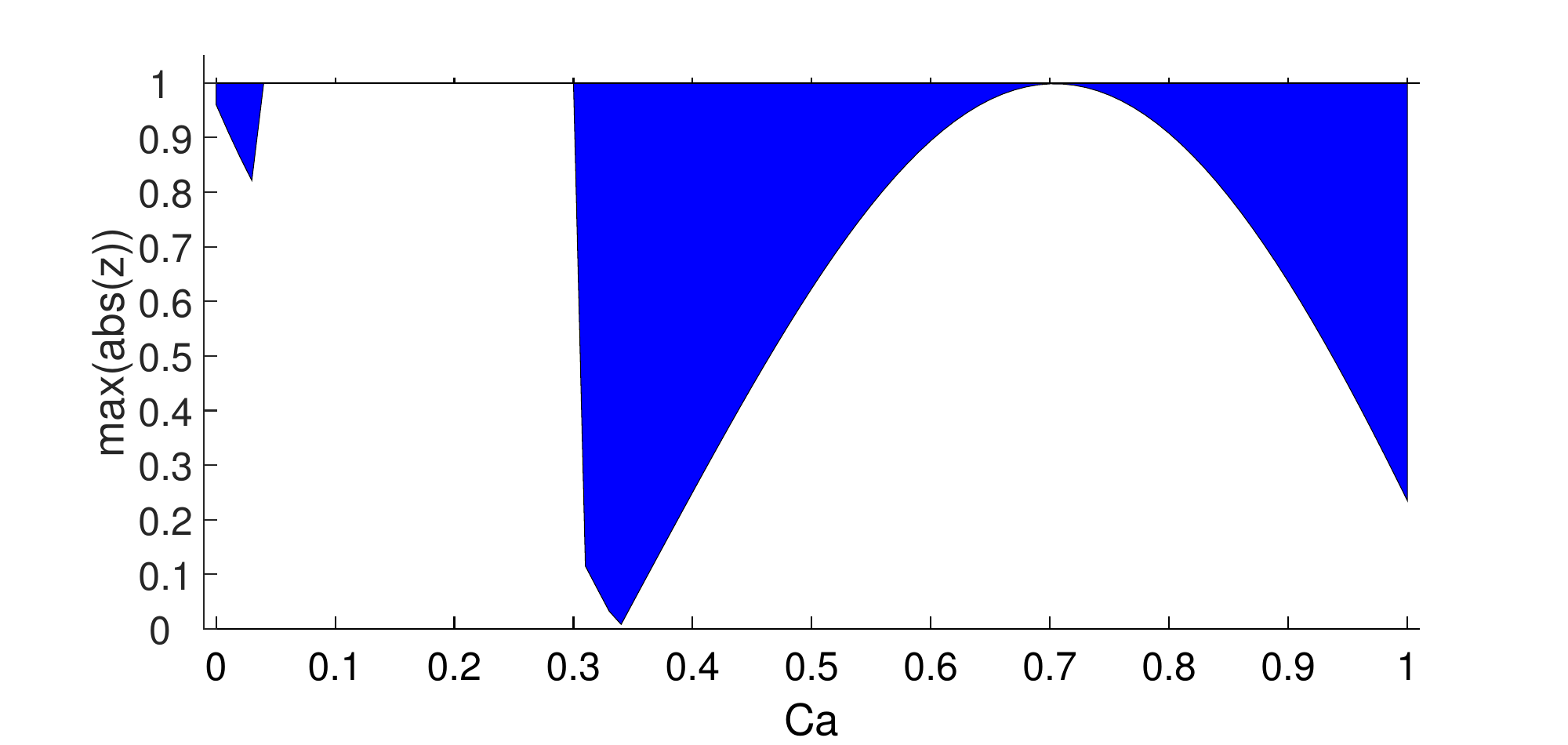}}
\subfigure[Fifth order scheme with $\alpha=5.12$]{
		\includegraphics[width=0.5\textwidth]{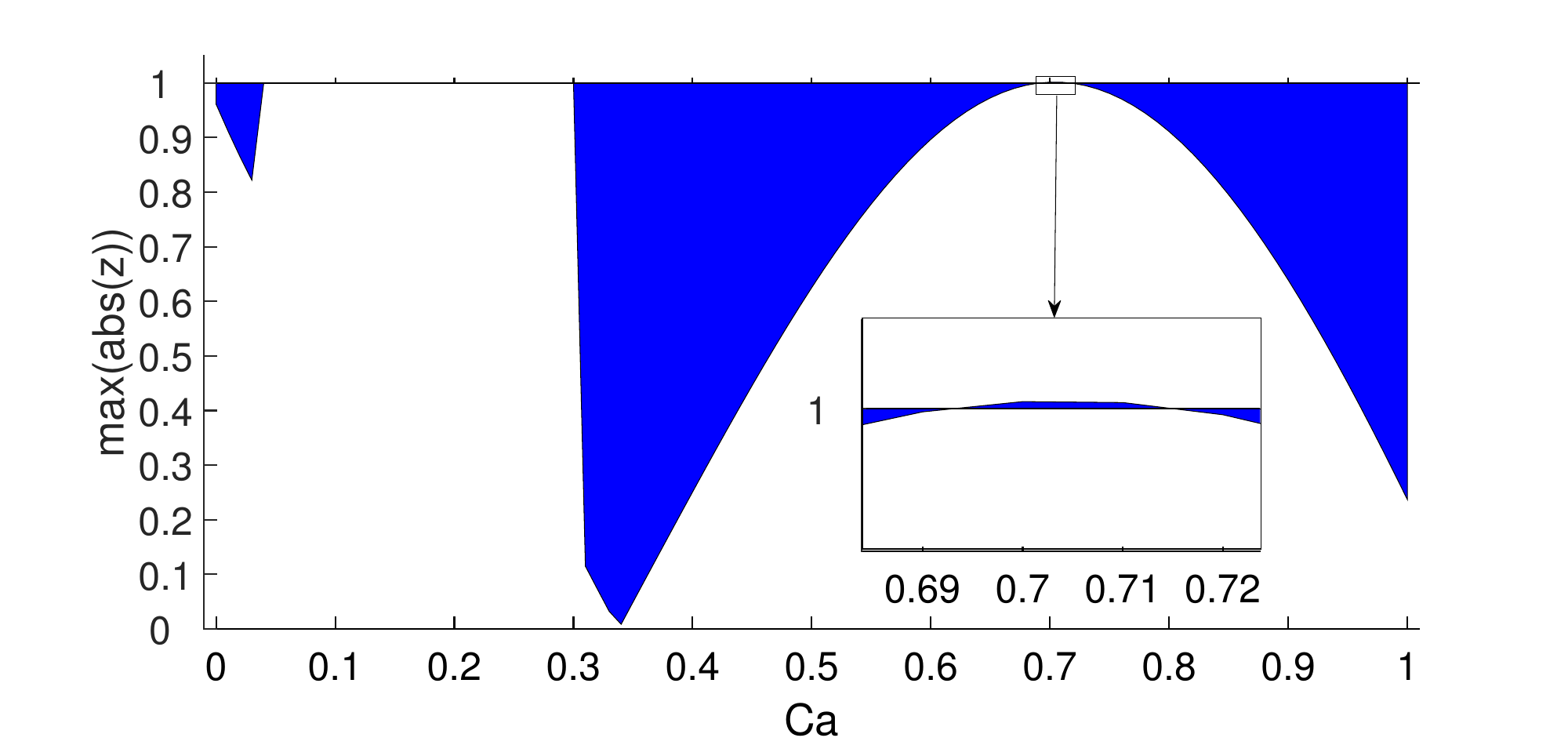}}
\subfigure[Seventh order scheme with $\alpha=1.33$]{
		\includegraphics[width=0.5\textwidth]{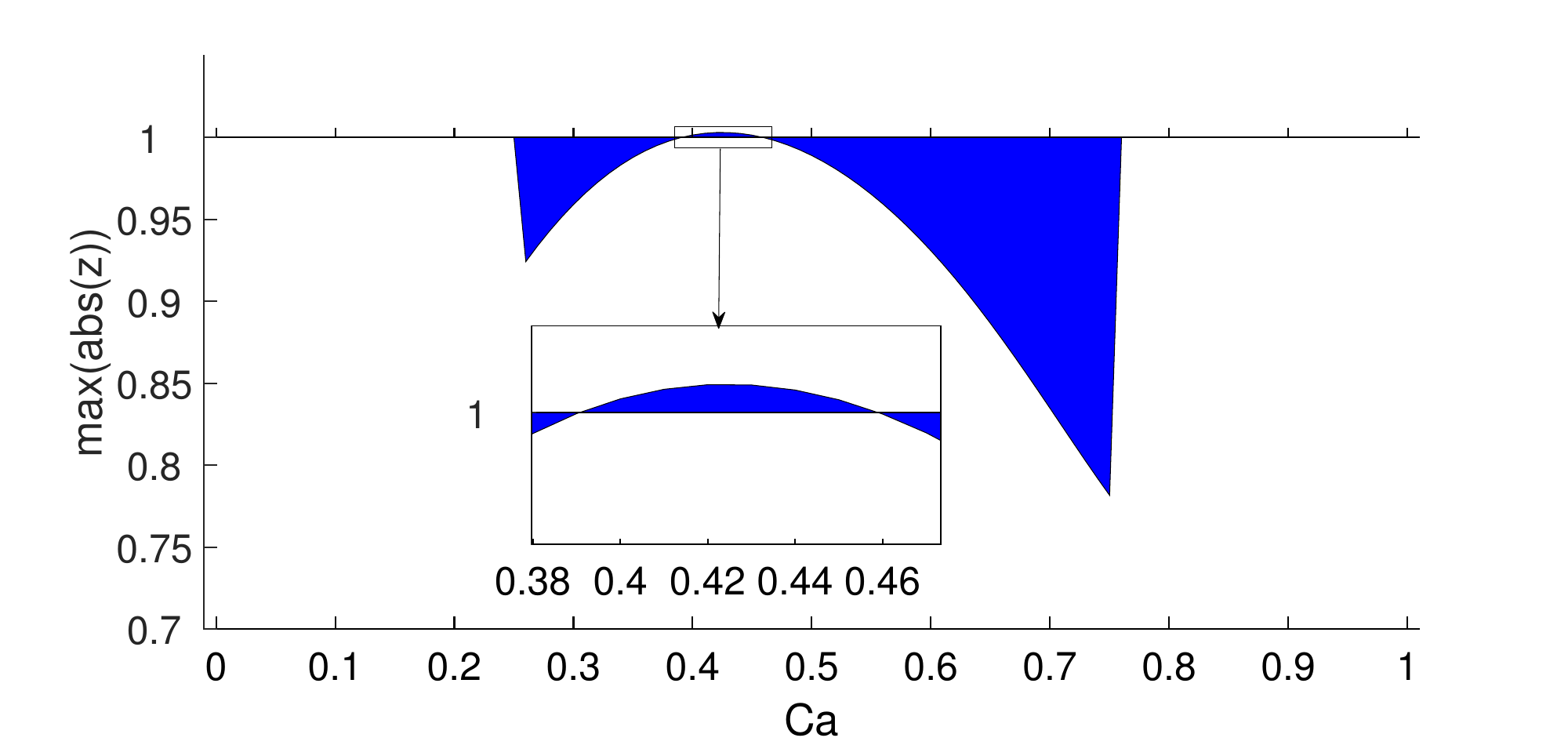}}	
\subfigure[Seventh order scheme with $\alpha=1.34$]{
		\includegraphics[width=0.5\textwidth]{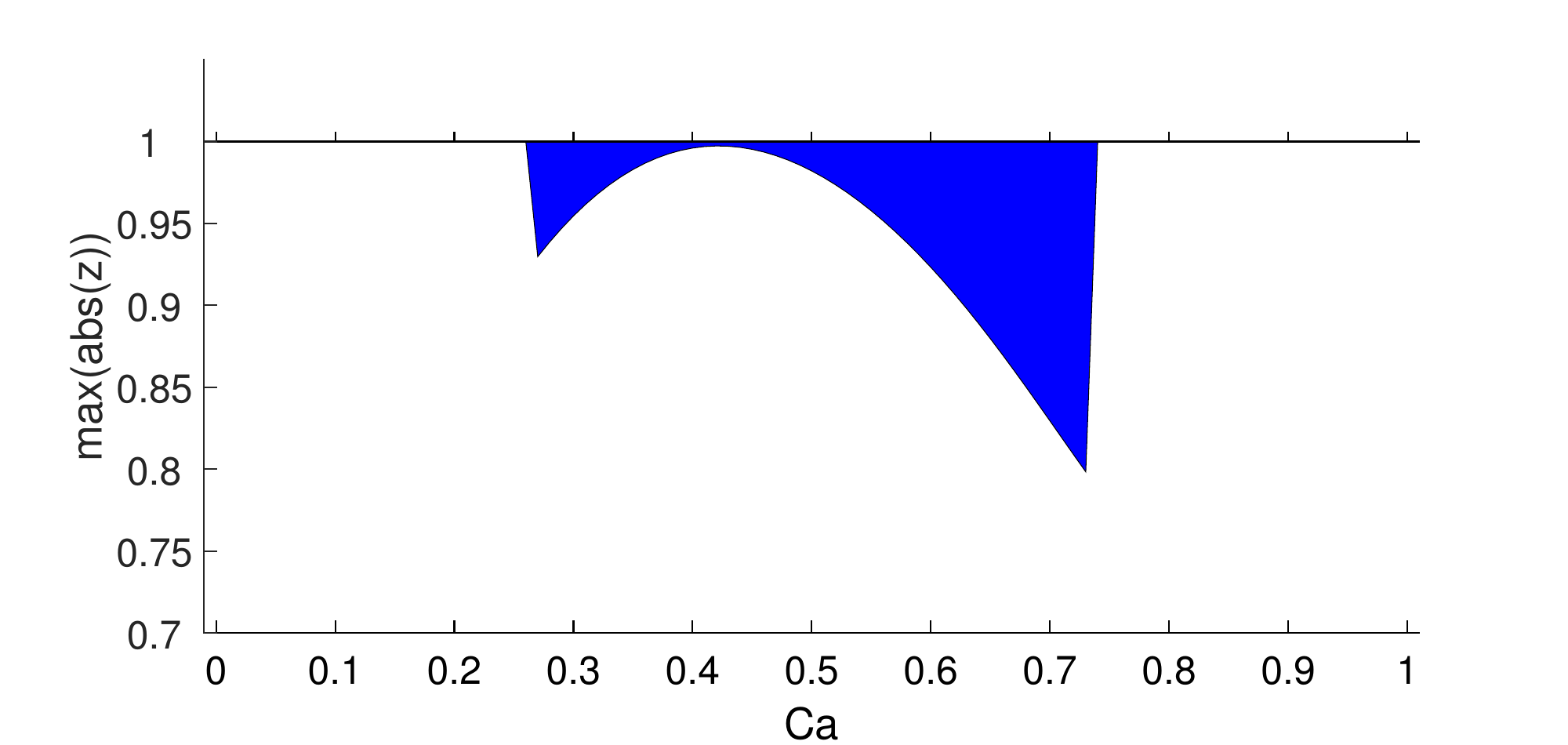}}
\subfigure[Seventh order scheme with $\alpha=1.99$]{
		\includegraphics[width=0.5\textwidth]{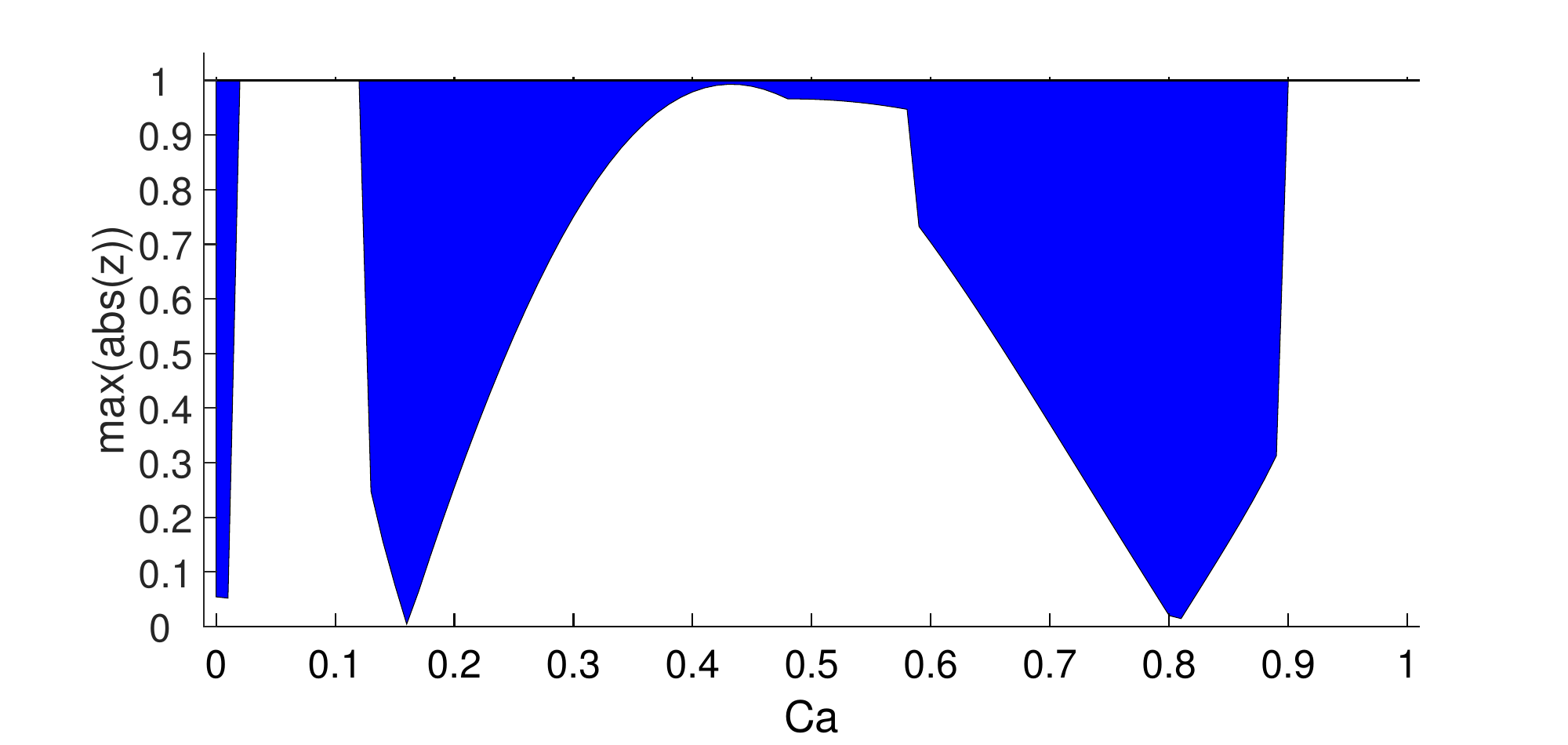}}
\subfigure[Seventh order scheme with $\alpha=2.00$]{
		\includegraphics[width=0.5\textwidth]{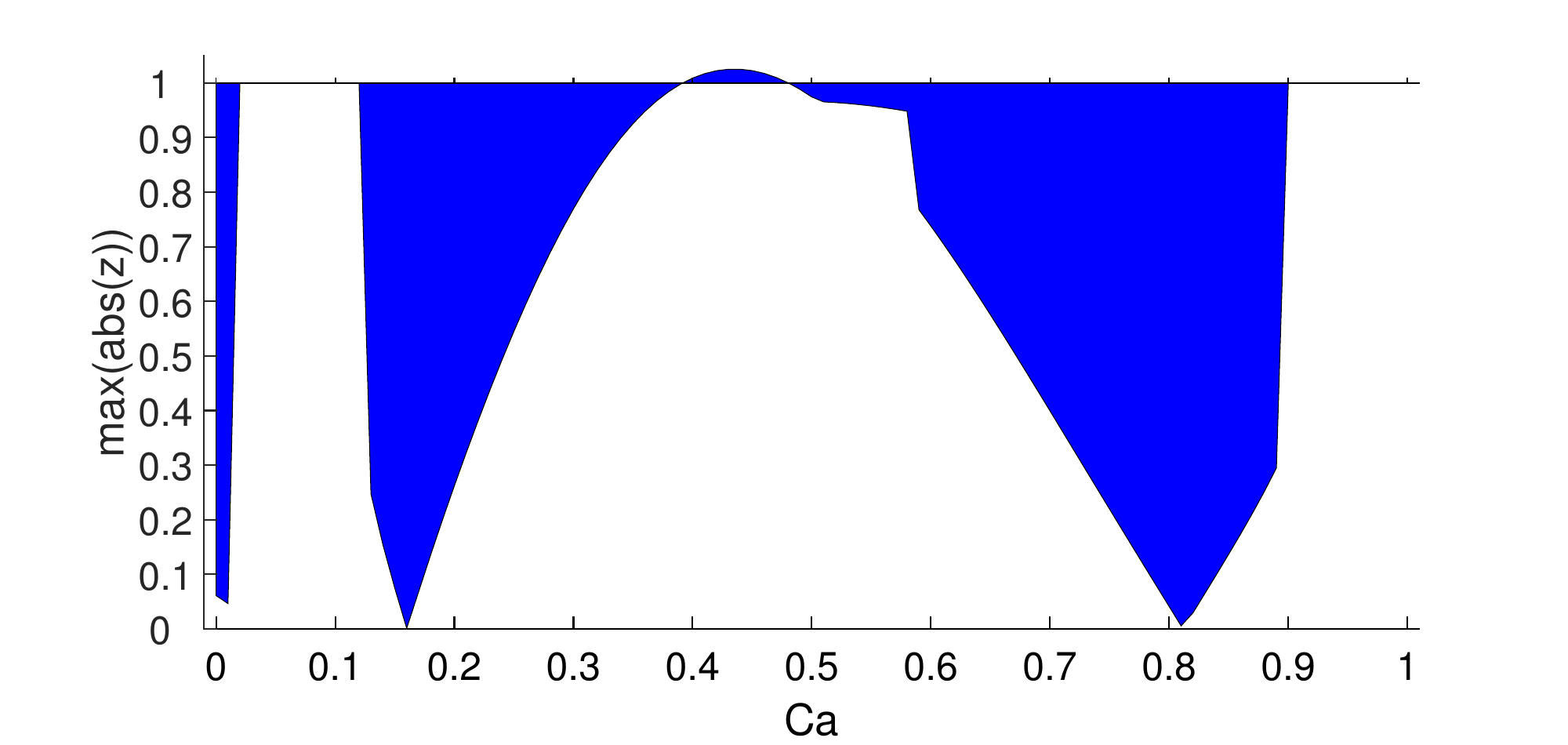}}
\caption{The result of linear stability analysis with $k_d=2$. The horizontal axis represents $C_a$ and the vertical axis represents the largest $|z(\mu)|$.
}
\label{fig:stabilitykd2}
\end{figure}

\begin{figure}[htb!]
\subfigure[Nineth order scheme with $\alpha=1.28$]{
		\includegraphics[width=0.5\textwidth]{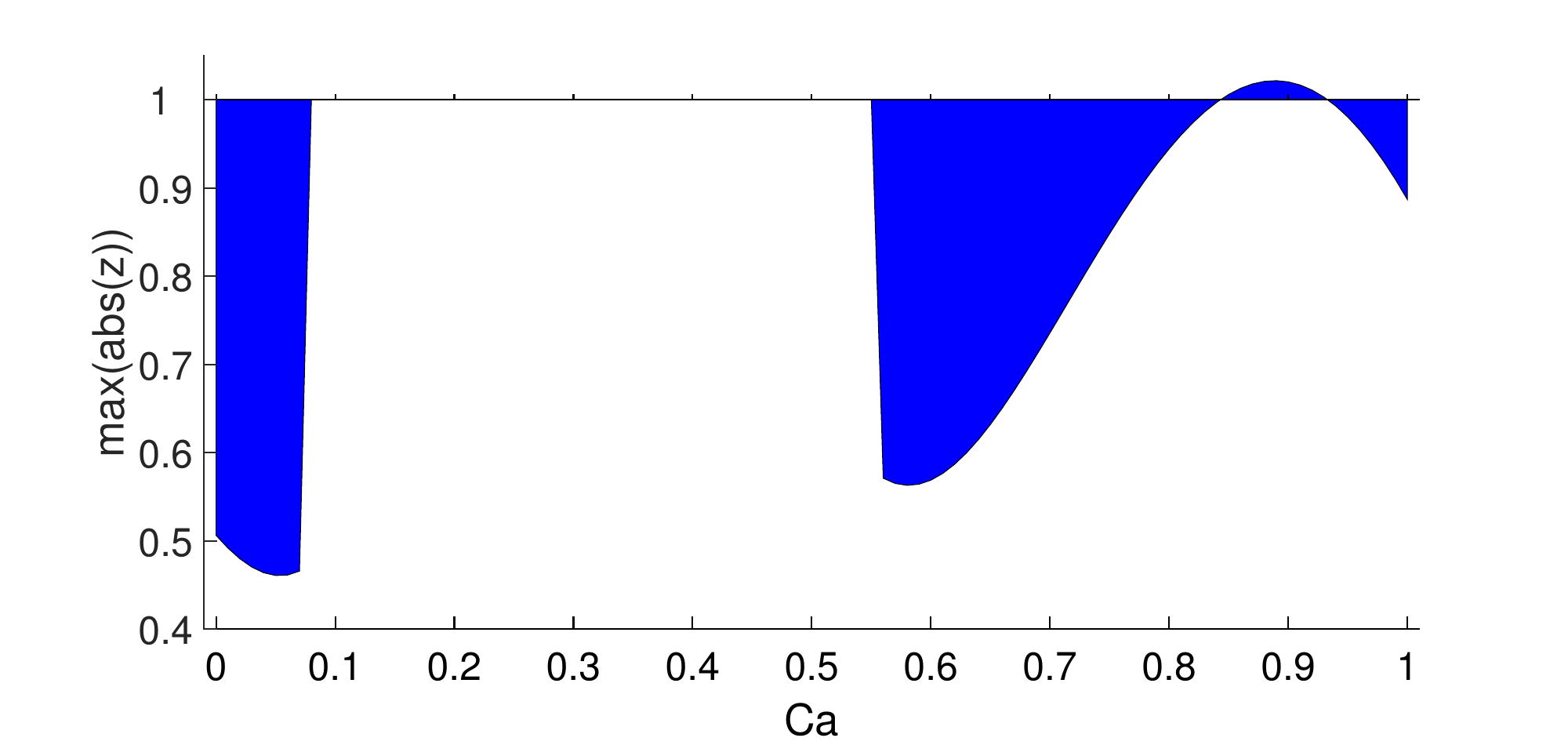}}
\subfigure[Nineth order scheme with $\alpha=1.29$]{
		\includegraphics[width=0.5\textwidth]{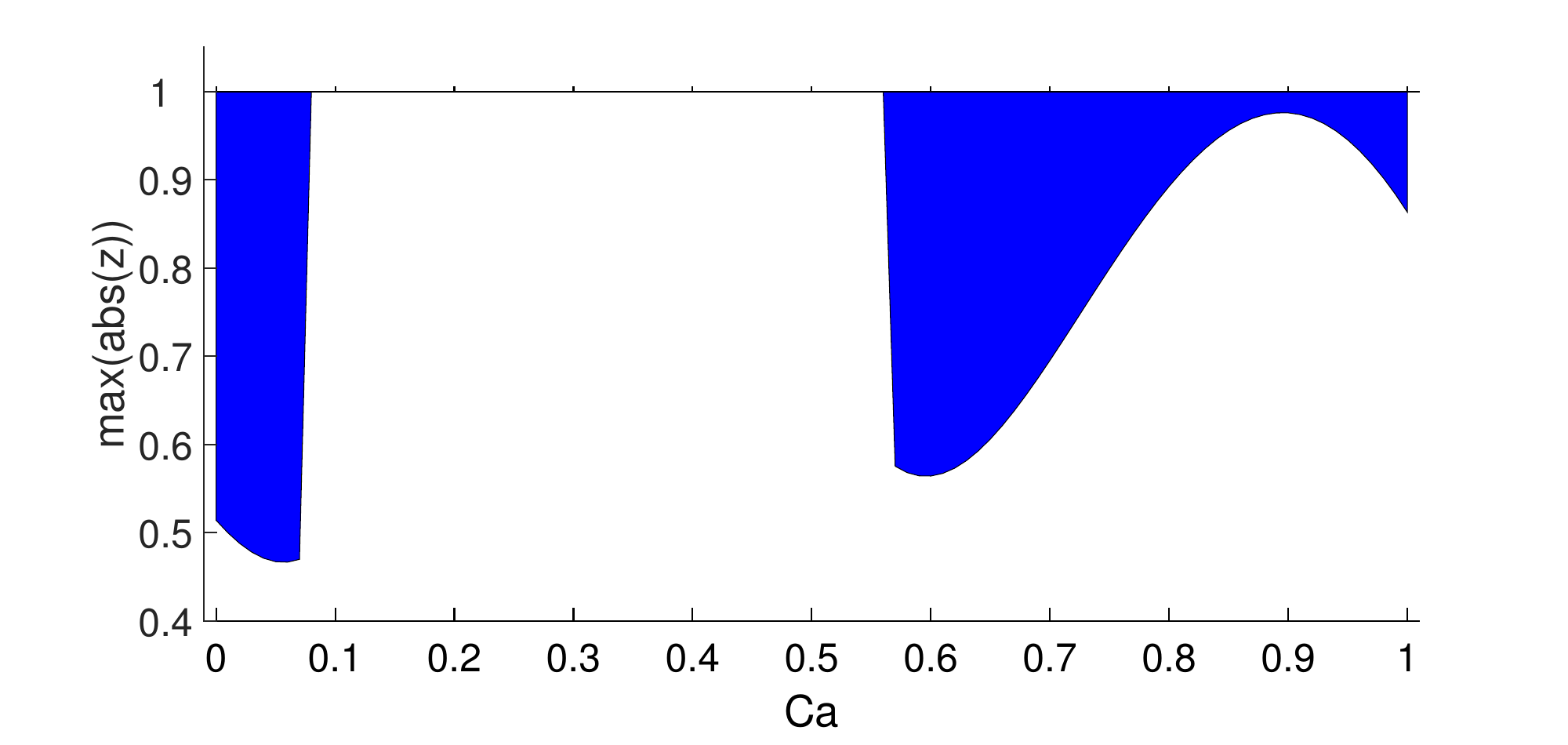}}
\subfigure[Nineth order scheme with $\alpha=2.43$]{
		\includegraphics[width=0.5\textwidth]{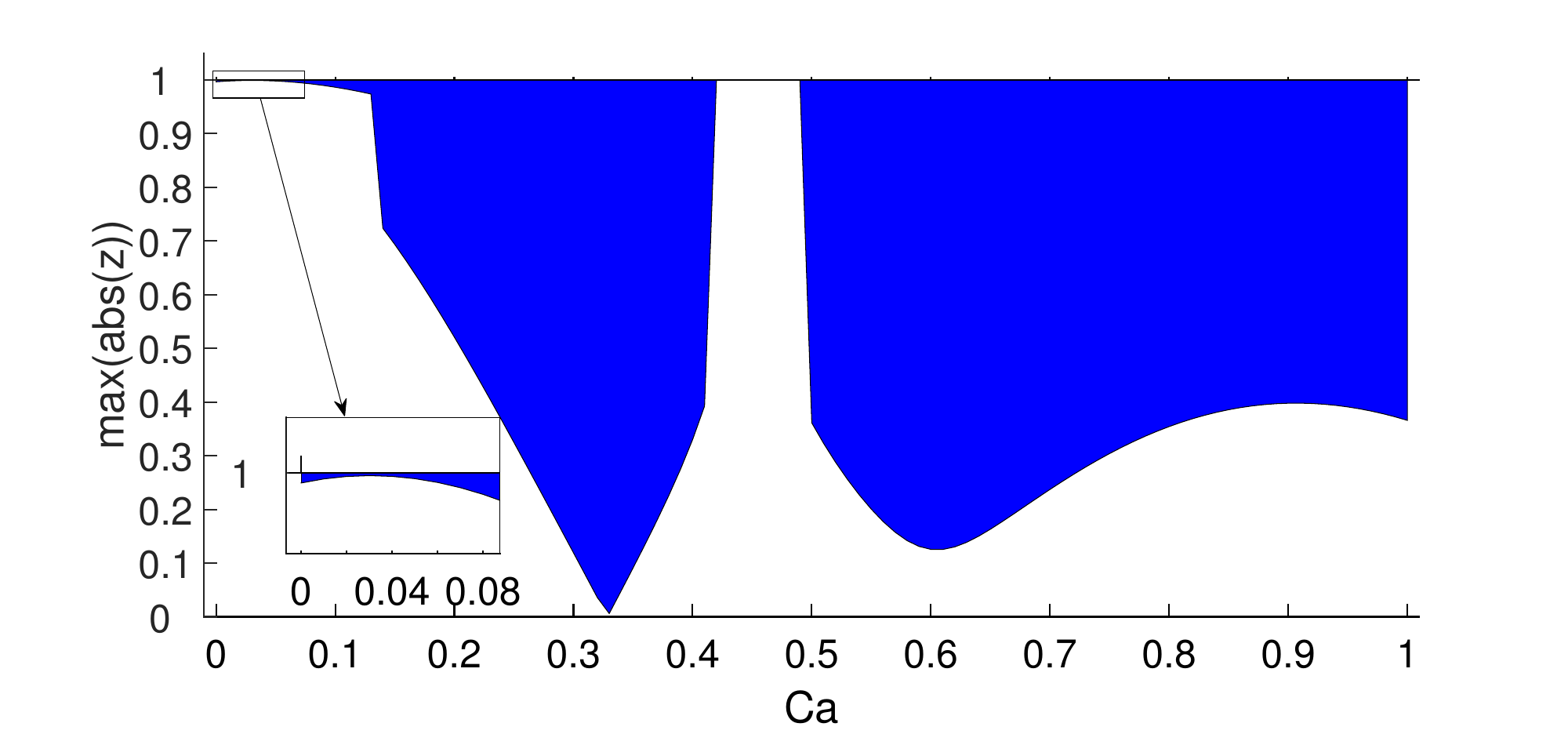}}
\subfigure[Nineth order scheme with $\alpha=2.44$]{
		\includegraphics[width=0.5\textwidth]{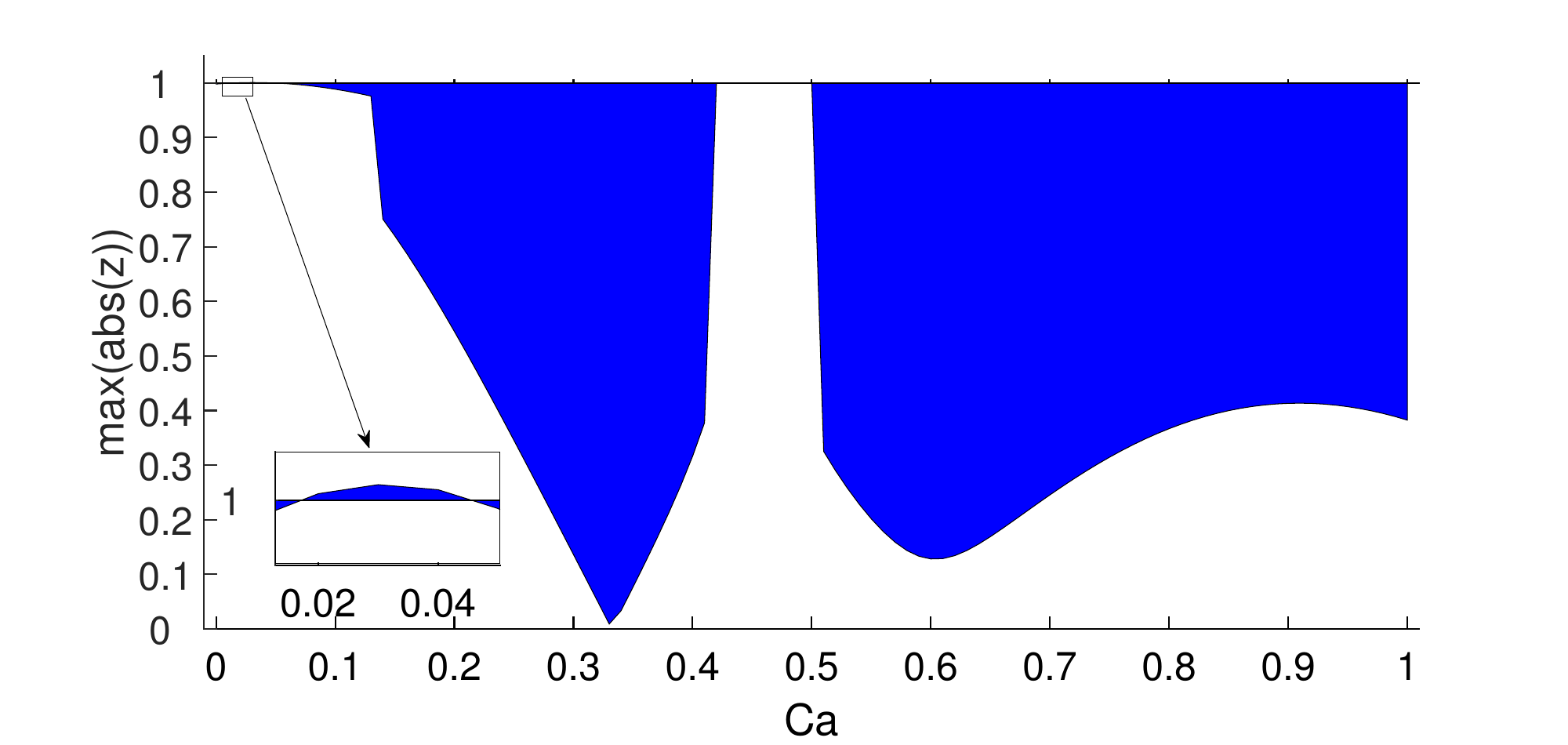}}		
\subfigure[Eleventh order scheme with $\alpha=1.41$]{
		\includegraphics[width=0.5\textwidth]{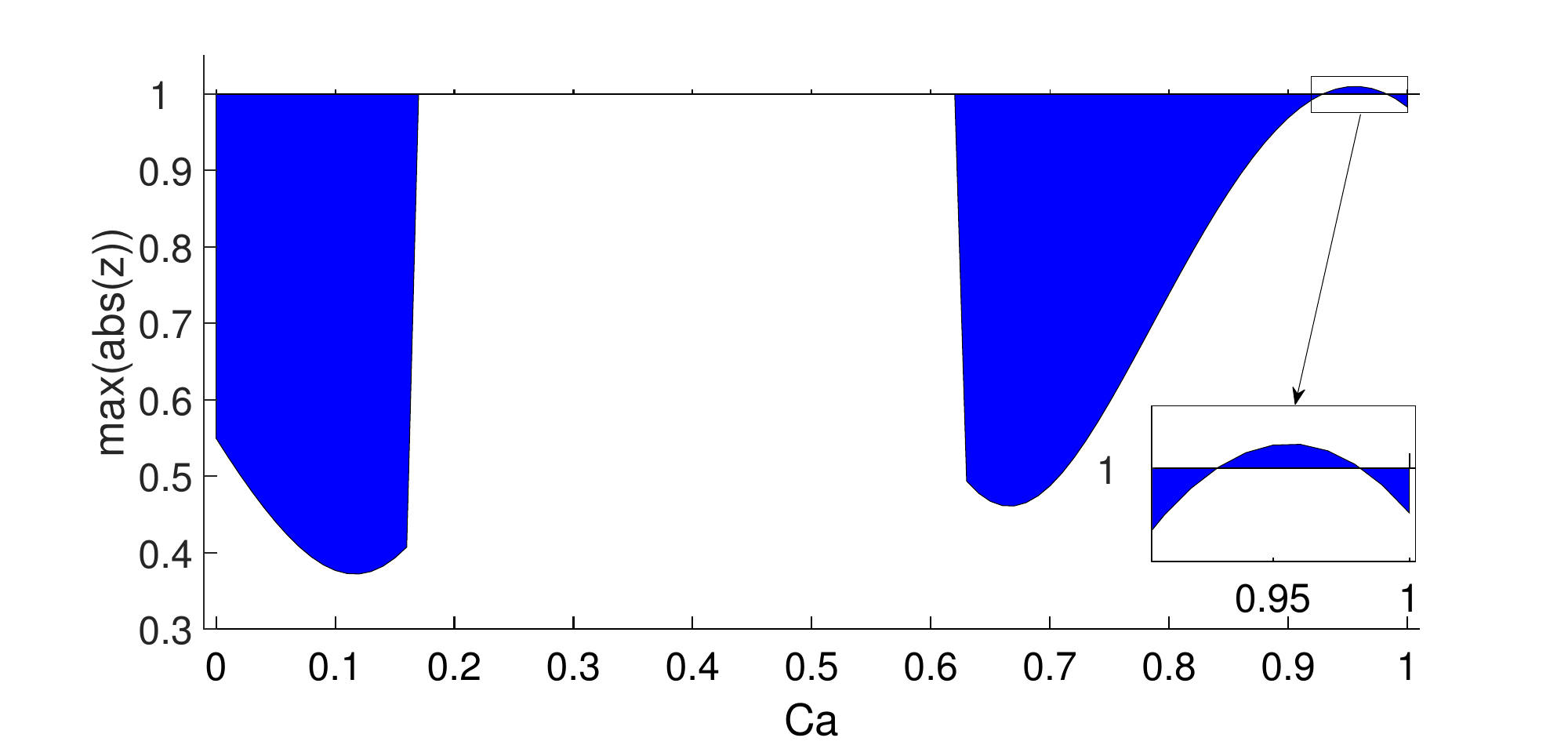}}	
\subfigure[Eleventh order scheme with $\alpha=1.42$]{
		\includegraphics[width=0.5\textwidth]{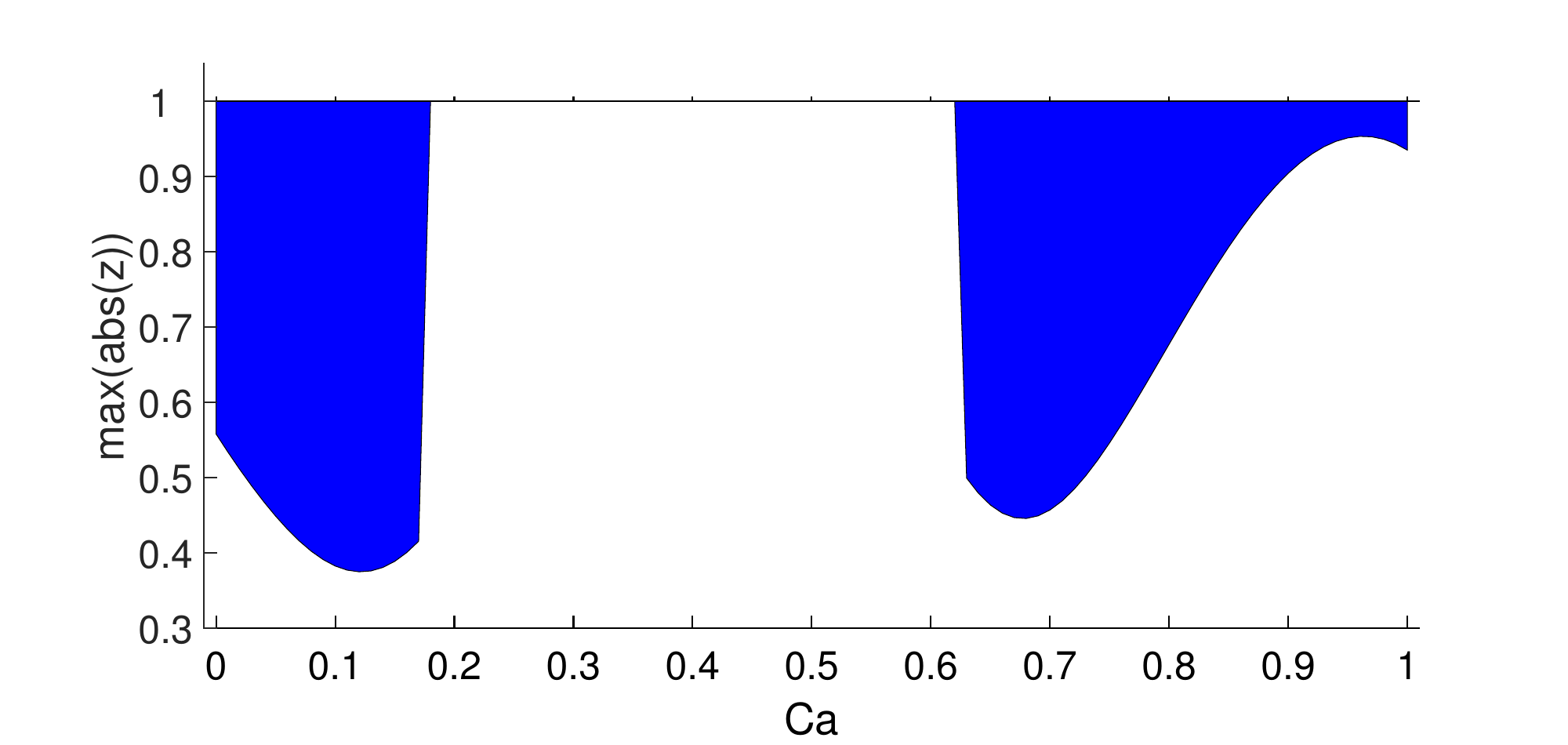}}
\subfigure[Eleventh order scheme with $\alpha=1.70$]{
		\includegraphics[width=0.5\textwidth]{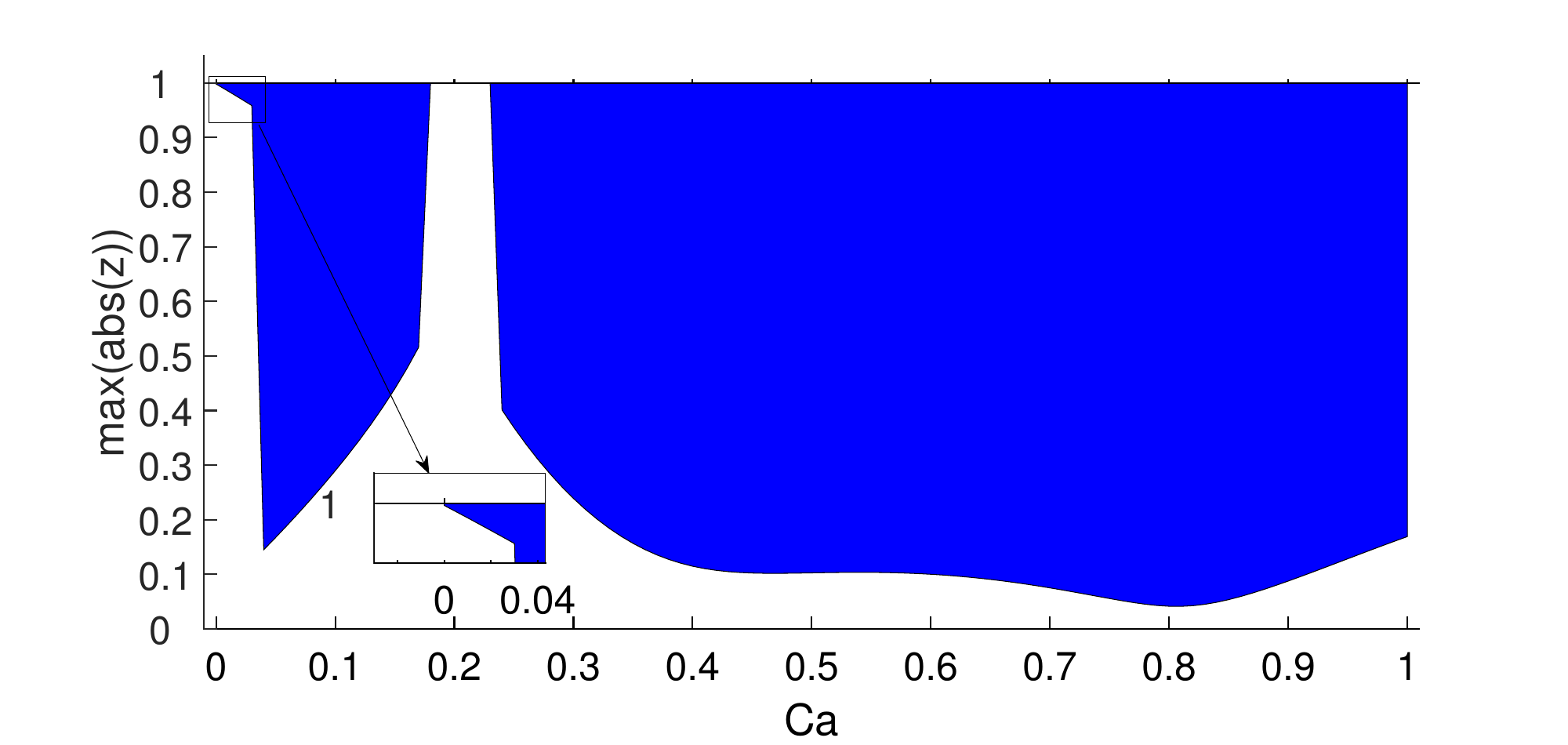}}
\subfigure[Eleventh order scheme with $\alpha=1.71$]{
		\includegraphics[width=0.5\textwidth]{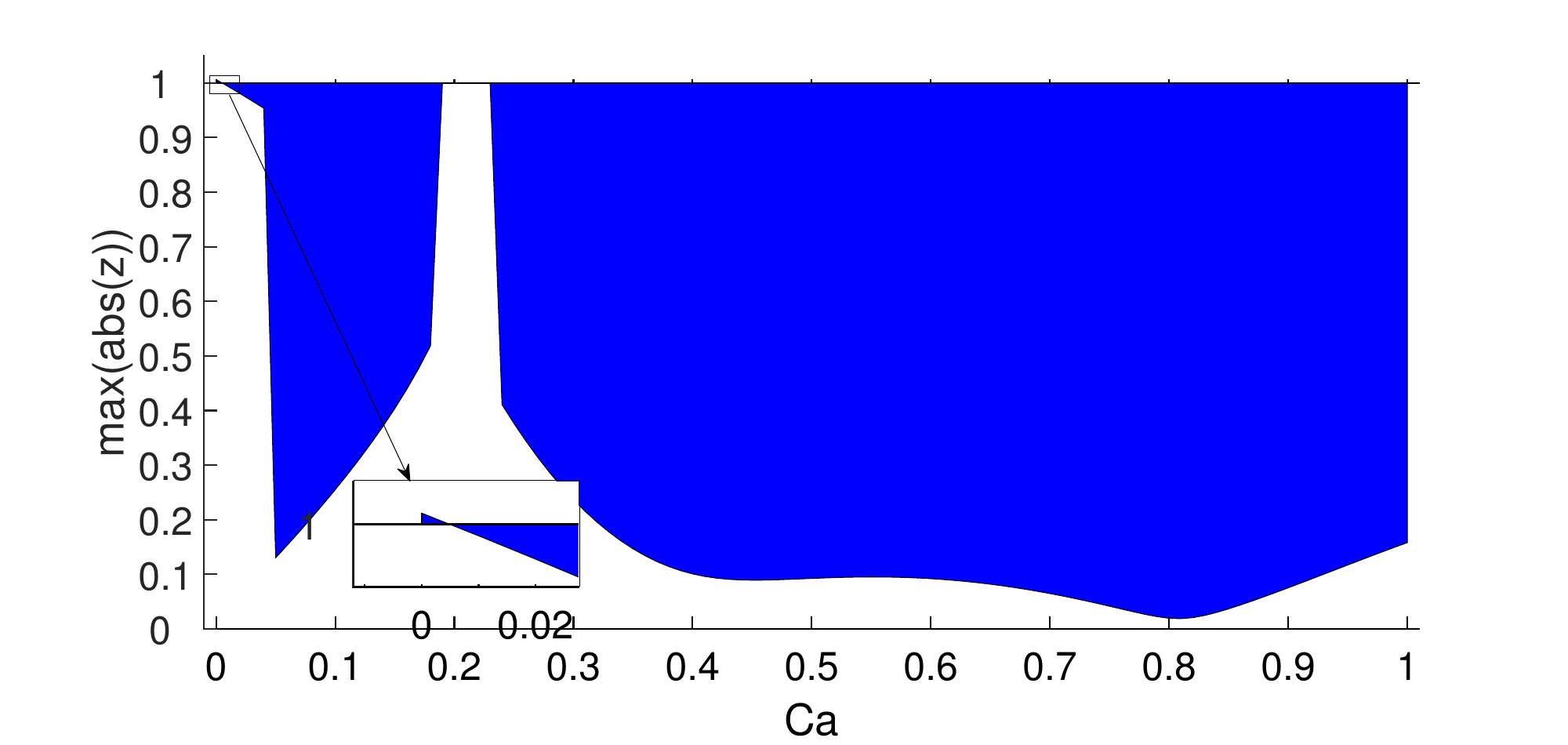}}
\caption{The result of linear stability analysis with $k_d=3$. The horizontal axis represents $C_a$ and the vertical axis represents the largest $|z(\mu)|$.
}
\label{fig:stabilitykd3}
\end{figure}

\begin{figure}[htb!]
\subfigure[Thirteenth order scheme with $\alpha=1.48$]{
		\includegraphics[width=0.5\textwidth]{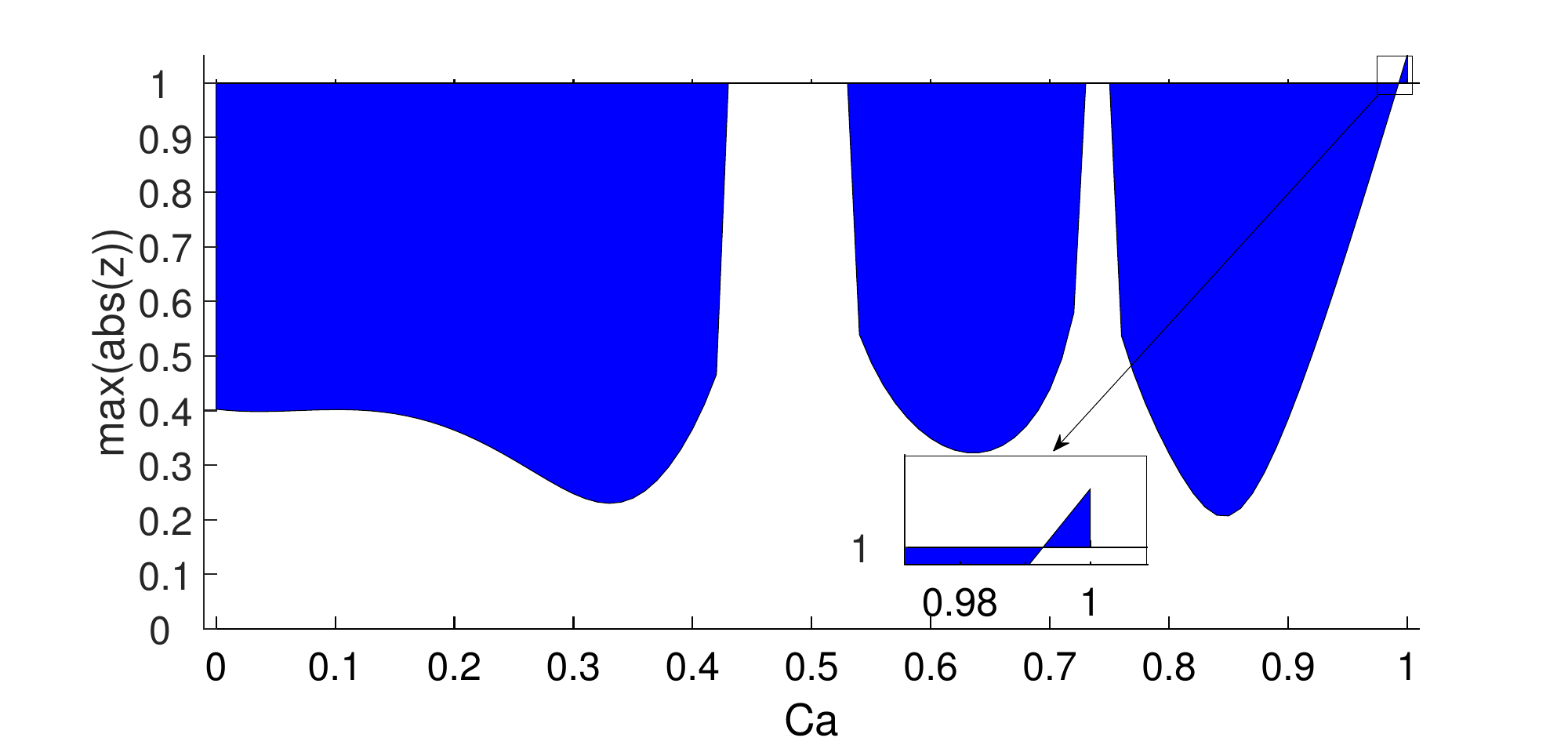}}
\subfigure[Thirteenth order scheme with $\alpha=1.49$]{
		\includegraphics[width=0.5\textwidth]{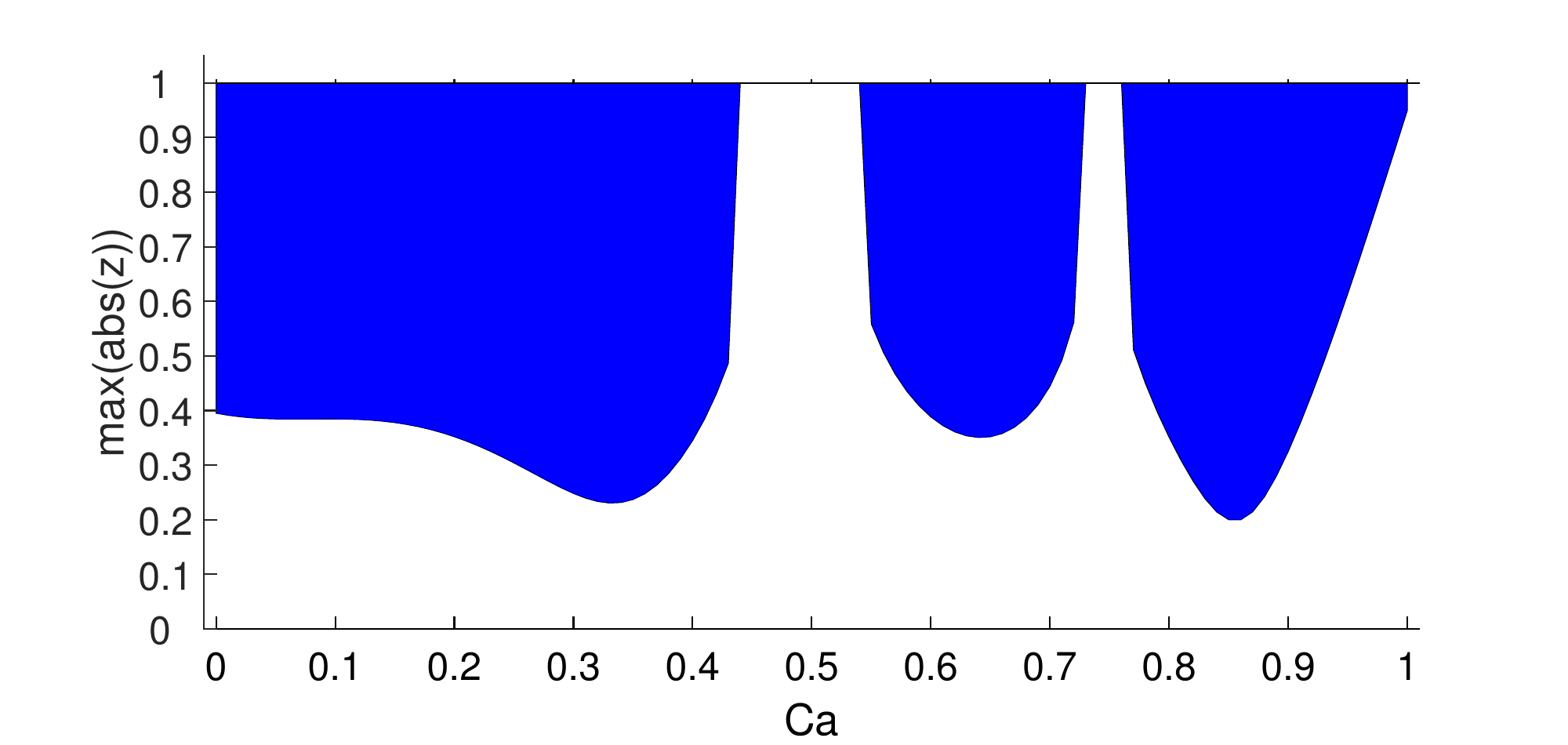}}
\subfigure[Thirteenth order scheme with $\alpha=2.08$]{
		\includegraphics[width=0.5\textwidth]{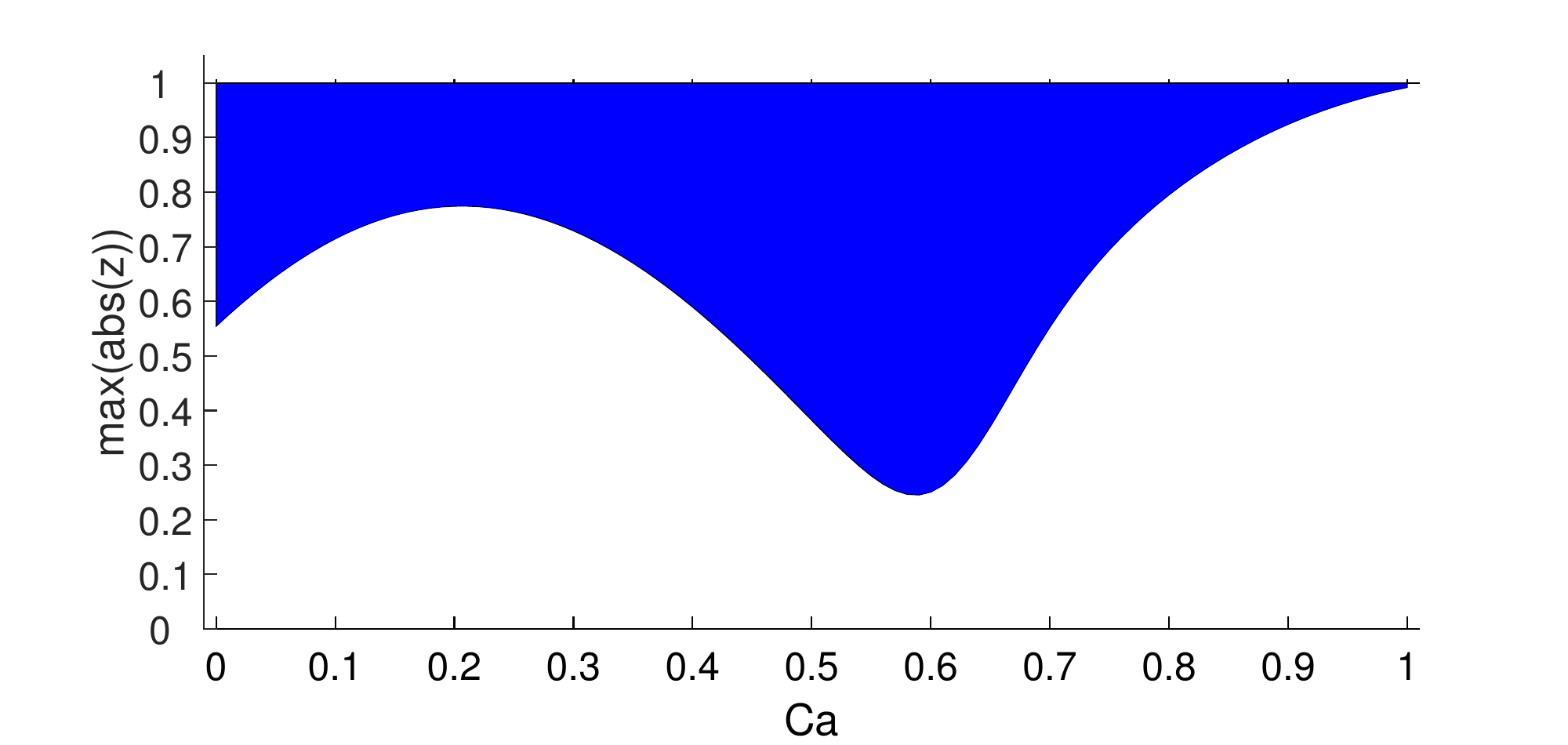}}
\subfigure[Thirteenth order scheme with $\alpha=2.09$]{
		\includegraphics[width=0.5\textwidth]{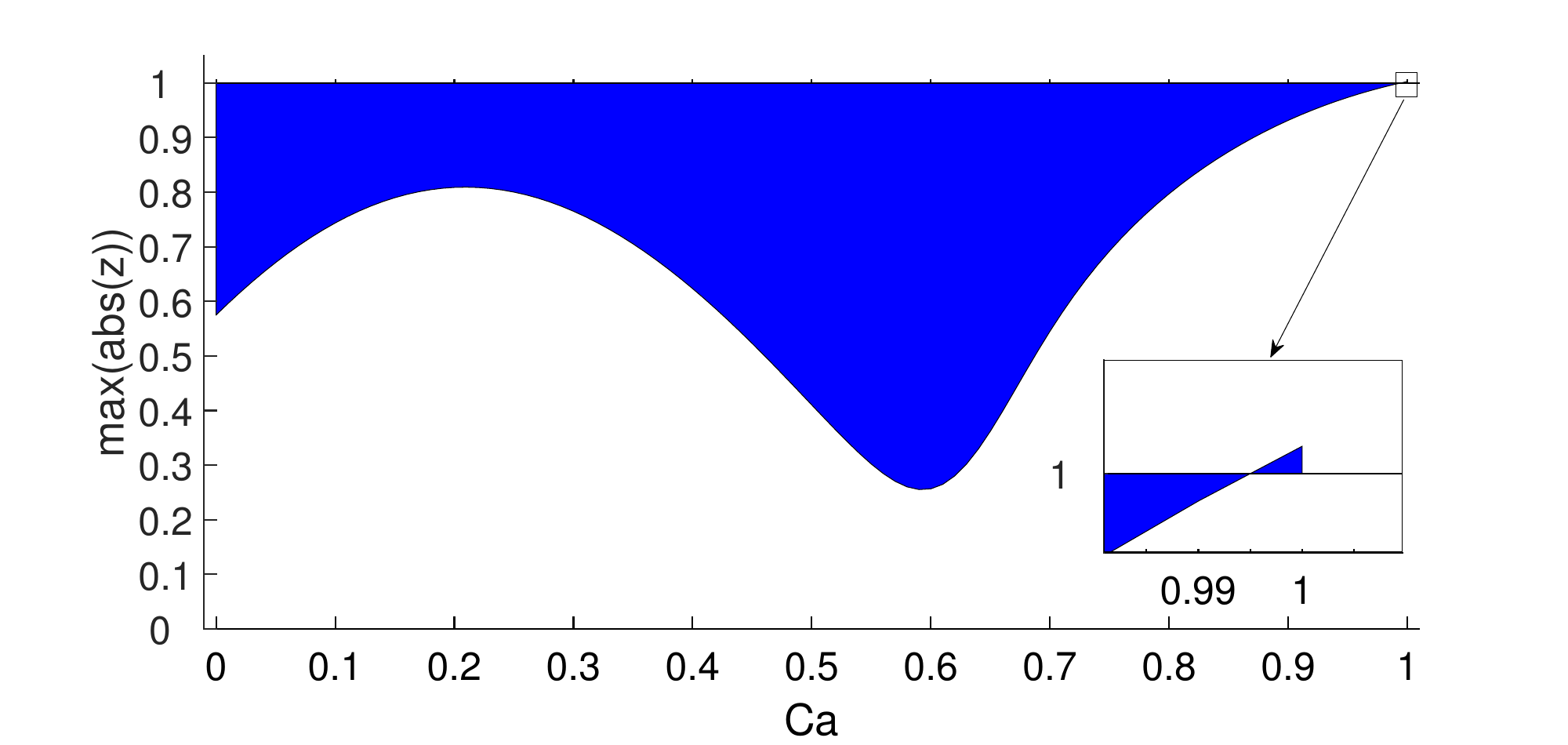}}
\caption{The result of linear stability analysis with $k_d=4$. The horizontal axis represents $C_a$ and the vertical axis represents the largest $|z(\mu)|$.
}
\label{fig:stabilitykd4}
\end{figure}

The linear stability analysis results of the new SILW method with different internal schemes and different $k_d$ are shown in Figures \ref{fig:stabilitykd2}-\ref{fig:stabilitykd4}. The parameter $\alpha$ is taken as the critical value between stable and unstable thresholds. These figures verify the correctness and the optimality of the $\alpha$ range given in Table \ref{TAB:NEWSILW_STABILITY}. 

\end{appendices}
%

\end{document}